\documentclass[preprint,12pt,3p]{elsarticle}
\usepackage{amsfonts} 
\usepackage{amsmath}
\usepackage{amssymb}
\usepackage{amsthm}
\usepackage{graphicx} 
\usepackage[caption=false]{subfig}
\usepackage{xcolor} 
\usepackage{float} 
\usepackage[utf8]{inputenc}
\graphicspath{ {figures/} }
\usepackage{array}
\usepackage{empheq}
\usepackage{hyperref}
\usepackage{mathtools}
\DeclarePairedDelimiter{\ceil}{\lceil}{\rceil}

\usepackage{epstopdf}

\journal{Chaos, Solitons and Fractals}

\begin{document}

\begin{frontmatter}

\title{Distinguished correlation properties of Chebyshev dynamical systems and their generalisations}

\author[label]{Jin Yan}
\ead{j.yan@qmul.ac.uk} 
\author[label]{Christian Beck}
\ead{c.beck@qmul.ac.uk}
\address[label]{School of Mathematical Sciences, Queen Mary University of London, Mile End Road, London E1 4NS, UK}

\begin{abstract}
                                   We show that, among all smooth one-dimensional maps conjugated to an $N$-ary shift (a Bernoulli shift of $N$ symbols), Chebyshev maps are distinguished in the sense that they have least higher-order correlations. We generalise our consideration and study a family of shifted Chebyshev maps, presenting analytic results for two-point and higher-order correlation functions. We also review results for the eigenvalues and eigenfunctions of the Perron-Frobenius operator of $N$-th order Chebyshev maps and their shifted generalisations. The spectrum is degenerate for odd $N$. Finally, we consider coupled map lattices (CMLs) of shifted Chebyshev maps and numerically investigate zeros of the temporal and spatial nearest-neighbour correlations, which are of interest in chaotically quantized field theories.
\end{abstract} 


\end{frontmatter}

\section{Introduction}
Maps conjugated to a Bernoulli shift are a standard paradigm for modelling chaotic dynamics \cite{BoyGora, AHCB2001, CB1991, CBblue, CBred, froyland2001, VepstasBern}. It is useful to introduce a generating partition and symbolic dynamics for such maps \cite{CBred, symbolic}. In the symbol space, dynamics generated by these mappings corresponds to a shift of symbols, and from a statistical point of view the sequence of symbols is statistically independent, thus implying strong mixing properties for the map under consideration. 

One may, however, ask further questions for maps conjugated to a Bernoulli shift, as some properties depend on the particular way how the map is conjugated to the shift, i.e. they depend on the function underlying the topological conjugation. An interesting question is about the structure of higher-order correlation functions of the iterates of the map \cite{AHCB2001, CB1991}. Given a map $T: X \to X$ with iterates $x_{n+1} = T(x_n)$ these higher-order correlation functions are defined as $\langle x_{n_1} x_{n_2} \cdots x_{n_r}\rangle$, where $\langle \cdots \rangle$ denotes the expectation with respect to the (natural) invariant measure of the map. We assume that the average $\langle x_n \rangle$ is zero (if it is non-zero it can be just subtracted from the iterates). An interesting question is: which smooth map conjugated to a Bernoulli shift of $N$ symbols with average $\langle x_n \rangle =0$ is the ``most random" one, in the sense that it has the largest number of tuples $(n_1, \ldots , n_r)$ such that the higher-order correlation function is exactly zero?  

In this paper we show that the answer is given by Chebyshev maps of $N$-th order. These are conjugated to the shift of $N$ symbols by means of a cosine function, and have been subject of many previous papers \cite{adler, geisel, dettmann, groote, williams}. We will show that any other conjugating function produces more higher-order correlations. So for example, the binary shift map, $T(x) =2x \mod 1$ (with subtracted mean), has more non-vanishing higher-order correlations than the second-order Chebyshev map $T(x)=2x^2-1$, which is topologically conjugated via a cosine function. We show that all higher-order correlations can be analytically understood by studying the solutions of a certain set of diophantine equations, which can be solved by a graph-theoretical method, for general $N$.

We will generalize our considerations to shifted Chebyshev maps, defined as $T(x) = \cos (N \arccos x + a)$, which are conjugated by the same cosine function. We will evaluate higher-order correlations in full generality for any $a \in [-\pi/2, 0]$, and show that again $a = 0$ yields the smallest possible skeleton of non-zero higher-order correlations. We will find a suitable topological conjugation for shifted Chebyshev maps, for particular values of $a$ and $N$, which relates their dynamics to that of ordinary ($a = 0$) Chebyshev maps.

Investigating the higher-order correlation structure is significant especially when generating ``noise” in a stochastic differential equation by a smooth deterministic chaotic dynamics \cite{williams, CBGR, melbourne, mackey}. Clearly, the mathematical construction of Gaussian white noise that drives a stochastic differential equation possesses no correlations at all in time, but one may ask what type of higher-order correlations are generated by a smooth deterministic chaotic system at a microscopic level. Since the dynamics is deterministic and discrete there must be correlations even for maps conjugated to a Bernoulli shift: only the symbols that are shifted are statistically independent, but the iterates itself are not. We show that Chebyshev systems have least non-vanishing higher-order correlations when calculated with respect to their invariant measure, and are in this sense as close to white noise as possible for a smooth one-dimensional chaotic dynamics. 

It is important to emphasise smoothness here. For example, a random number generator is not a smooth function of its seed variable. 
The question of the ultimate source of the noise in stochastically quantized field theories was discussed in \cite{CBblue}: one may assume that there is always a deterministic dynamics at the smallest scales (say, the Planck scale), since by definition the smallest scale cannot contain additional degrees of freedom that are just effectively described by a random process. For this reason, in \cite{CBblue, groote, chaotic-string, nonlin} chaotically quantized field theories (sometimes also called ``chaotic strings") were studied, which do possess a chaotic dynamics generating the ``noise" of the path integral approach in a deterministic way on the smallest scale. States of zero spatial nearest-neighbour correlations were identified as physical states in this approach. Our consideration here shows that Chebyshev maps are the most distinguished candidates for such a fundamental noise dynamics at the Planck scale, with a minimum possible skeleton of higher-order correlations.

This paper is organised as follows. In section 2 we introduce shifted Chebyshev maps, which contain the ordinary Chebyshev maps as a special case ($a = 0$). We discuss the general behaviour and the invariant measures of these maps and provide some examples. In section 3 we review results on the spectrum of the Perron-Frobenius operator for Chebyshev maps, and present some new results
on a complete set of eigenfunctions for $a$ different from 0. In section 4 we derive the diophantine equations describing the higher-order correlation structure of Chebyshev maps, both for $a = 0$ and for general $a$. We show that for any other conjugating function than the cosine there are more solutions to these equations, thus leading to more non-vanishing higher-order correlations. In section 5 we consider coupled map lattices (CMLs) of shifted Chebyshev maps, which are of relevance in chaotically quantized 
field theories. We present numerical results on spatial and temporal two-point correlation functions. The shape depends both on the coupling parameter $c$ as well as on the shift parameter $a$. We discuss possible physical applications for these types of CMLs
in terms of generating the ``noise" in chaotically quantized field theories. Finally, we present our conclusions in section 6.

\section{The shifted Chebyshev maps}
\subsection{Definition}
\begin{figure}[H]
\centering 
\subfloat[\footnotesize $N = 2, a = 0$]{
\includegraphics[width = 0.24\textwidth]{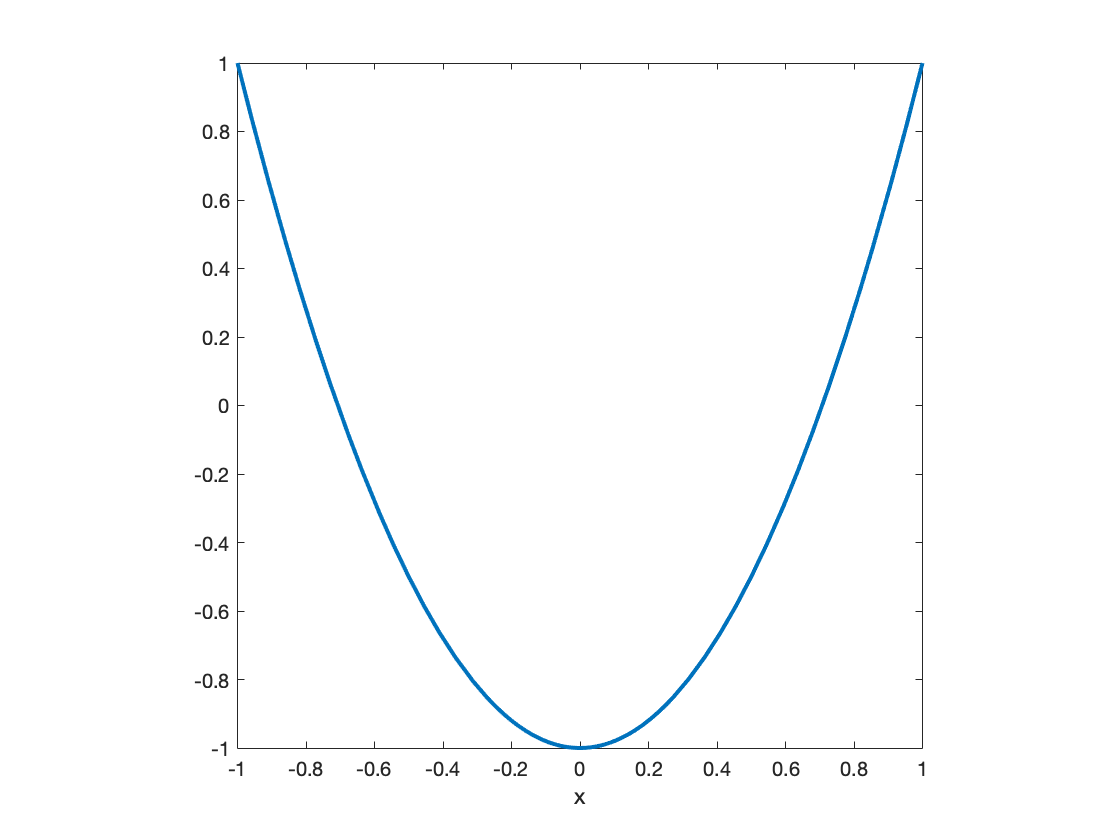}
}
\subfloat[\footnotesize $N = 3, a = 0$]{
\includegraphics[width = 0.24\textwidth]{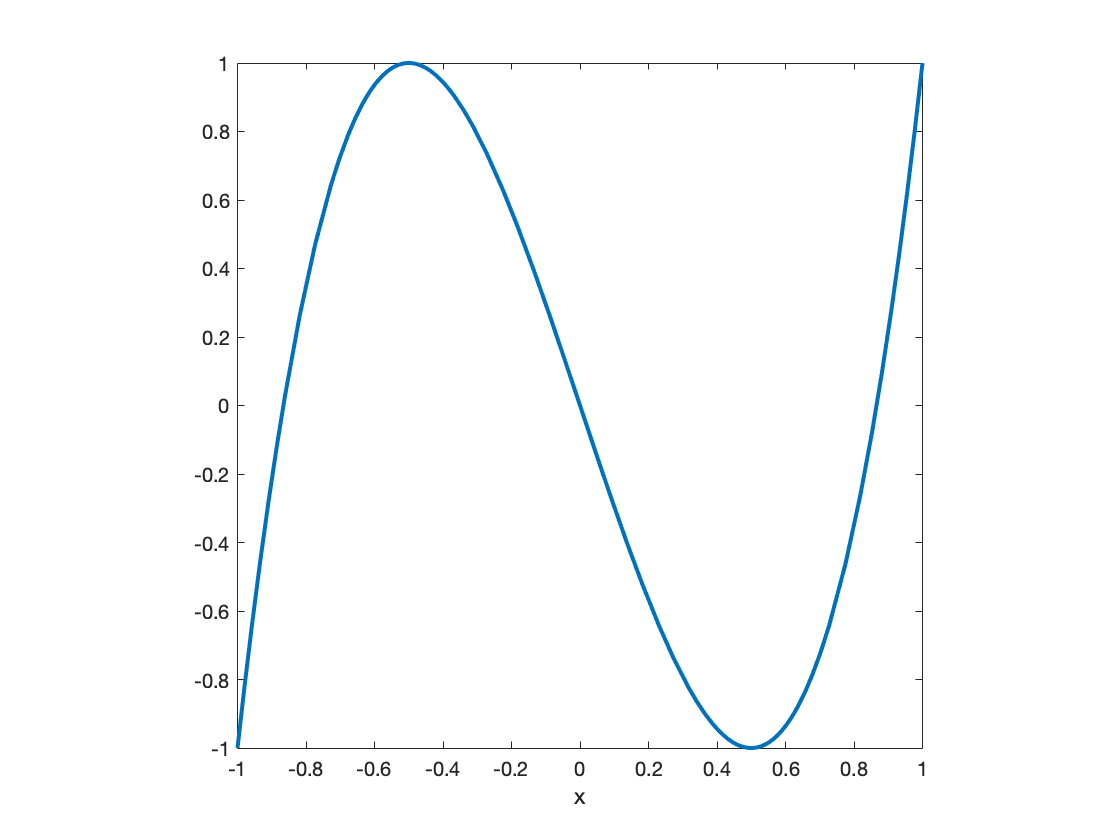}
}
\subfloat[\footnotesize $N = 4, a = 0$]{
\includegraphics[width = 0.24\textwidth]{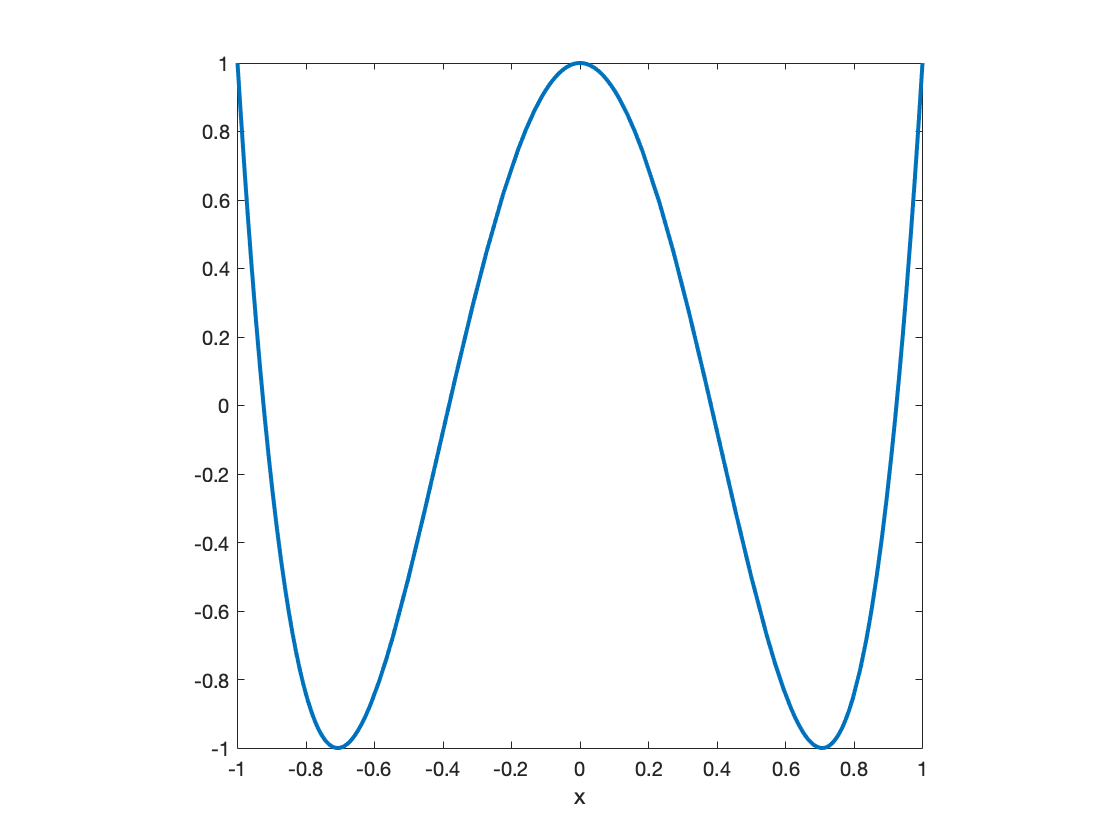}
}
\subfloat[\footnotesize $N = 7, a = 0$]{
\includegraphics[width = 0.24\textwidth]{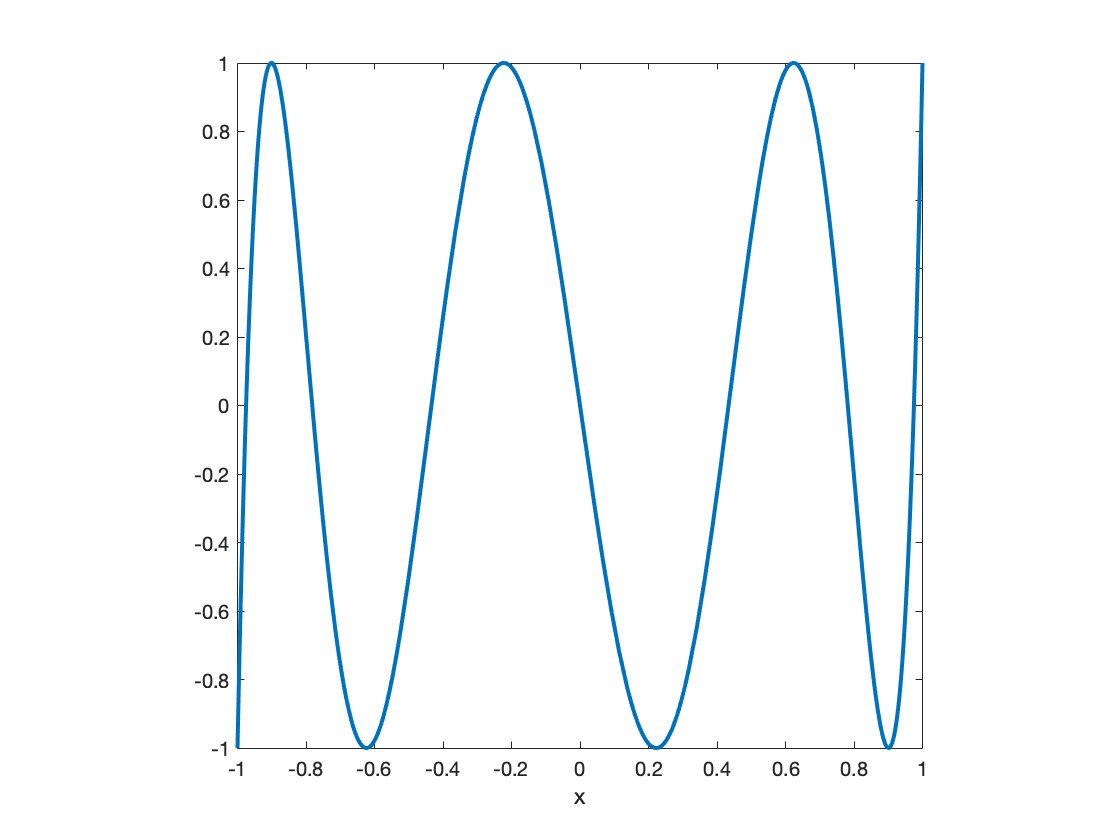}
}
\quad

\subfloat[\footnotesize $N = 2, a = -\pi/4$]{
\includegraphics[width = 0.24\textwidth]{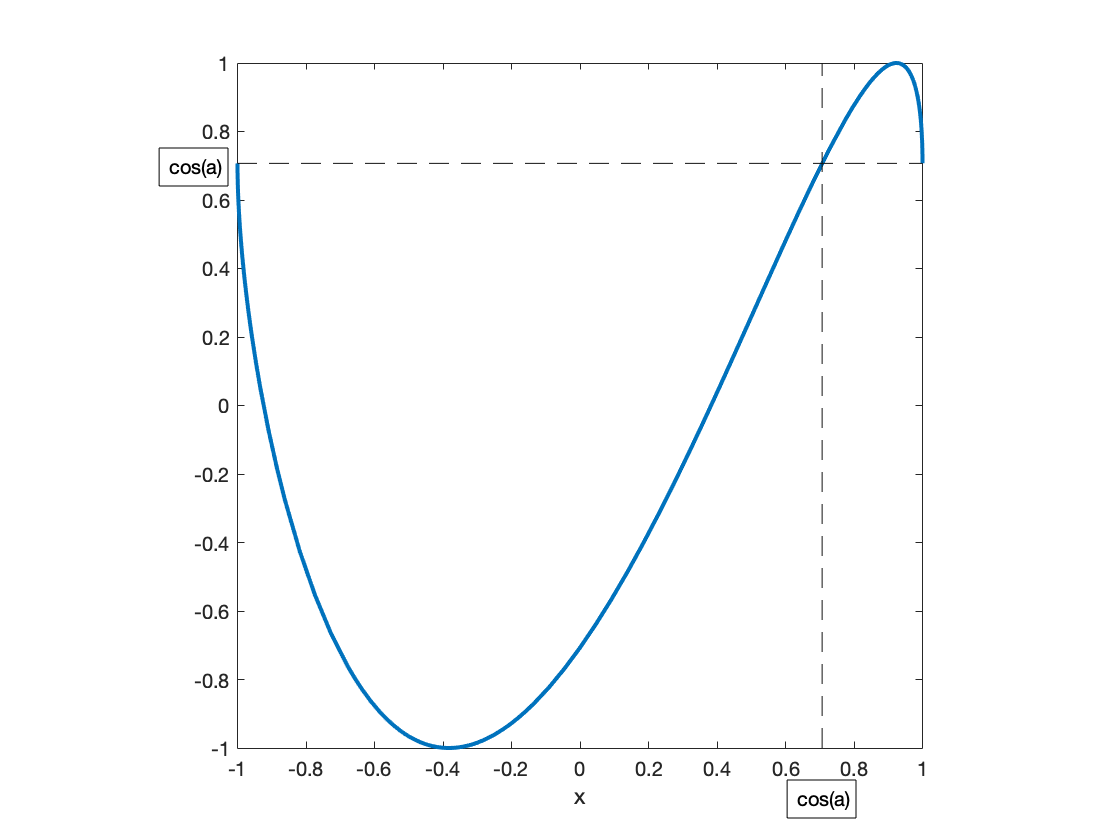}
\label{N2-a-pi4(1)}}
\subfloat[\footnotesize $N = 3, a = -\pi/4$]{
\includegraphics[width = 0.24\textwidth]{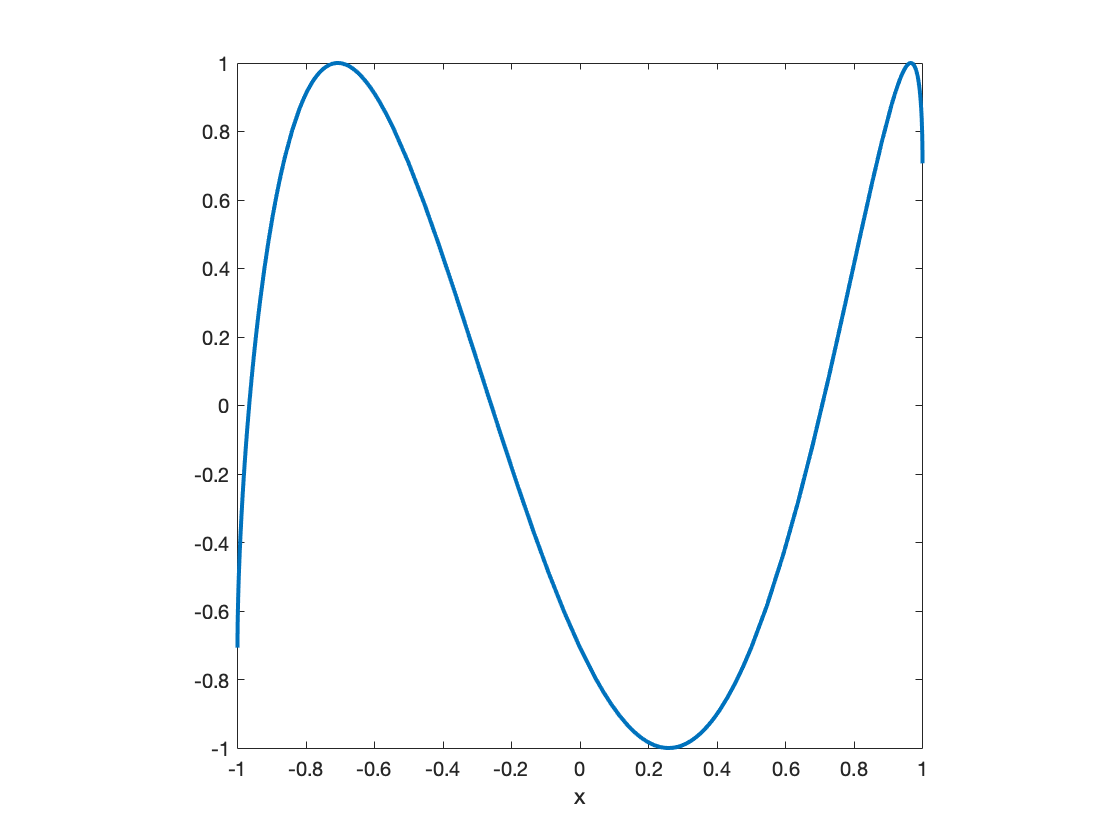}
}
\subfloat[\footnotesize $N = 4, a = -\pi/4$]{
\includegraphics[width = 0.24\textwidth]{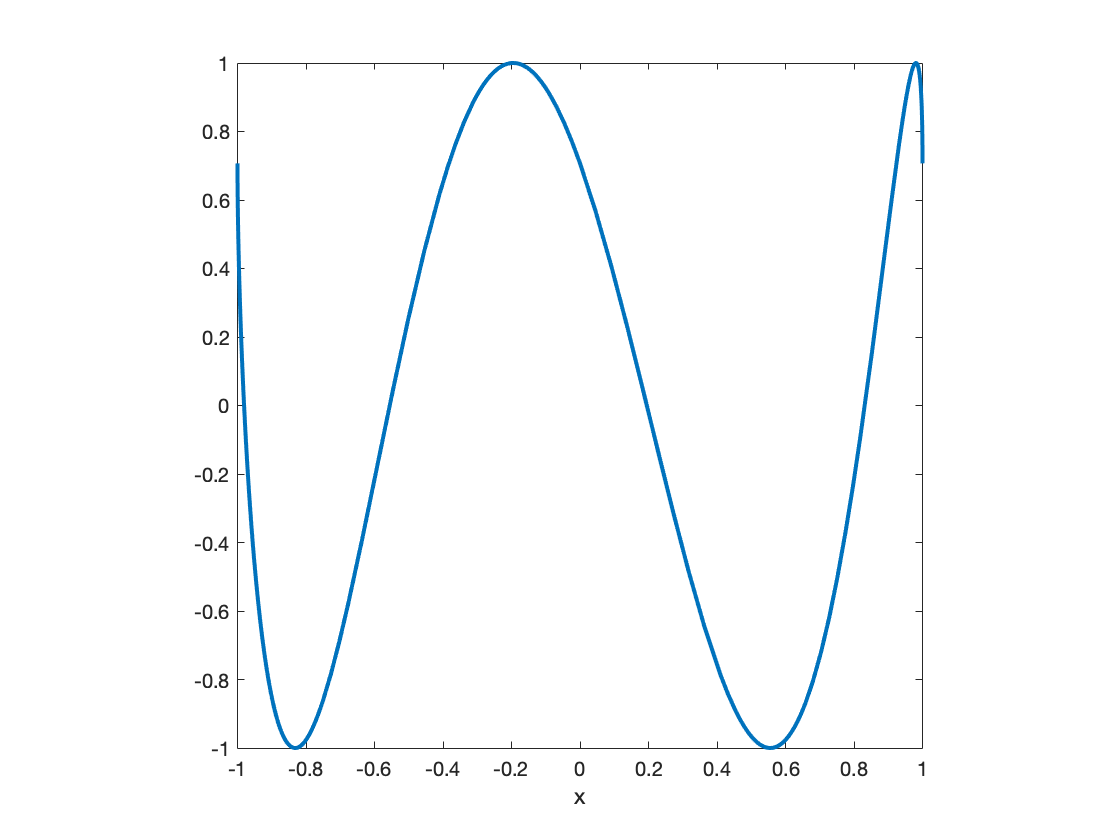}
}
\subfloat[\footnotesize $N = 7, a = -\pi/4$]{
\includegraphics[width = 0.24\textwidth]{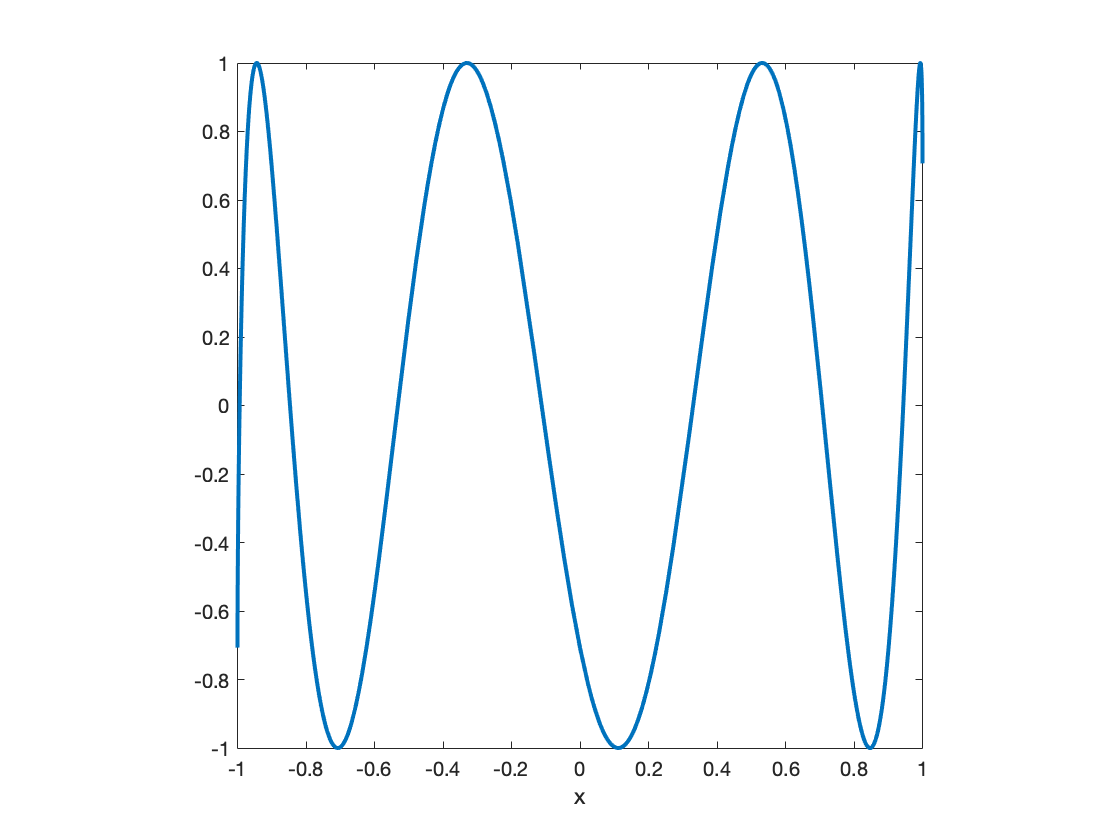}
}
\quad

\subfloat[\footnotesize $N = 2, a = -\pi/2$]{
\includegraphics[width = 0.24\textwidth]{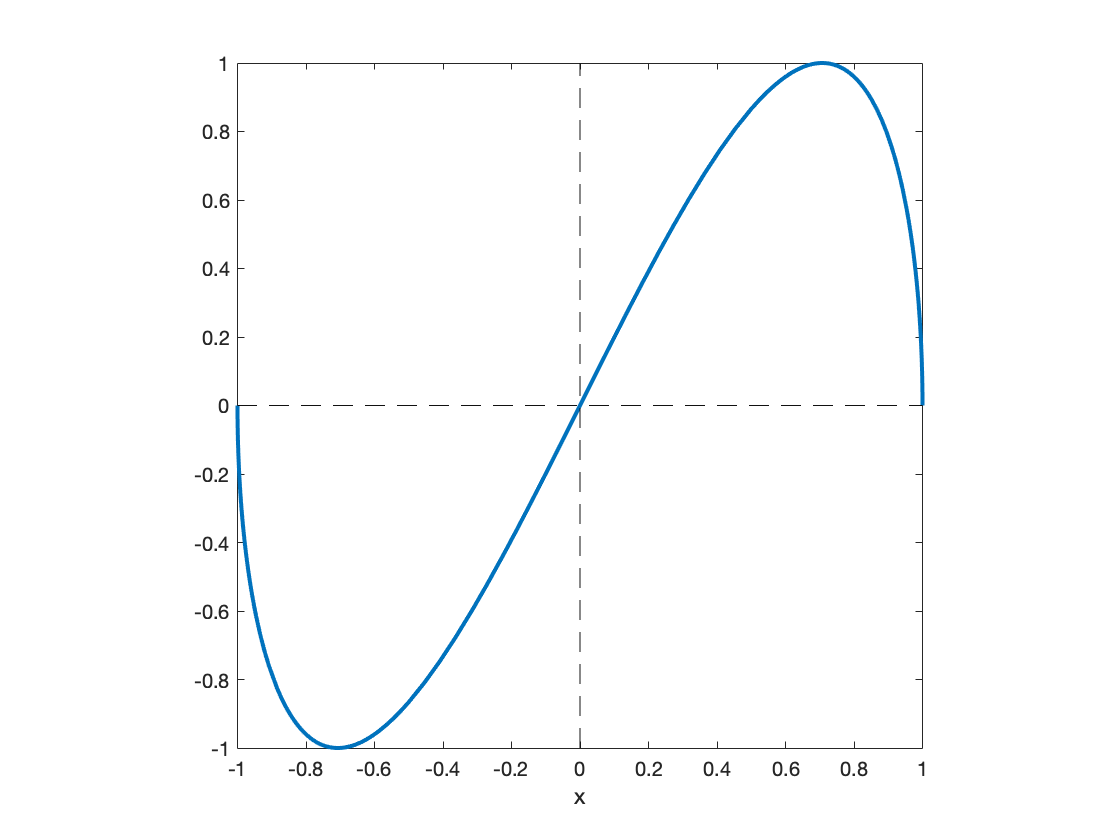}
\label{N2-a-pi2(1)}}
\subfloat[\footnotesize $N = 3, a = -\pi/2$]{
\includegraphics[width = 0.24\textwidth]{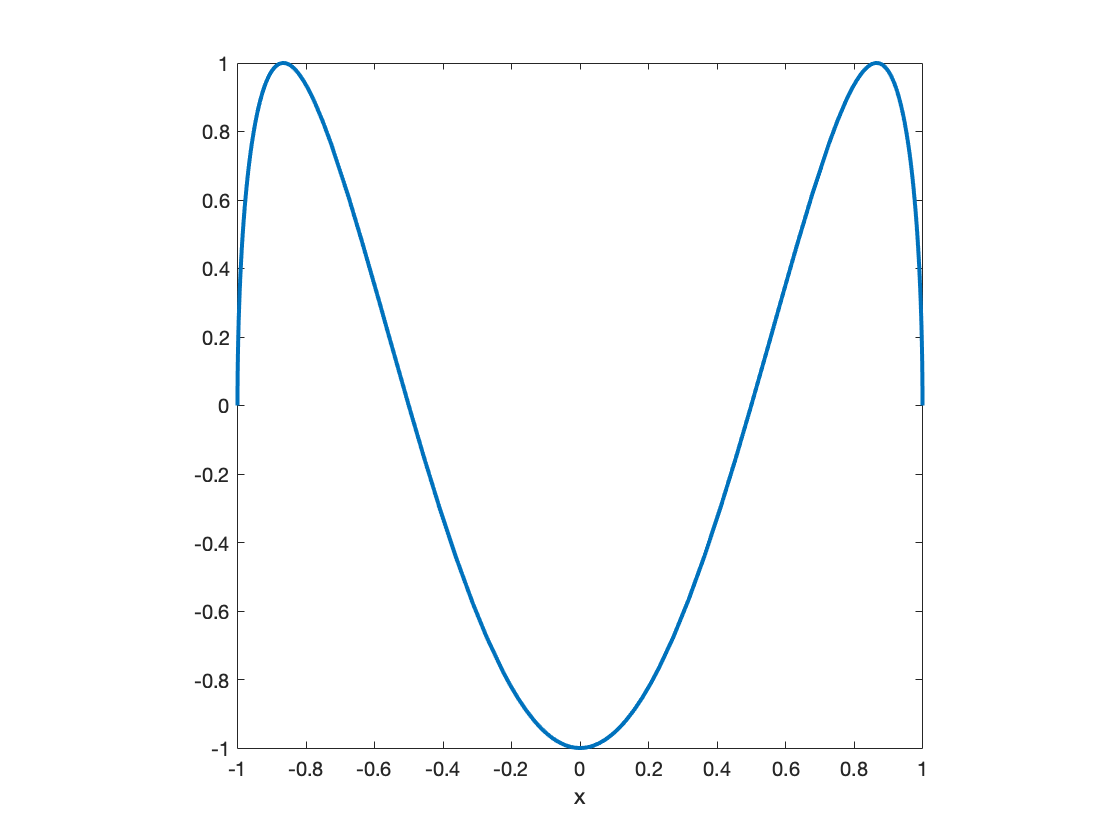}
}
\subfloat[\footnotesize $N = 4, a = -\pi/2$]{
\includegraphics[width = 0.24\textwidth]{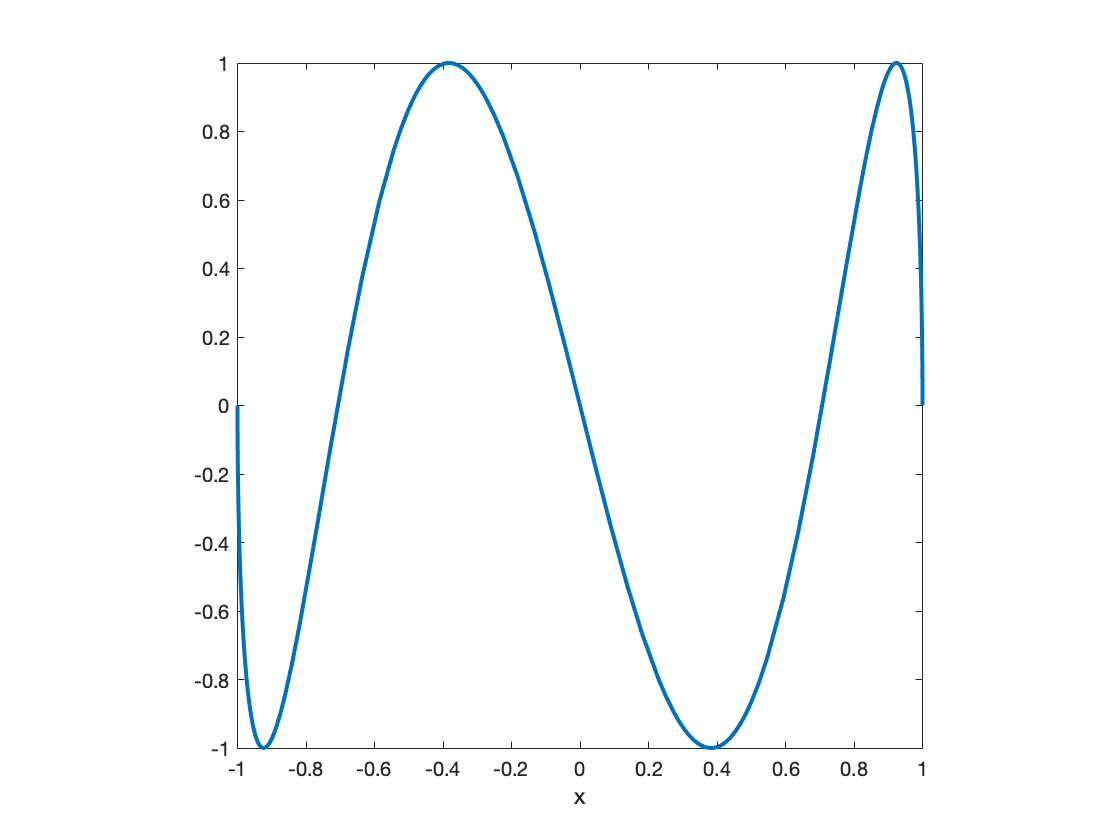}
}
\subfloat[\footnotesize $N = 7, a = -\pi/2$]{
\includegraphics[width = 0.24\textwidth]{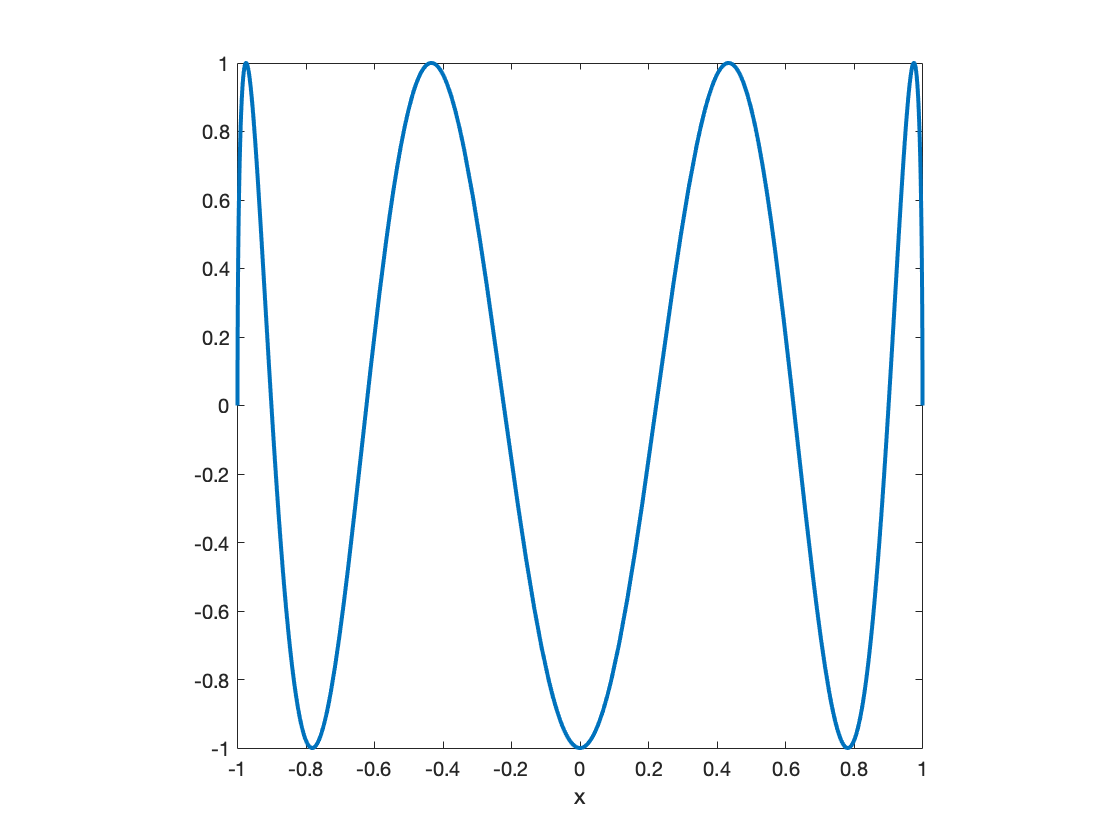}
}
\caption{Graphs of the shifted Chebyshev maps $T_{N, a}(x)$, with $N$ and $a$ indicated in each caption.}
\label{graphs-shifted-cheby}
\end{figure}

Define a discrete-time dynamical system $x_{n+1} = T_{N, a}(x_n)$, $n = 0, 1, ...$, as 
\begin{equation}
T_{N, a}(x) := \cos (N\arccos x + a), \quad x \in [-1, 1],
\label{def-shifted-cheby}%
\end{equation}
with $N = 2, 3, ...$ and $a \in \left[-\frac{\pi}{2}, 0\right]$. $T_{N, a}$ is called the \textit{shifted Chebyshev maps of order $N$}. For $a = 0$ we have \textit{ordinary} Chebyshev maps. Some graphs are shown in Fig.\ref{graphs-shifted-cheby}. 

Notice that for $N = 2$ and $a \neq 0$, the map can be decomposed into two independent (ergodic) components, $I_1 = [-1, \cos a]$ and $I_2 = [\cos a, 1]$, which does not happen for $N > 2$, cf. Figs.\ref{N2-a-pi4(1)}, \ref{N2-a-pi2(1)}. 

\subsection{Topological conjugation}
Consider a change of variables, $x_0 = \cos(\pi u_0) =: \cos \left( u_0' - \frac{a}{N-1}\right)$, $u_0 \in (0, 1)$, then
\begin{equation*}
\begin{split}
x_n &= \cos\left(N^n\pi u_0 + N^{n-1}a + N^{n-2}a + ... + Na + a\right)\\
&= \cos\left(N^n\pi u_0 + \frac{N^n - 1}{N - 1}a\right)\\
&= \cos\left( N^n \left(\pi u_0 + \frac{a}{N-1}\right) - \frac{a}{N-1}\right)\\
&= \cos\left( N^n u_0' - \frac{a}{N-1}\right).
\end{split}
\end{equation*}
So essentially the dynamics can be reduced to a (transformed) $N$-ary shift, at each time step shifting the variable $u_0' := \pi u_0 + \frac{a}{N-1} = \arccos (x_0) + \frac{a}{N-1}$ by one digit in its $N$-ary representation: 
$u_n' = N^n u_0'$, $u_0' \in \left( \frac{a}{N-1}, \frac{a}{N-1} + \pi \right)$. 

\begin{figure}[H]
\centering
\subfloat[\footnotesize $T_{3, a}$]{
\includegraphics[width = 0.48\textwidth]{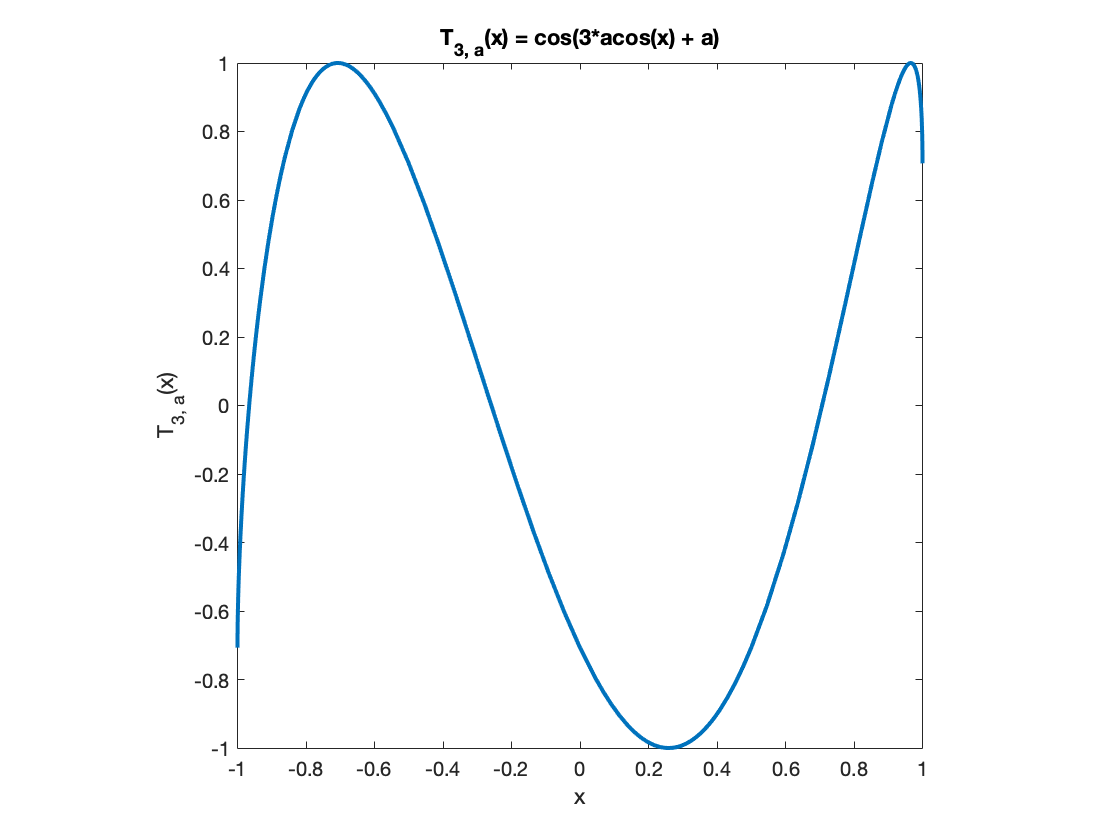}
}
\subfloat[\footnotesize $g_{3, a}$]{
\includegraphics[width = 0.48\textwidth]{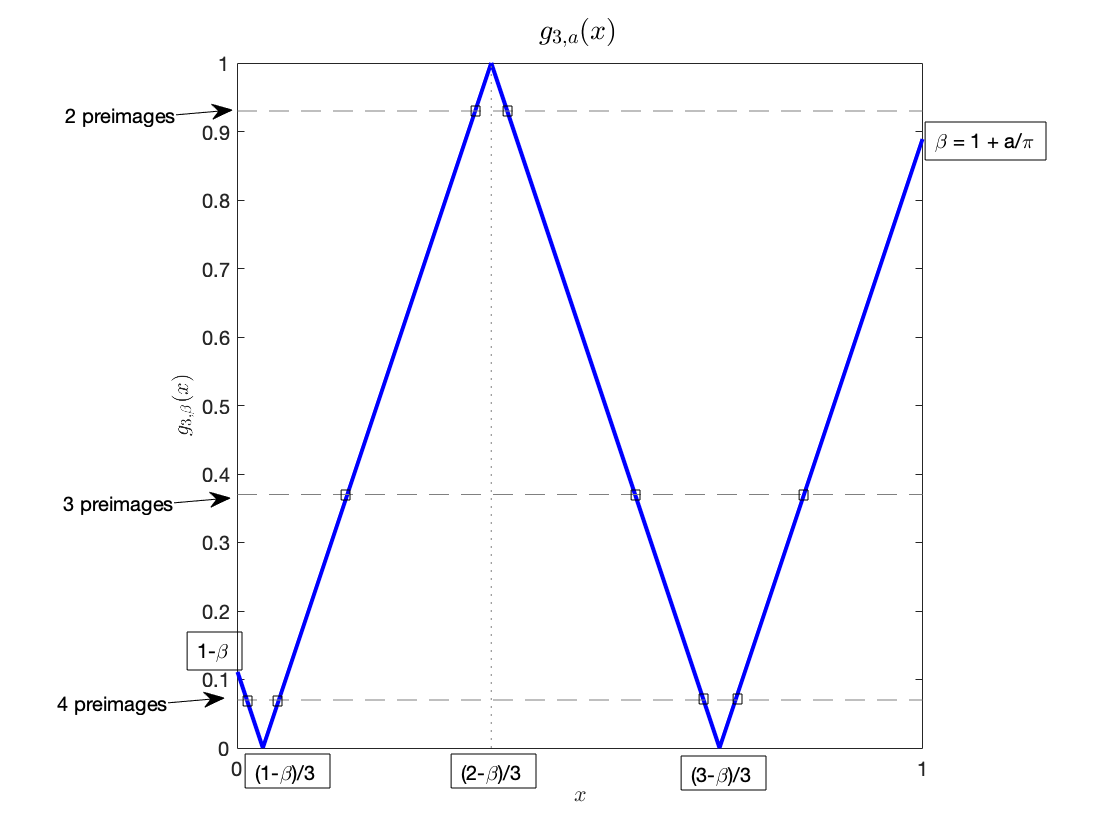} 
\label{counter-eg-of-full-branch}}
\caption{$N = 3$: shifted Chebyshev map $T_{N, a}$ and its conjugated piecewise-linear map $g_{N, a}$.}
\label{illu-shifted-N3-conj}
\end{figure}

In fact, it can be shown (see \ref{A1}) that for all $N \in \mathbb{N}_{\geq 2}$, $T_{N, a}$ is topologically conjugated to a piecewise-linear map $g_{N, a}$ via the conjugacy $h: [-1, 1] \rightarrow [0, 1]$
\begin{equation*}
h(x) = \frac{1}{\pi}\text{arccos}(-x), \quad x \in [-1, 1] 
\end{equation*}
such that $h\circ T_{N, a} = g_{N, a}\circ h$. The number of branches of $g_{N, a}$ depends on the order $N$, with the corresponding shifted amount being $\beta := 1 + \frac{a}{\pi}$. Fig.\ref{illu-shifted-N3-conj} illustrates such a conjugation for $N = 3$. One has 
\begin{equation*}
g_{3, a}(y) = \begin{cases}
-3y + (1 - \beta), \quad y \in \left[ 0, \frac{1 - \beta}{3}\right)\\
3y - (1 - \beta), \quad y \in \left[ \frac{1 - \beta}{3}, \frac{2 - \beta}{3}\right)\\
-3y + (3 - \beta), \quad y \in \left[ \frac{2 - \beta}{3}, \frac{3 - \beta}{3}\right)\\
3y - (3 - \beta), \quad y \in \left[ \frac{3 - \beta}{3}, 1\right].
\end{cases}
\end{equation*}

\subsection{Invariant densities}
Since the invariant density of a piecewise-linear map is much simpler to find than that of a general shifted Chebyshev map, for a monotonic coordinate transformation we have $\rho_T(x)dx = \rho_g(y)dy$, so that the invariant density of $T_{N, a}$ is given by (we drop the subscripts of $T$ and $g$, for simplicity)
\begin{equation*}
\rho_T(x) = \rho_{g}\frac{dh}{dx} = \frac{1}{\pi\sqrt{1 - x^2}}\rho_g, \quad x \in [-1, 1]. 
\end{equation*}

\begin{figure}[H]
\centering 
\subfloat[\footnotesize $N = 2, a = 0$]{
\includegraphics[width = 0.24\textwidth]{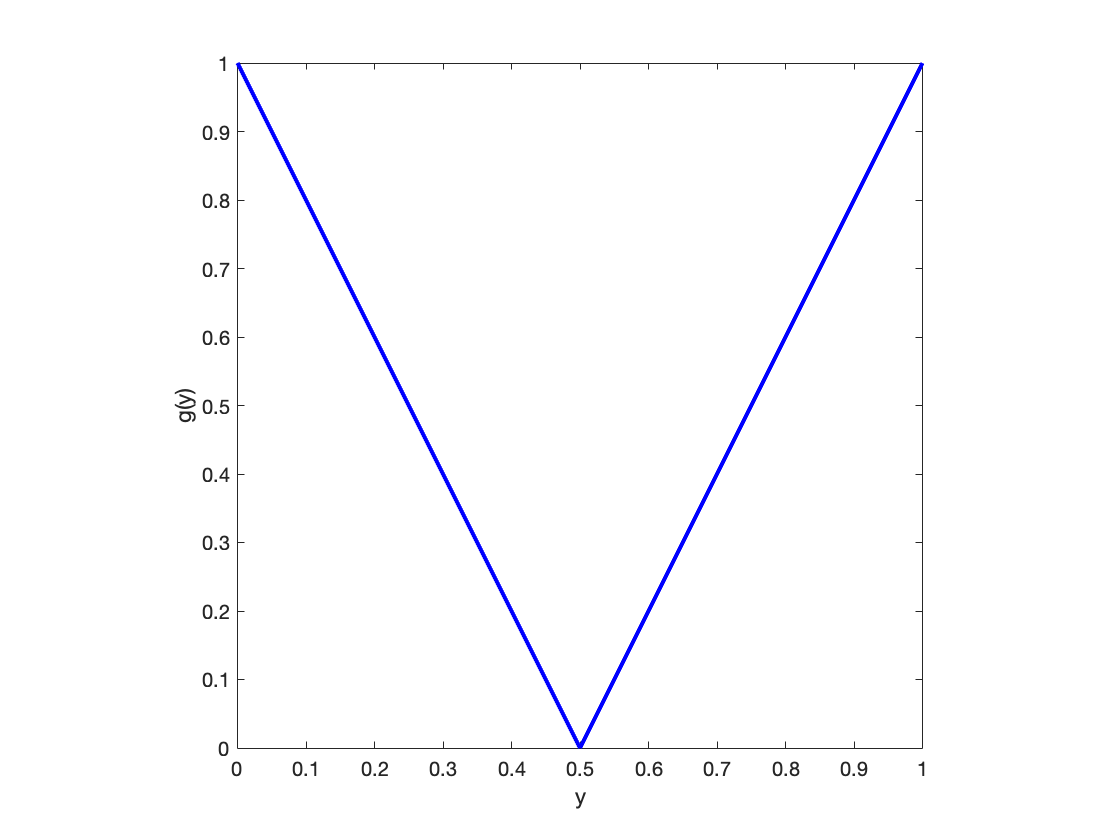}
}
\subfloat[\footnotesize $N = 3, a = -\pi/2$]{
\includegraphics[width = 0.24\textwidth]{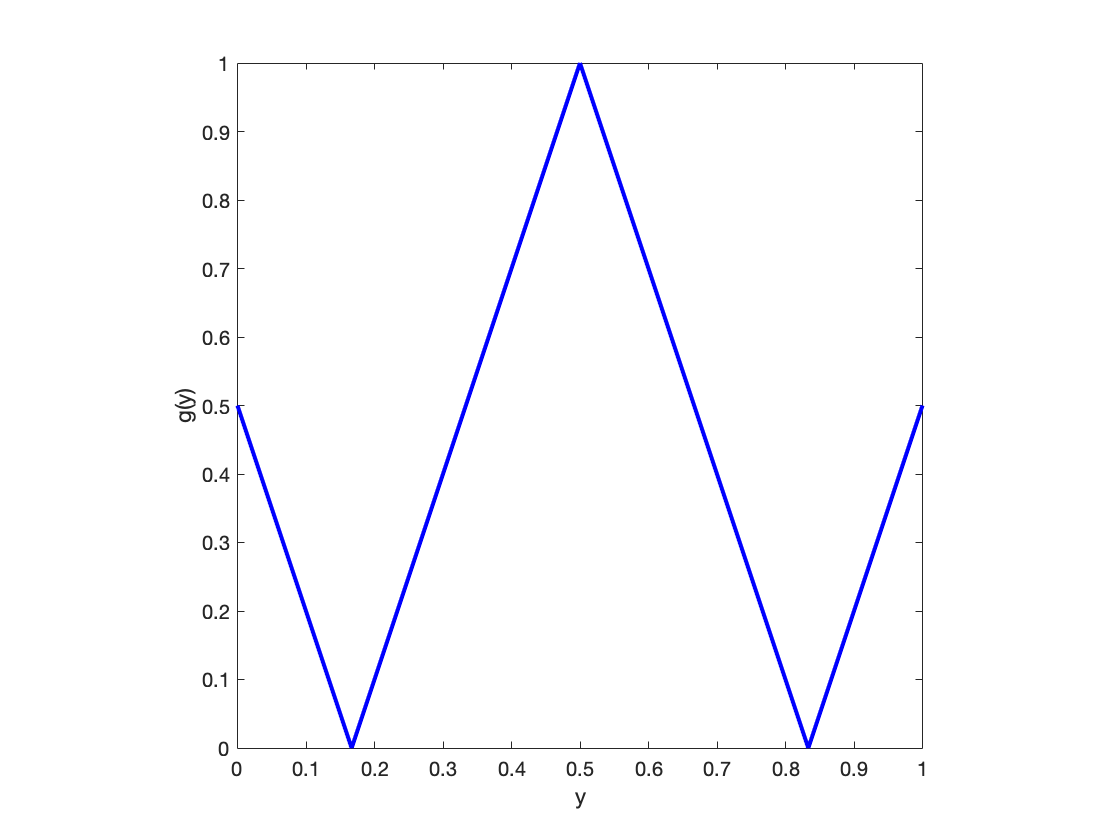}
}
\subfloat[\footnotesize $N = 3, a = -\pi/3$]{
\includegraphics[width = 0.24\textwidth]{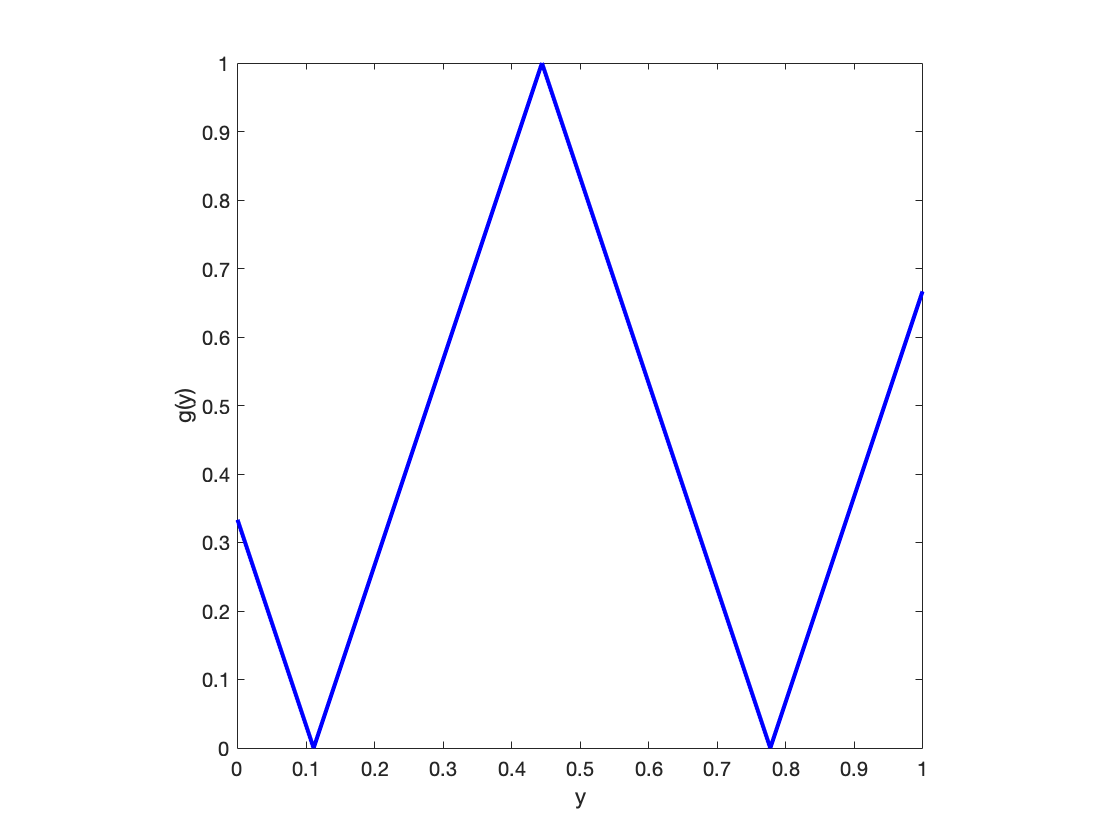}
}
\subfloat[\footnotesize $N = 3, a = -\pi/9$]{
\includegraphics[width = 0.24\textwidth]{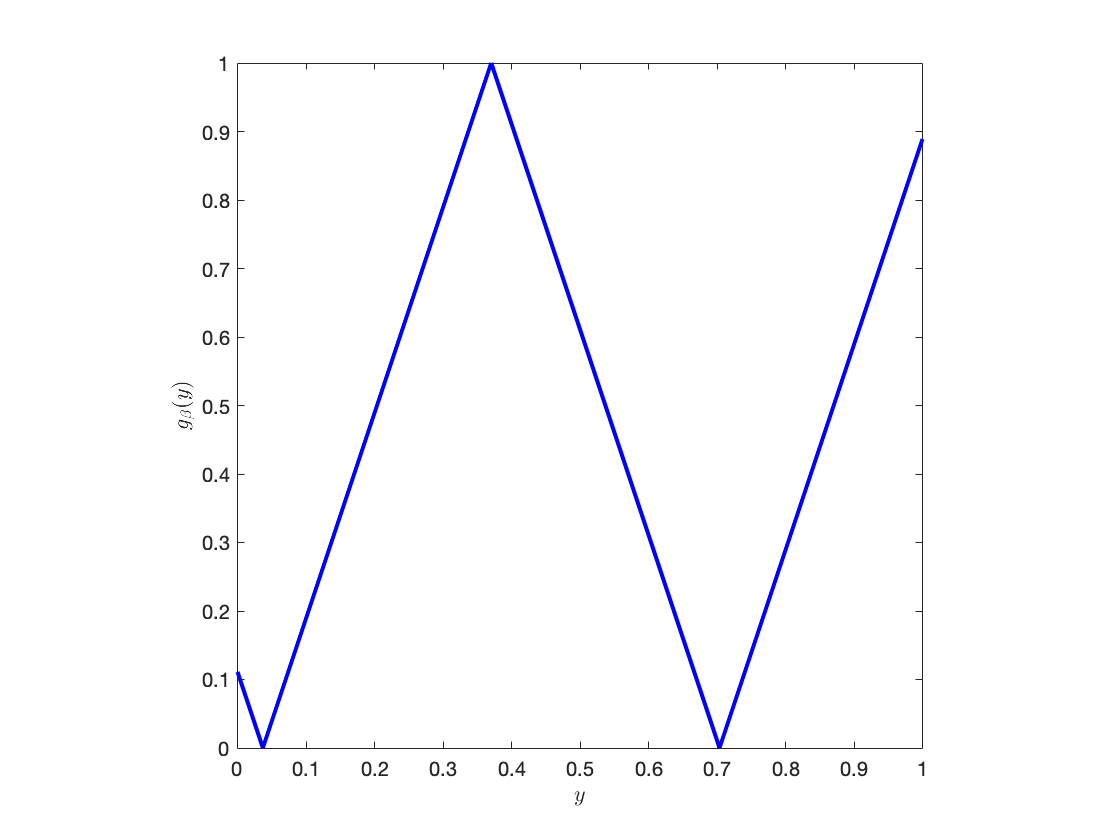}
}
\quad 

\subfloat[\footnotesize $N = 2, a = 0$]{
\includegraphics[width = 0.24\textwidth]{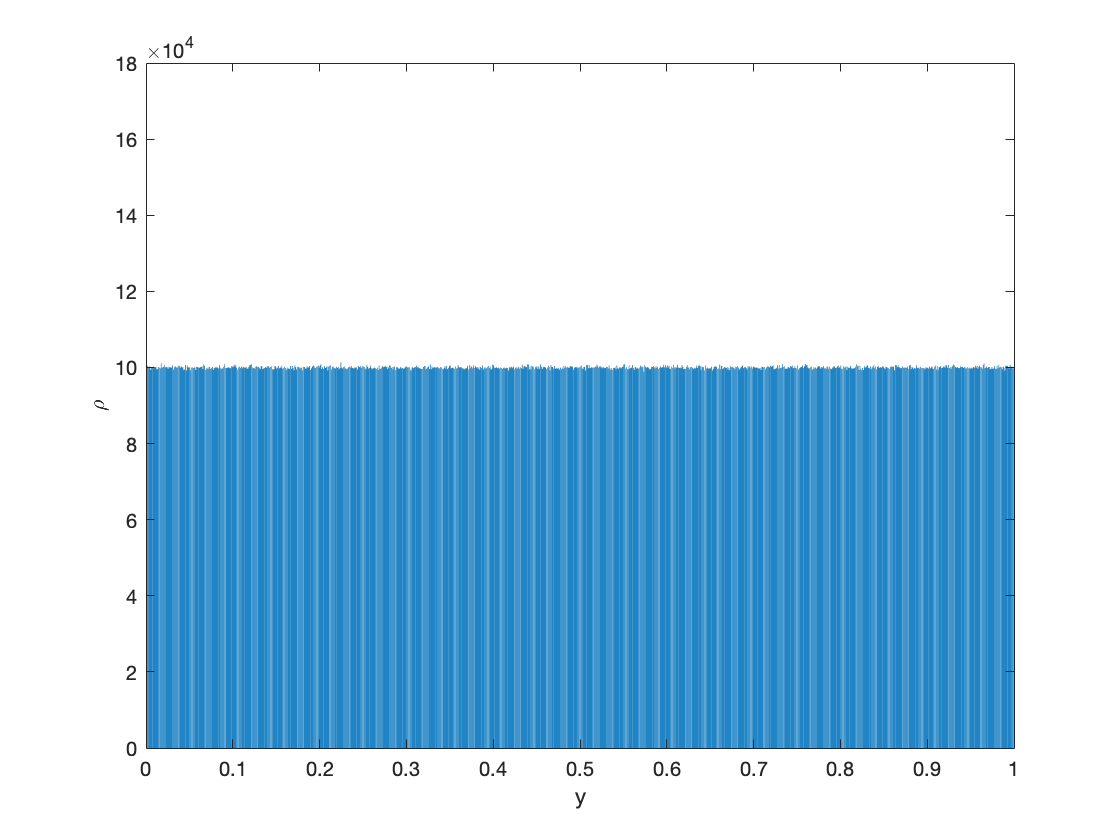}
}
\subfloat[\footnotesize $N = 3, a = -\pi/2$]{
\includegraphics[width = 0.24\textwidth]{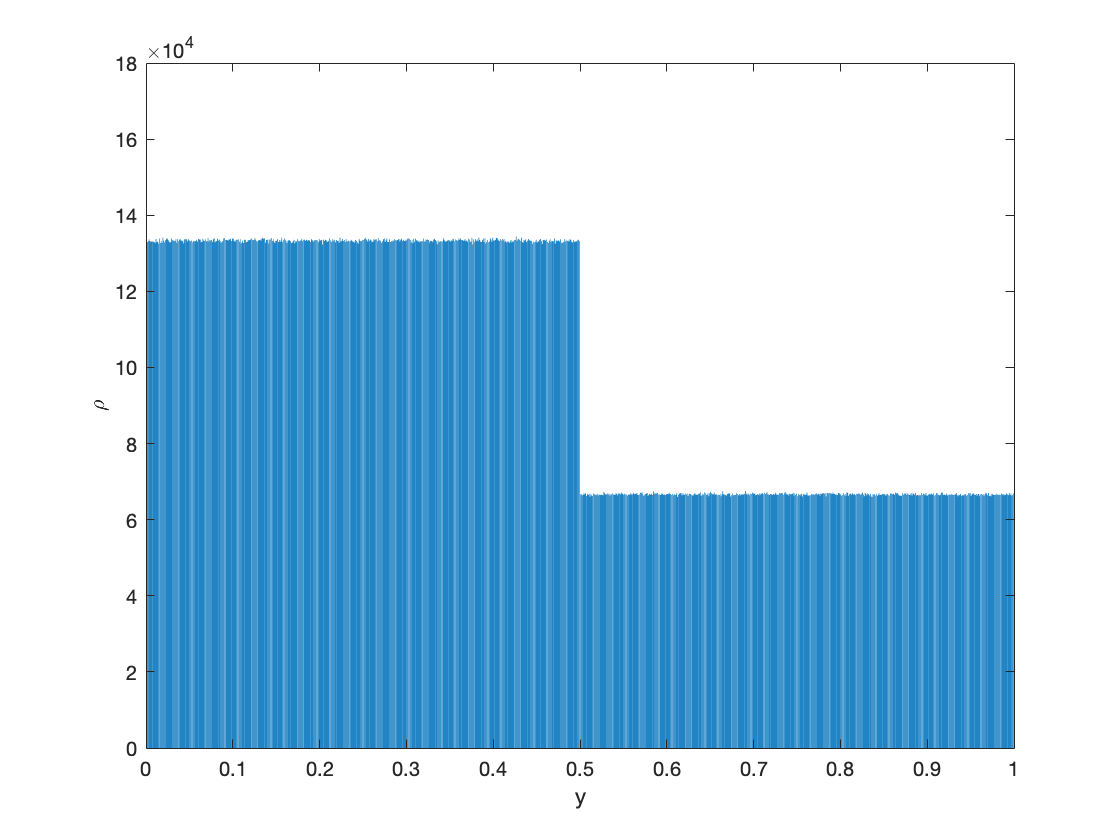}
}
\subfloat[\footnotesize $N = 3, a = -\pi/3$]{
\includegraphics[width = 0.24\textwidth]{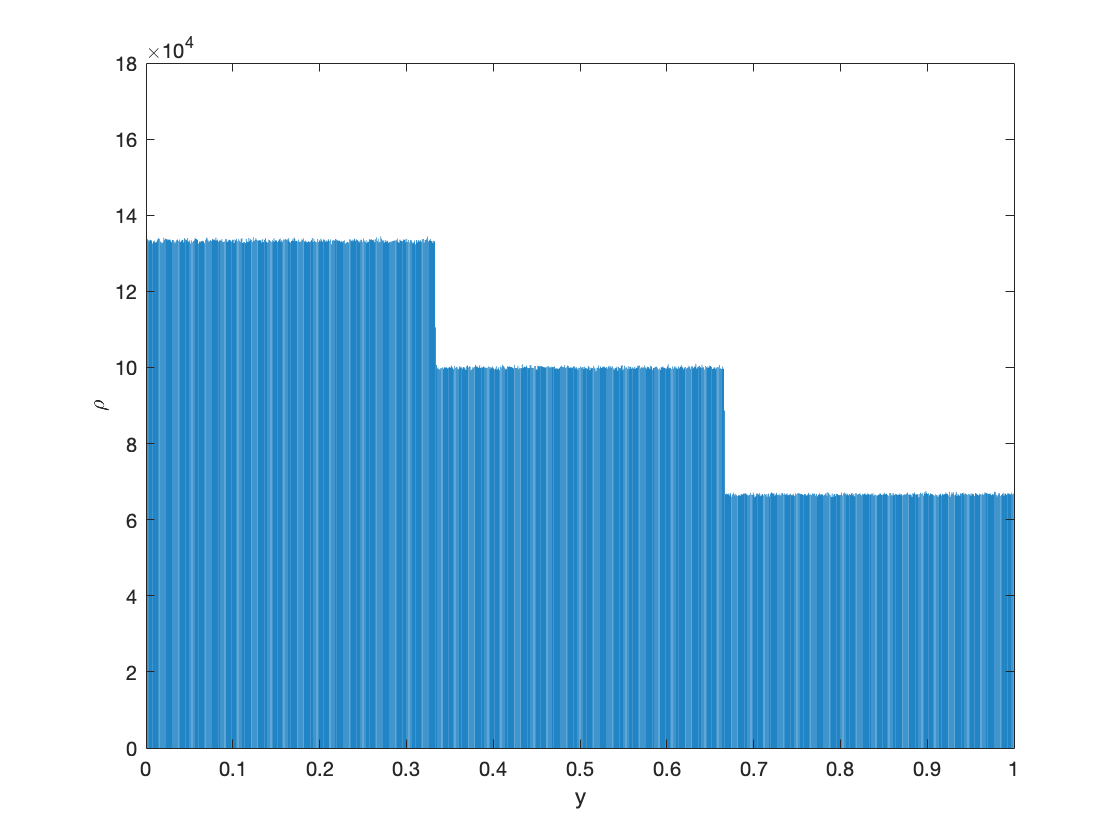}
}
\subfloat[\footnotesize $N = 3, a = -\pi/9$]{
\includegraphics[width = 0.24\textwidth]{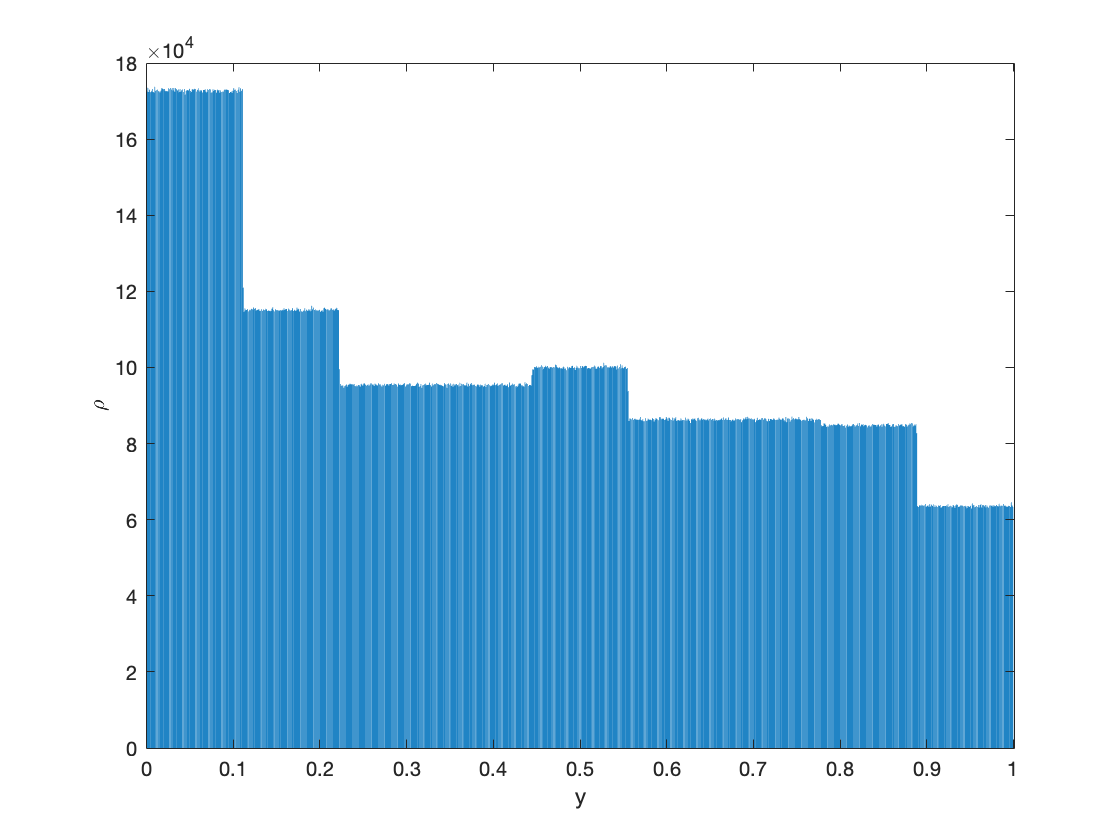}
\label{pwlinear-N3-a-pi9-invden}}
\quad 

\subfloat[\footnotesize $N = 2, a = 0$]{
\includegraphics[width = 0.24\textwidth]{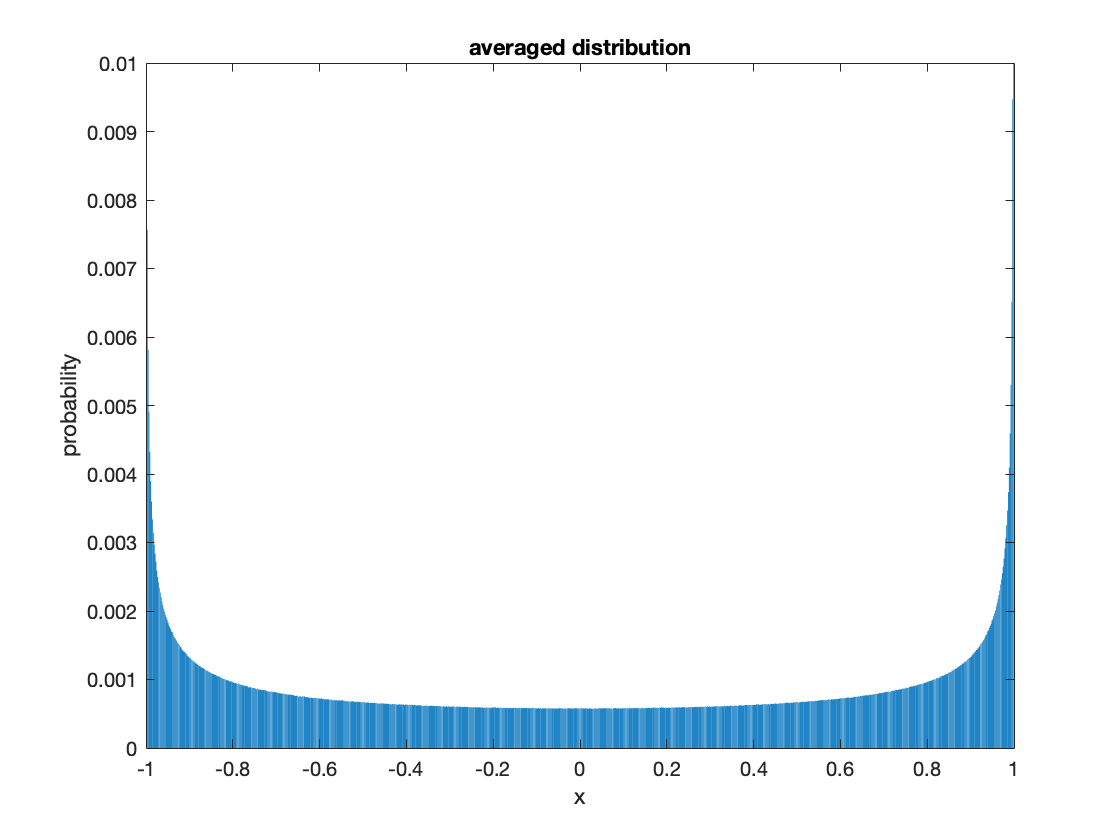}
}
\subfloat[\footnotesize $N = 3, a = -\pi/2$]{
\includegraphics[width = 0.24\textwidth]{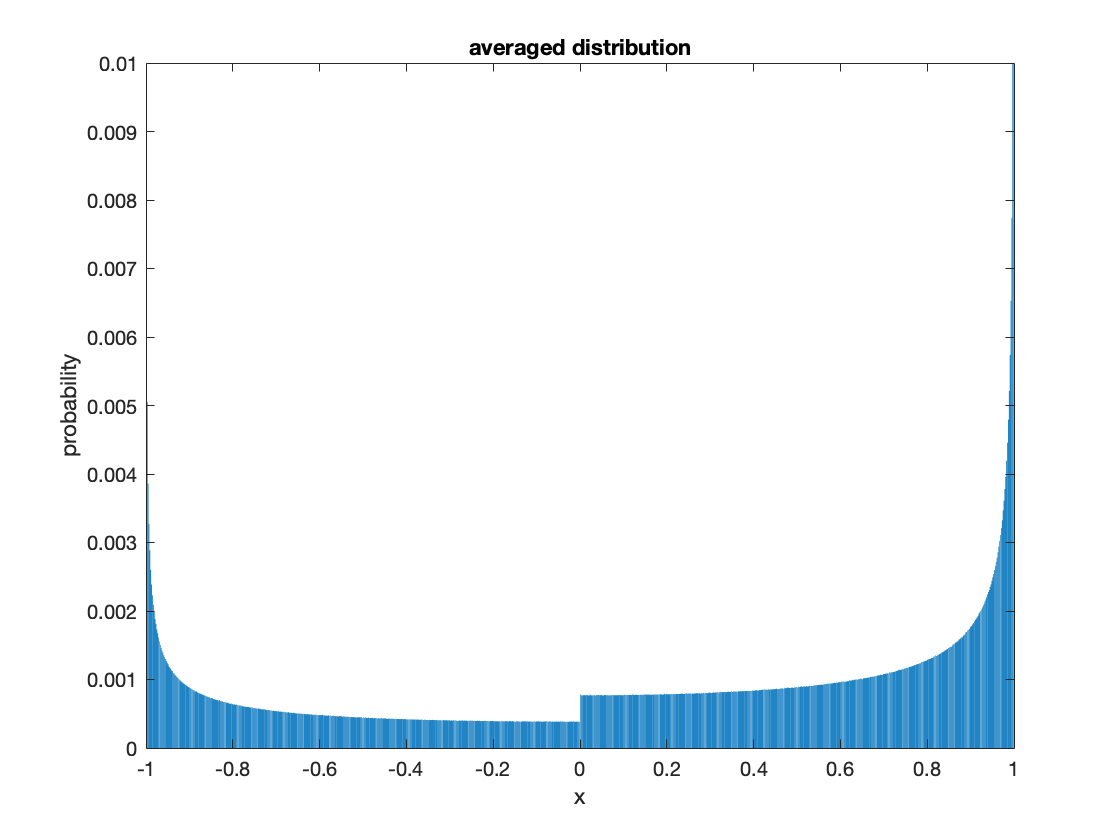}
}
\subfloat[\footnotesize $N = 3, a = -\pi/3$]{
\includegraphics[width = 0.24\textwidth]{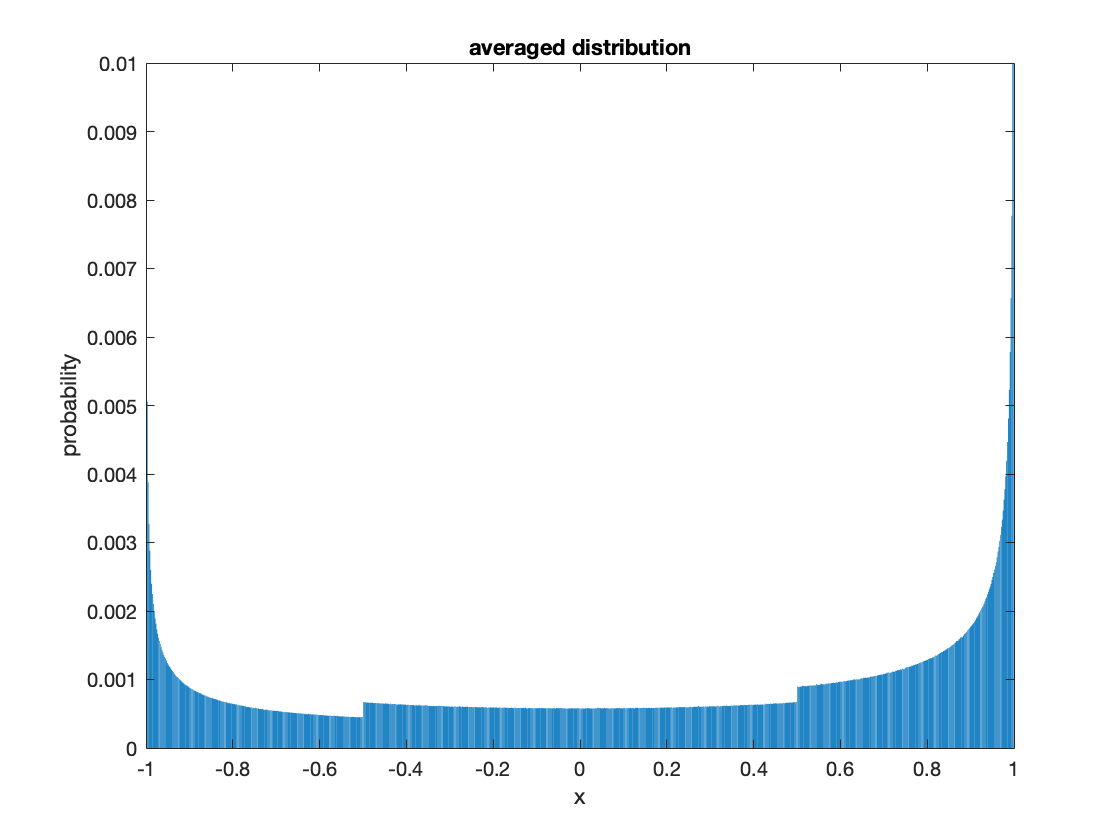}
}
\subfloat[\footnotesize $N = 3, a = -\pi/9$]{
\includegraphics[width = 0.24\textwidth]{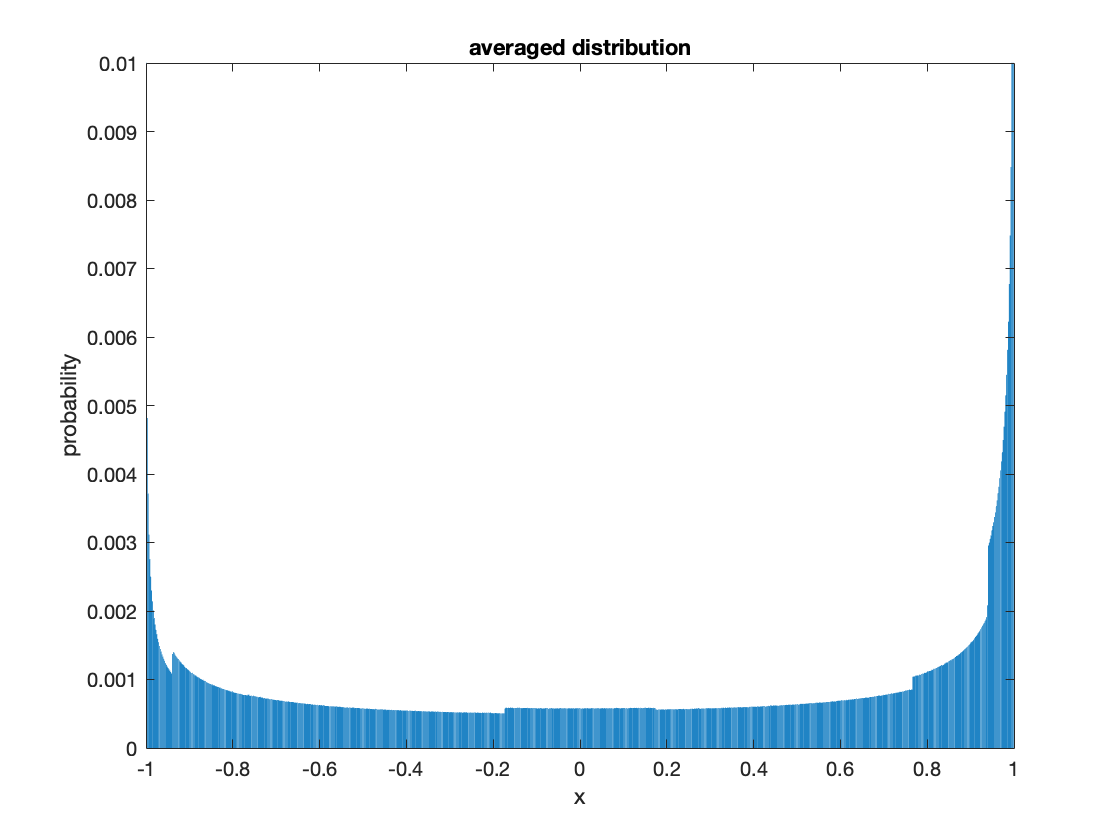}
}
\caption{Examples of different invariant densities of $T_{N, a}$. First row: graphs of the conjugated piecewise-linear maps $g_{N, a}$; second row: their invariant densities $\rho_g$; last row: invariant densities $\rho_T$ of the corresponding shifted Chebyshev maps. Values of $N$ and $a$ are indicated in each caption. Histograms (each with $1,000$ bins) are obtained from simulations: distribution averaged over $10,000$ trajectories each evolving for $11,000$(-1,000, to eliminate the initial transient stage) time steps; verified by analytical results.} 
\label{3egs-inv-den}
\end{figure}

For the simplest case, where $a = 0$, $T_{N, 0}$ is the ordinary Chebyshev map, whose conjugated piecewise-linear map is \textit{of full-branch}\footnote{A piecewise linear map is \textit{of full-branch} if, after a suitable  translation along the $x$-axis, each branch is mapped to the whole interval $[0, 1]$, making the Lebesgue measure preserved. A counterexample is shown in Fig.\ref{counter-eg-of-full-branch} for which the Lebesgue measure is not invariant. } and preserves the Lebesgue measure, $\rho_g = 1$, so we simply have $\rho_{T}(x) = \frac{1}{\pi \sqrt{1 - x^2}}$, $\forall x \in [-1, 1]$. See the first column in Fig.\ref{3egs-inv-den}. 

Furthermore, it can be easily argued that for 
\begin{equation}
\text{even $N$ with all $a \in \left[-\frac{\pi}{2}, 0\right]$, or odd $N$ with $a = 0$}, 
\tag{$\star$}
\label{star-cond}
\end{equation}
the conjugated piecewise-linear maps are of full-branch so $\rho_g = 1$ and the invariant density for $T_{N, a}$ is given by 
\begin{equation*}
\rho_{T}(x) = \frac{1}{\pi \sqrt{1 - x^2}}, \quad x \in [-1, 1].
\end{equation*}
Otherwise, $\rho_g$ is a non-trivial stepwise-constant function. Three examples are shown in the second, third and last column in Fig.\ref{3egs-inv-den}. The density in such cases  can be calculated by finding the corresponding Markov partition and the transfer matrix, see \ref{Amarkov} for a simple example ($N = 3, a = -\frac{\pi}{9}$). 

\section{Eigenfunctions of Perron-Frobenius (PF) operator for Chebyshev maps}
In this section we describe properties of the eigenfunctions and eigenvalues of the transfer operator for Chebyshev maps, for both the original maps and the shifted ones. These eigenfunctions are important when characterising non-equilibrium properties of the dynamics, i.e. the way how an arbitrary initial density approaches the invariant density.
If we write the initial distribution as a linear combination of eigenfunctions, then the time evolution of the density becomes simple, as the application of the (linear) Perron-Frobenius operator just corresponds to multiplying the eigenfunctions in the expansion of the initial density with the corresponding eigenvalues \cite{CBred, lasota-mackey}.

\subsection{Ordinary Chebyshev maps}
Recall that the Perron-Frobenius (PF) operator $\mathcal{L}$ for a map describes the time evolution of a set of points characterised by a density function \cite{BoyGora, CBred, BaladiGreenBook, pollicott}. Eigenfunctions $\rho$ and the corresponding eigenvalues $\lambda$ satisfy $\mathcal{L}\rho = \lambda \rho$, so that the invariant density $\rho^*$ is the eigenfunction that corresponds to the unit eigenvalue: $\mathcal{L}\rho^* = \rho^*$. For a one-dimensional discrete dynamical system $f$ we simply have 
\begin{equation*}
\mathcal{L}\rho(y) = \sum_{x \in f^{-1}(y)} \frac{\rho(x)}{|f'(x)|}.
\end{equation*}

It is known \cite{VepstasBern} that eigenfunctions of the PF operator for the binary shift are given by the Bernoulli polynomials $\rho^{(n)}(x) = B_n(x)$ with eigenvalues $\lambda^{(n)} = 2^{-n}$. More generally, eigenfunctions and eigenvalues for an $N$-ary shift are $\rho^{(n)}(x) = B_n(x)$, $\lambda^{(n)} = N^{-n}$ ($N \in \mathbb{N}_{\geq 2}$). 

\begin{figure}[H]
\centering
\subfloat[\footnotesize multi-upside-down tent (for even $N$)]{
\includegraphics[width = 0.48\textwidth]{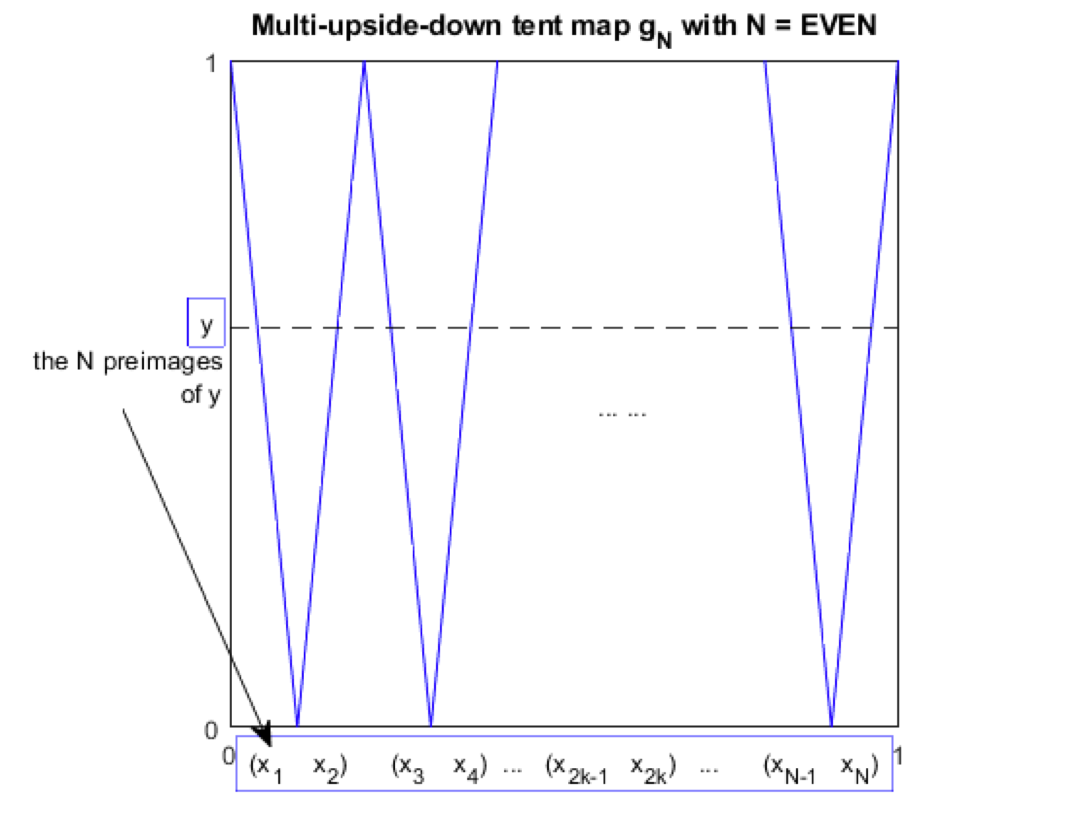}
\label{multi-even}}
\subfloat[\footnotesize multi-tent (for odd $N$)]{
\includegraphics[width = 0.48\textwidth]{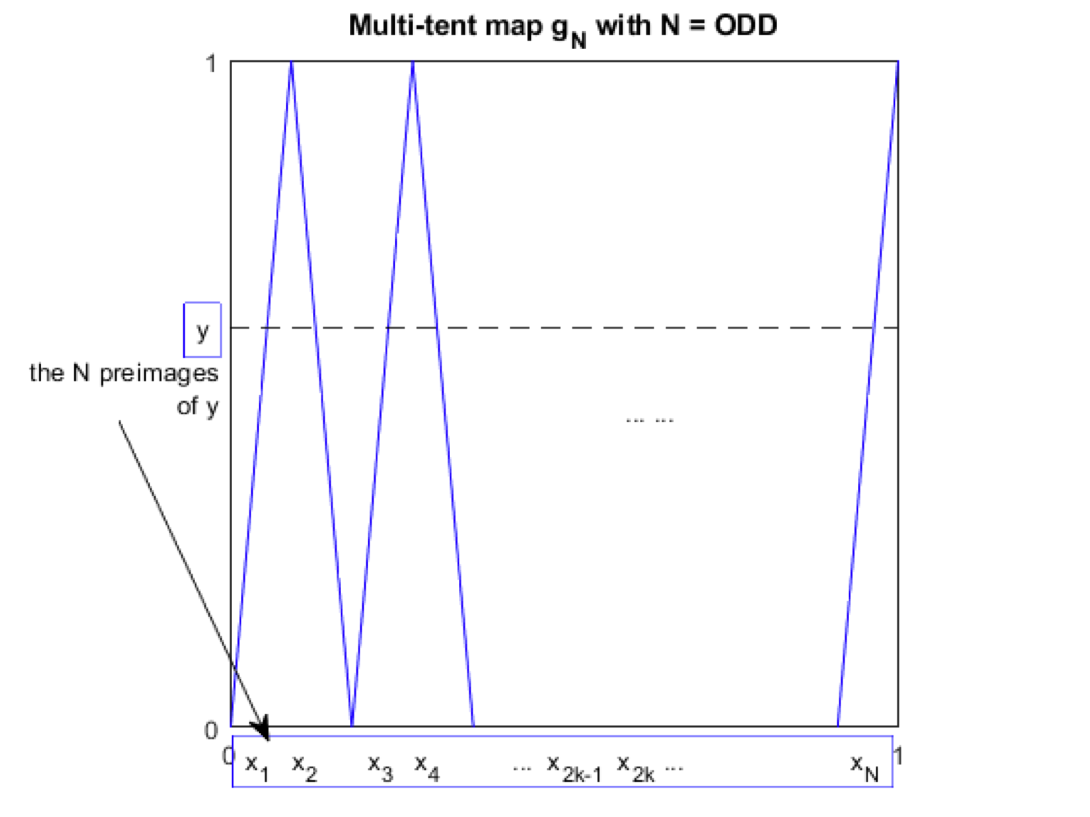}
\label{multi-odd}}
\caption{Graphs of piecewise-linear maps to which the ordinary Chebyshev maps conjugated.}
\label{multi2}
\end{figure}

We now use this result to find eigenfunctions for the Chebyshev maps. First consider the ordinary Chebyshev maps $T_N = T_{N, 0}$. Via the conjugacy mentioned in the previous section one sees that, for even $N = 2, 4, 6, ...,$ the piecewise-linear function $g_N$ consists of a multiple of upside-down tent maps (see Fig.\ref{multi-even}), while for odd $N = 3, 5, 7, ...$ the map $g_N$ consists of a multiple of tent maps (see Fig.\ref{multi-odd}); let us call them \textit{multi-upside-down tent} and \textit{multi-tent maps}, respectively. One can show (details in \ref{A2}), by the Multiplication Theorem and symmetry properties for Bernoulli and Euler polynomials \cite{Doerfle1985}, that the eigenfunctions are given by ($N \in \mathbb{N}_{\geq 2}$, $n \in \mathbb{N}_0$)
\begin{itemize}
\item for multi-upside-down tents $g_{N}$ (even $N$): $\rho^{(n)}_{g_{N}}(x) = B_{2n}\left(\frac{x}{2} + \frac{1}{2}\right)$, with $\lambda^{(n)} = N^{-2n}$;
\vspace{-0.3cm}
\item for multi-tents $g_{N}$ (odd $N$): $\rho^{(n, 1)}_{g_{N}}(x) = B_{2n}(x)$, $\rho^{(n, 2)}_{g_N}(x) = E_{2n-1}(x)$, with $\lambda^{(n)} = N^{-2n}$,
\end{itemize}
where $B_n$ and $E_n$ are Bernoulli and Euler polynomials, defined by their generating functions $G_B(x, t) = \frac{te^{xt}}{e^t - 1} = \sum_{n = 0}^{\infty}B_n(x)\frac{t^n}{n!}$ and $G_E(x, t) = \frac{2e^{xt}}{e^t + 1} = \sum_{n = 0}^{\infty}E_n(x)\frac{t^n}{n!}$, respectively. Notice that the eigenfunctions are independent of the order $N$ of the map.

By a topological coordinate transformation, $\rho_{T_N}(x) = h'(x)\rho_{g_N}(h(x))$, where $h(x) = \frac{1}{\pi}\arccos (-x)$, we get the following result: 
the eigenfunctions of the PF operator $\mathcal{L}$ for the ordinary Chebyshev maps $T_N = T_{N, 0}$ are given by ($N \in \mathbb{N}_{\geq 2}$, $n \in \mathbb{N}_0$)
\begin{equation}
\begin{split}
\mathcal{L}\rho_{T_N}^{(n)}(x) &= \lambda^{(n)}\rho_{T_N}^{(n)}(x) \\
\lambda^{(n)} &= N^{-2n},\\
\rho_{T_N}^{(n)}(x) &= \begin{cases}
\frac{1}{\pi \sqrt{1 - x^2}}B_{2n}\left(\frac{1}{2\pi}\arccos (-x) + \frac{1}{2}\right), \quad \text{if $N$ is even}; \\
\begin{cases}
\rho_{T_N}^{(n, 1)} = \frac{1}{\pi \sqrt{1 - x^2}}B_{2n}\left(\frac{1}{\pi}\arccos (-x)\right)\\
\rho_{T_N}^{(n, 2)} = \frac{1}{\pi \sqrt{1 - x^2}}E_{2n-1}\left(\frac{1}{\pi}\arccos (-x)\right), 
\end{cases}
\quad \text{if $N$ is odd; }
\end{cases}\\
&\text{where $B_n(x)$ and $E_n(x)$ are Bernoulli and Euler polynomials, respectively.}
\end{split}
\label{efun-TN}%
\end{equation}

Note that for odd $N$ the eigenfunctions are degenerate, i.e. two independent sets of eigenfunctions exist.
One can easily check that when $n = 0$ we have $\lambda^{(0)} = 1$ and $\rho^{(0)}_{T_N} = \frac{1}{\pi\sqrt{1 - x^2}}$, which is the invariant density of the ordinary Chebyshev maps for all $N$. 

It is intriguing that all eigenvalues are real, and the eigenfunctions are orthogonal polynomials. This reminds us of quantum mechanics and the corresponding properties of the eigenfunctions of the Schr\"odinger operator, although clearly the Perron-Frobenius operator is not a Hermitean operator.

\subsection{Shifted Chebyshev maps}
First consider maps within the category \eqref{star-cond}:

(i) any (even or odd) $N$ with $a = 0$: These are the ordinary Chebyshev maps, and we have already shown above that eigenfunctions are given by Bernoulli or/and Euler polynomials, with eigenvalues $N^{-2n}$; 

(ii) even $N$ with $a \in \left[ -\frac{\pi}{2}, 0\right)$: 
Assume the eigenvalues are the same as in (i), and let $N = 2q$ ($q \in \mathbb{N}$) and consider $a = -\frac{\pi}{m}$ with $m \in \mathbb{N}_{\geq 2}$. It can be verified that $T_{2q, -\frac{\pi}{m}}$ is topologically semi-conjugated to its corresponding ordinary map $T_{2q, 0}$ via the semi-conjugacy function $h_1 (x) = -T_{m, 0}(x)$, $x \in [-1, 1]$ and, moreover,
\begin{equation}
(-T_{m, 0}) \circ T_{2q, -\frac{\pi}{m}} = T_{2q, 0} \circ (-T_{m, 0}) = T_{2qm, 0}, 
\label{semi-conj1}%
\end{equation}
i.e.\ the combined action is equivalent to a single ordinary Chebyshev map (with a multiplicative order $N' = 2qm$).

From $h_1'(x) = -\frac{m\sin (m \arccos x)}{\sqrt{1 - x^2}}$ and eq.\eqref{efun-TN}, also taking into account that there are $m$ preimages of the semi-conjugacy $h_1$, we get 
\begin{equation*}
\rho_{T_{2q, -\frac{\pi}{m}}}^{(n)}(x) = \frac{|h'_1(x)|}{m}\rho_{T_{2q, 0}}^{(n)}(h_1(x)) = \frac{1}{\pi \sqrt{1 - x^2}}B_{2n}\left( \frac{m}{2\pi}\arccos x + \frac{1}{2}\right). 
\end{equation*}

The situation becomes more complicated for maps outside the category \eqref{star-cond}, and we will briefly discuss this case in \ref{AshiftedEfun}. 
Also, in the case that $a$ is not a rational multiple of $\pi$, the form of eigenfunctions is more complicated and in general exhibits some fractal features.

\section{Higher-order correlations of shifted Chebyshev maps} 
\subsection{Definition}
Generally, the \textit{$r$th-order correlation function} for a given map $x_{n+1}=T(x_n)$ is defined as 
\begin{equation*}
\langle x_{n_1} x_{n_2} \cdots x_{n_r} \rangle = \int_{-1}^1 \rho_T(x_0)x_{n_1} x_{n_2}\cdots x_{n_r}dx_0, 
\end{equation*}
where the average is taken over all initial points weighted with respect to the invariant density $\rho_T$. In the case of shifted Chebyshev maps $T$ in category \eqref{star-cond}, we obtain after some calculation 
\begin{equation}
\langle x_{n_1}\cdots x_{n_r} \rangle = 2^{-r} \sum_{\sigma} \left[ \exp\left(ia \sum_{l=1}^r \sigma_l \frac{N^{n_l} - 1}{N-1}\right) \cdot \delta \left( \sum_{l = 1}^r \sigma_l N^{n_l}, 0\right) \right], 
\label{shifted-rth-corr}%
\end{equation}
where the sum over $\sigma$ is a summation over all possible spin configurations $\sigma := (\sigma_1, \sigma_2, ..., \sigma_r)$, $\sigma_l \in \{ -1, +1\}$, and $\delta(x, 0)$ is the Kronecker delta defined as being 1 if $x=0$ and 0 else. Details in \ref{A3}.

\subsection{Ordinary Chebyshev maps}
In this case
\begin{equation*}
T_{N, 0}(x) = \cos (N \text{arccos}x), \quad x \in [-1, 1]
\end{equation*}
which can be written as polynomials: 
\begin{equation*}
\begin{split}
T_{2, 0}(x) &= 2x^2 - 1\\
T_{3, 0}(x) &= 4x^3 - 3x\\
T_{4, 0}(x) &= 8x^4 - 8x^2 + 1\\
T_{5, 0}(x) &= 16x^5 - 20x^3 + 5x\\
\vdots
\end{split}
\end{equation*}
The correlation functions reduce to 
\begin{equation}
\langle x_{n_1} x_{n_2} \cdots x_{n_r}\rangle _{a = 0}= \frac{1}{2^r}\sum_{\sigma} \delta \left( \sum_{l = 1}^r \sigma_l N^{n_l}, 0 \right).
\label{a0-rth-corr}%
\end{equation}
The non-vanishing correlations of $T_{N, 0}$ correspond to tuples $(n_1, ..., n_r)$ that solve the diophantine equations 
\begin{equation}
\sum_{l = 1}^r \sigma_l N^{n_l} = 0
\label{Diophantine}%
\end{equation}
for $\sigma_l \in \{ -1, +1\}$. For a given $N$, for each $r = 2, 3, ...$, the tuples that solve the above equations can be represented by simple graphs ($N$-ary double forests), and this graph-theoretical method was introduced in \cite{AHCB2001, CB1991}. 
It turns out that for all odd $N$ with odd $r$, the correlation for $T_{N, 0}$ vanishes identically.

\subsection{Distinctive property of higher-order correlations for Chebyshev maps}
We now come to the central result of this paper, namely the fact that Chebyshev maps are distinguished as having a minimum set of higher-order correlations. 
First, let us discuss why higher-order correlation functions are relevant, and why, for example, the 2-point correlation
is not sufficient.

Consider a variable defined by a sum of iterates of a given map $T$, $y_n := \sum_{j = 1}^n x_j$, $x_j = T (x_{j-1})$; the $r$th moment is then given by
\begin{equation*}
\langle y_n^r \rangle = \langle \sum_{j_1 = 1}^n \cdots \sum_{j_r = 1}^n x_{j_1} \cdots x_{j_r}\rangle = \sum_{j_1 = 1}^n \cdots \sum_{j_r = 1}^n \langle x_{j_1} \cdots x_{j_r}\rangle. 
\end{equation*}
These types of sums are motivated by deterministic diffusion processes and generalised versions of Central limit theorems \cite{mackey, tirnakli, bountis, tsallis}.

If we know all the higher-order correlation functions of $T$, we will have the knowledge of all moments of $y_n$. 
For any smooth observable defined by a function $f$ of the variable $y_n$, we can write down its Taylor expansion as $f(y_n) = \sum_j b_j y_n^j$; its average is therefore given by 
\begin{equation*}
\langle f(y_n)\rangle = \sum_j b_j \langle y_n^j\rangle. 
\end{equation*}
Hence, provided all the moments of $y_n$ (and thus higher-order correlations of $x_n$) are known, we can calculate the expectation of this arbitrary observable. 

Now, consider an arbitrary map $W$ that is conjugated to an $N$-ary shift. Assume that the iterates $w_n$ of $W$ can be written as 
\begin{equation*}
w_n = f(N^n u)
\end{equation*}
with $w_0 = f(u)$, where $f$ is some smooth periodic function (for Chebyshev maps this is simply $f(u) = \cos \pi u$), indicating that it is conjugated to an $N$-ary shift dynamics. If $f$ has the Fourier representation
\begin{equation*}
f(u) = \sum_{k = -\infty}^{\infty} a_k e^{i\pi ku}, 
\end{equation*}
then the $r$th-order correlation of iterates of $W$ evaluates to 
\begin{equation*}
\langle w_{n_1} w_{n_2} \cdots w_{n_r}\rangle = \sum_{k_1 = -\infty}^{\infty} \sum_{k_2 = -\infty}^{\infty} \cdots \sum_{k_r = -\infty}^{\infty} a_{k_1}a_{k_2} \cdots a_{k_r} \delta (k_1 N^{n_1} + k_2 N^{n_2} + \cdots + k_r N^{n_r}, 0),
\end{equation*}
by a simple generalisation of the derivation that we presented in \ref{A3}. Non-zero correlations occur if a given tuple $(n_1, n_2, ..., n_r)$ solves any of the diophantine equations
\begin{equation*}
\sum_{j = 1}^r k_j N^{n_j} = 0, \quad k_j \in \mathbb{Z}. 
\end{equation*}
These equations in general have much more solutions than those for Chebyshev maps, as in the latter case each coefficient only takes one of the two integer values $k_j \in \{ -1, +1\}$,
whereas in this more general case the $k_j$ can take on any integer values. If we regard the original Chebyshev correlations as being described by an Ising model (up-down spin configurations), then this would be a generalisation towards a Potts model (integer spin configurations) \cite{CBred}.

The number of tuples $(n_1, \ldots , n_r)$ with non-zero correlations is minimised when the underlying map is from the Chebyshev family, in which case  only the coefficients $a_{\pm 1}$ are non-zero ($a_{-1} = a_{1} = \frac{1}{2}$) in the Fourier representation of the conjugating function $f(u) = \cos \pi u = \frac{1}{2}(e^{-i\pi u} + e^{i\pi u}) = a_{-1}e^{-i \pi u} + a_{1}e^{i \pi u}$. Hence, Chebyshev maps can be regarded as producing a ``minimum skeleton" of higher-order correlations.
They can serve as the most random-like deterministic system in this context, in the sense of a minimum set of correlations, or strongest possible similarity to white noise for a smooth deterministic dynamics.

\subsection{Example: two-point correlations} 
Our result of the previous subsection implies that the higher-order correlation functions for (ordinary) Chebyshev maps are identically equal to 0 for more tuples $(n_1, \ldots , n_r)$ than those for other $N$-ary shift systems, indicating that, although they are equally chaotic in the sense of topology (i.e., topologically conjugated to each other), one appears more random than the other in the sense of vanishing correlations. Let us now consider the special case $r=2$. It is known (see \ref{A4}) that the two-point correlation function for the binary shift map (normalised, i.e., with average $\frac{1}{2}$ subtracted) is given by
\begin{equation*}
\langle w_{n_1}w_{n_2}\rangle_{\text{binary}} = \frac{1}{12}\cdot \left(\frac{1}{2}\right)^{|n_1 - n_2|}, 
\end{equation*}
while for the $2$nd-order ordinary Chebyshev map $T_{2, 0}$ it vanishes immediately, as there is no solution to the corresponding diophantine equations, i.e. $\langle x_{n_1}x_{n_2}\rangle_{T_{2, 0}} = \frac{1}{2}\delta_{n_1, n_2}$. 

The two-point correlation function for a general $N$-ary shift (normalised) can be shown (\ref{A4}) to be
\begin{equation*}
\langle w_{n_1}w_{n_2}\rangle_{N-\text{ary}} = \frac{1}{6N}\cdot \left(\frac{1}{N}\right)^{|n_1 - n_2|}. 
\end{equation*}
Comparing with the two-point correlation function for the $N$th-order ordinary Chebyshev map $T_{N, 0}$, for which 
\begin{equation*}
\langle x_{n_1}x_{n_2}\rangle_{T_{N, 0}} = \frac{1}{2}\delta_{n_1, n_2}, \quad \forall N \geq 2, 
\end{equation*}
we see that for $N$-ary shifts it decreases exponentially in $N$ but it never vanishes identically, while for the Chebyshev maps the two-point correlation does not depend on $N$ and it attains zero whenever $n_1 \neq n_2$.  

\begin{figure}[H]
\centering
\subfloat[\footnotesize $\langle x_0 x_k\rangle_{T_{N, -\pi/2}}$, maximum is $\frac{\pi}{2}$ at $k = 0$ for all $N$]{
\includegraphics[width = 0.9\textwidth]{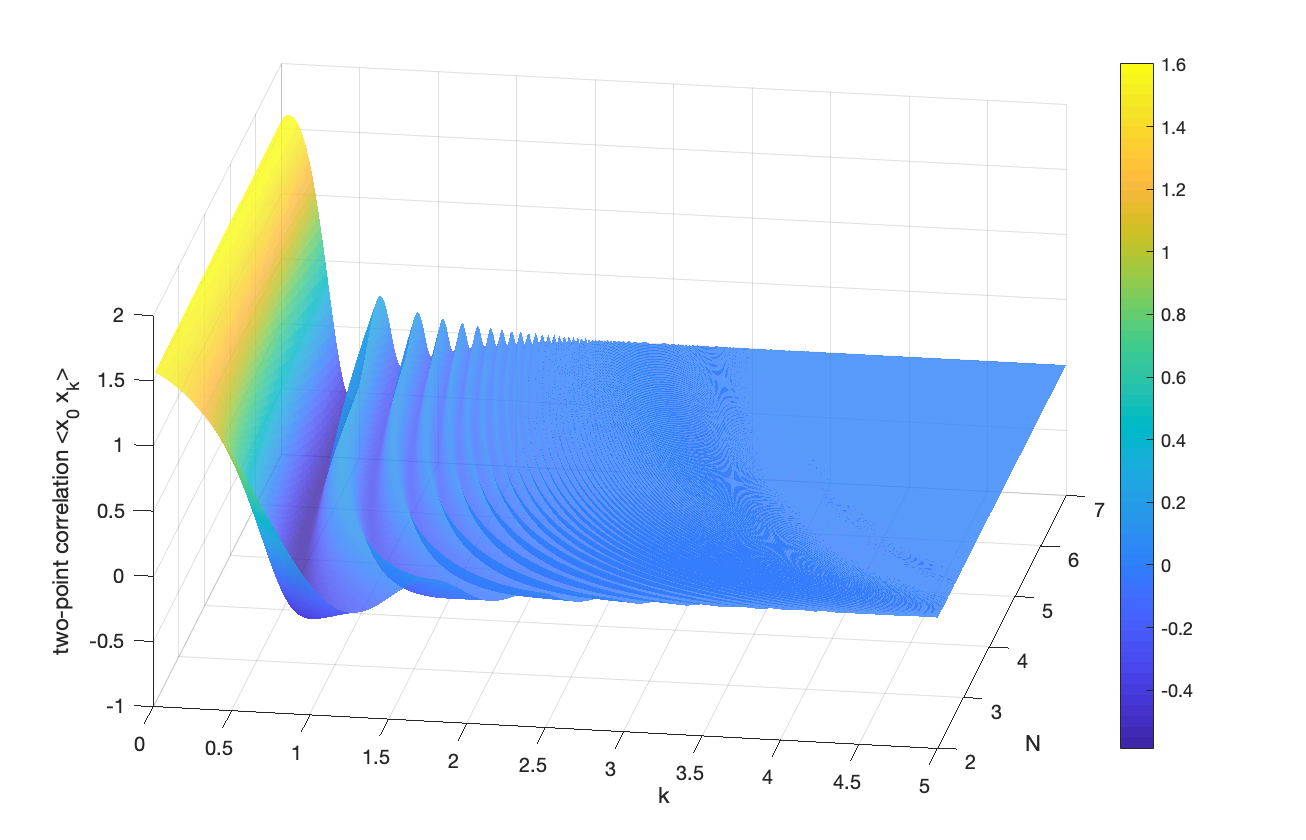}
}

\subfloat[\footnotesize $N = 2$]{
\includegraphics[width = 0.32\textwidth]{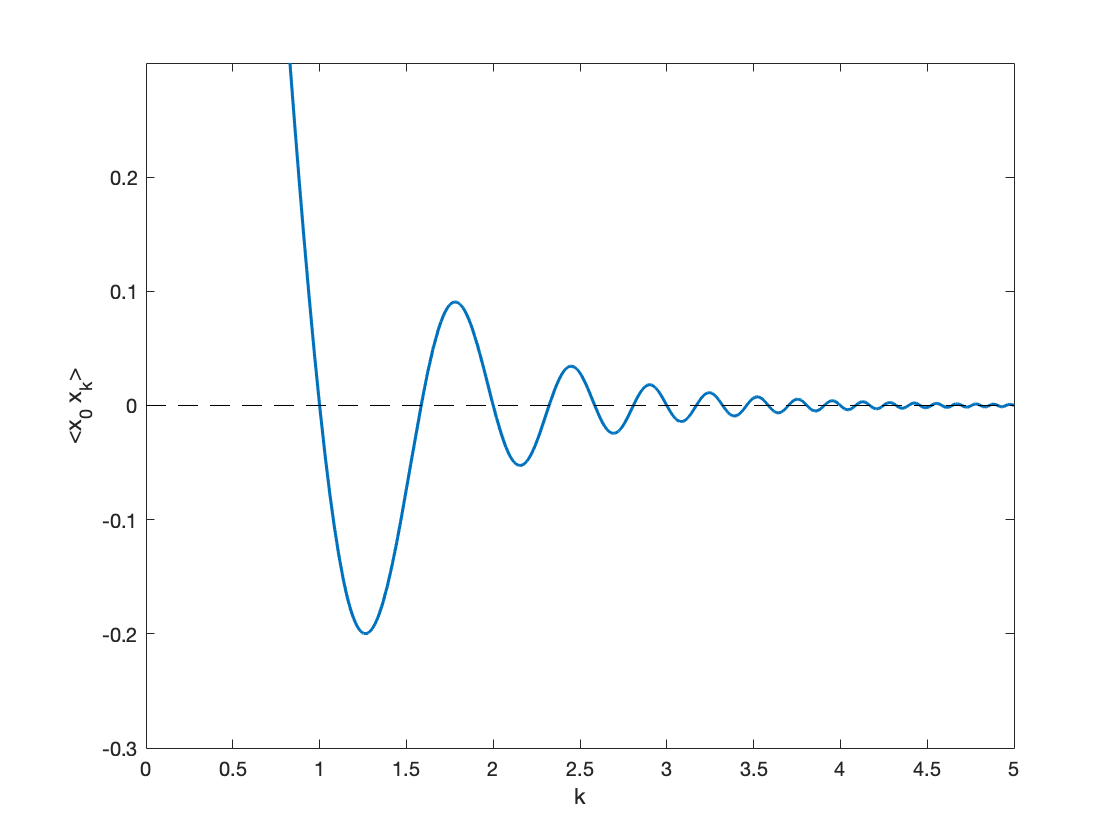}
}
\subfloat[\footnotesize $N = 3$]{
\includegraphics[width = 0.32\textwidth]{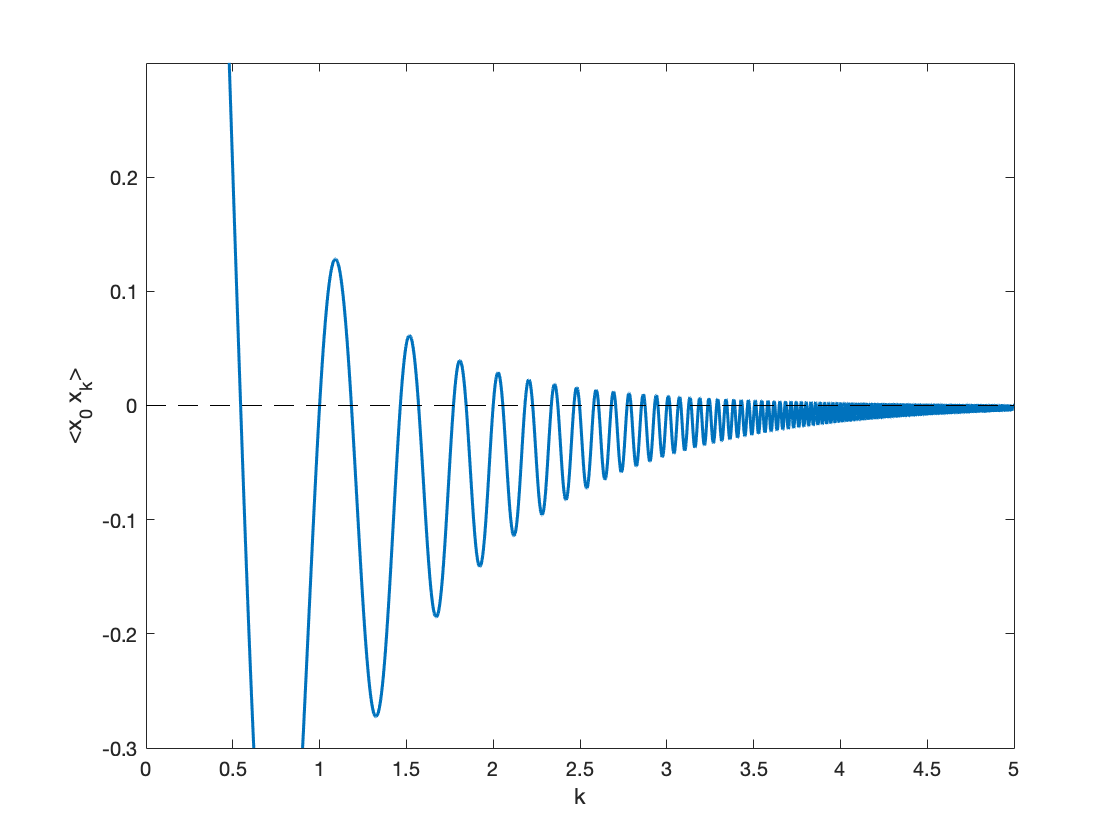}
}
\subfloat[\footnotesize $N = 4$]{
\includegraphics[width = 0.32\textwidth]{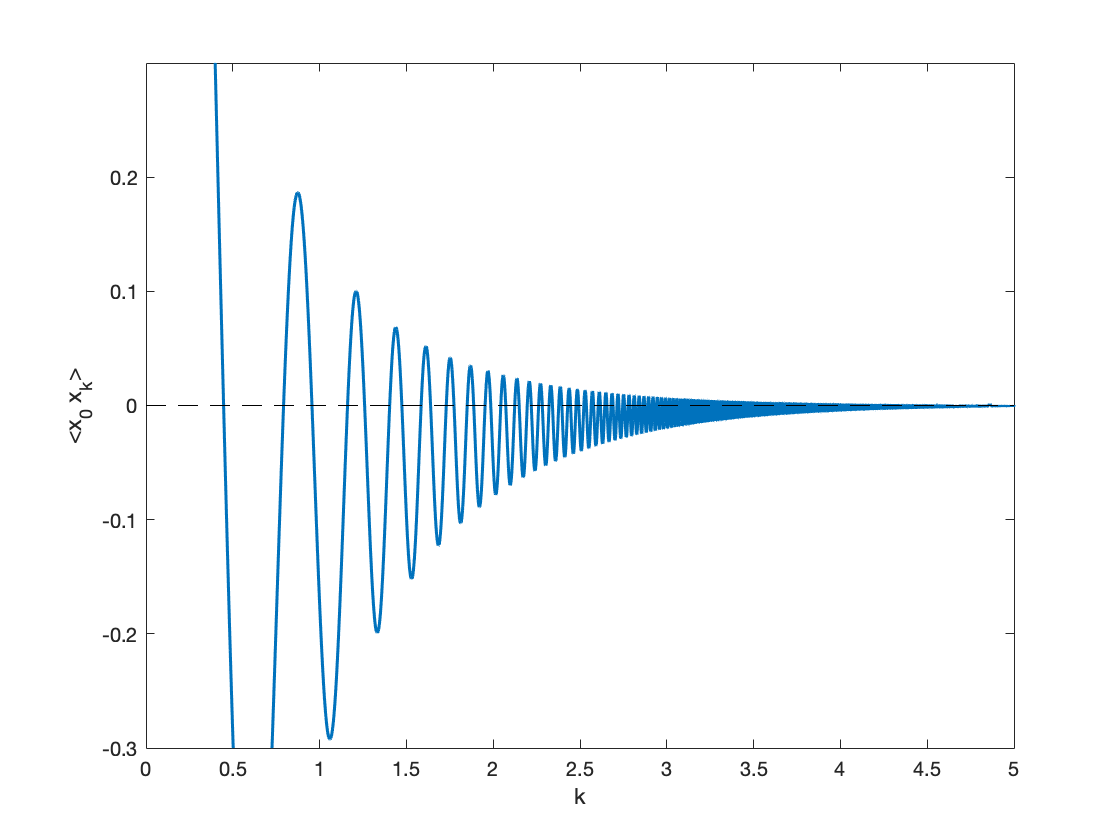}
}
\caption{Decay of two-point correlations for $T_{N, -\frac{\pi}{2}}$. Upper row: two-point correlation $\langle x_0 x_k\rangle_{T_{N, -\pi/2}}$ as a function of $(k, N) \in [0, 5]\times [2, 7]$. Lower row: for maps $T_{2, -\frac{\pi}{2}}$, $T_{3, -\frac{\pi}{2}}$, and $T_{4, -\frac{\pi}{2}}$, respectively. }
\label{two-pt-corr-sineCheby-N2N3N4}
\end{figure}

For shifted Chebyshev maps $T_{N, a}$ in the category \eqref{star-cond}, we can write down the two-point correlation functions explicitly (assuming stationarity)
\begin{equation*}
\begin{split}
&\langle x_{n}x_{n+k}\rangle_{T_{N, a}}\\
=& \langle x_0 x_k \rangle_{T_{N, a}}\\
=& \int_{-1}^1 \rho (x_0) x_0 T_{N, a}^{(k)}(x_0)dx_0\\
=& \int_{-1}^1 \frac{1}{\pi \sqrt{1 - \cos^2 \left(u_0' - \frac{a}{N-1}\right)}} \cos \left( u_0' - \frac{a}{N-1}\right) \cos \left( N^k u_0' - \frac{a}{N-1}\right) d\left[ \cos \left( u_0' - \frac{a}{N-1}\right)\right]\\
=& \frac{1}{2} \left[ \frac{\sin \left( \frac{2a}{N-1} - N^k \pi\right) + \sin \left( \frac{2a}{N-1}\right)}{N^k + 1} - \frac{\sin \left( N^k \pi\right)}{N^k - 1}\right].
\end{split}
\end{equation*}
One can check that for

i) $a = 0$, the ordinary Chebyshev maps, we have identically vanishing correlation for all $k \neq 0$; otherwise we have the second moment $\langle x_0^2\rangle_{T_{N, 0}} = \int_{-1}^1\rho(x)x^2dx = \int_{-1}^1 \frac{x^2}{\pi \sqrt{1 - x^2}} dx = \frac{1}{2}$. 

ii) $a = -\frac{\pi}{2}$, we have, for $k = 0$ 
\begin{equation*}
\langle x_0^2\rangle_{T_{N, -\pi/2}} = \int_0^{\pi} \cos^2 \left( u - \frac{a}{N-1}\right) du = \frac{\pi}{2}, 
\end{equation*}
and for $k > 0$
\begin{equation*}
\langle x_0 x_k \rangle_{T_{N, -\pi/2}} = -\frac{1}{2}\left[ \frac{\sin \left( \frac{\pi}{N-1} + N^k \pi\right) + \sin \left( \frac{\pi}{N-1}\right)}{N^k + 1} + \frac{\sin \left( N^k \pi\right)}{N^k - 1}\right], 
\end{equation*}
which exhibits an oscillating exponential decay in $k$, cf. Fig.\ref{two-pt-corr-sineCheby-N2N3N4}. 

From the results above we conclude again that already at the 2-point level, the ordinary Chebyshev maps have the lowest correlations among all the maps in the shifted Chebyshev family, and in fact
also the lowest correlations as compared to any other map conjugated to an $N$-ary shift.

\section{Coupled map lattices (CMLs) of shifted Chebyshev maps}
\subsection{CMLs of two sites}
We now study spatially coupled systems \cite{dettmann, groote, chaotic-string, kaneko} and investigate how the correlation patterns are modified. As the simplest model, consider a (periodic) lattice with just two sites that consists of two coupled 2nd-order ordinary Chebyshev maps: 
\begin{equation*}
\begin{split}
x_{n+1}^{(1)} &= (1-c)T_{2, 0}(x_n^{(1)}) + cT_{2, 0}(x_n^{(2)})\\
x_{n+1}^{(2)} &= (1-c)T_{2, 0}(x_n^{(2)}) + cT_{2, 0}(x_n^{(1)}).
\end{split}
\end{equation*}
The superscript, 1 or 2, labels the spatial lattice position and the subscript $n \in \mathbb{N}_0$ denotes a discrete time step; $c \in [0, 1]$ is the coupling strength, and $T_{2, 0}(x) = \cos (2\arccos x) = 2x^2 - 1$, $x \in [-1, 1]$ is the 2nd-order ordinary Chebyshev map. 

For the uncoupled case ($c = 0$) the invariant density of the system is just the direct product of the two individual densities; as $c$ increases the invariant density will gradually shrink to 
a support given by the diagonal, that is, a total synchronisation state $x^{(1)} = x^{(2)}$ is approached. See Fig.\ref{2cml-invden-plots} below. 
\begin{figure}[H]
\centering 
\subfloat[\footnotesize $c = 0.027$]{
\includegraphics[width = 0.24\textwidth]{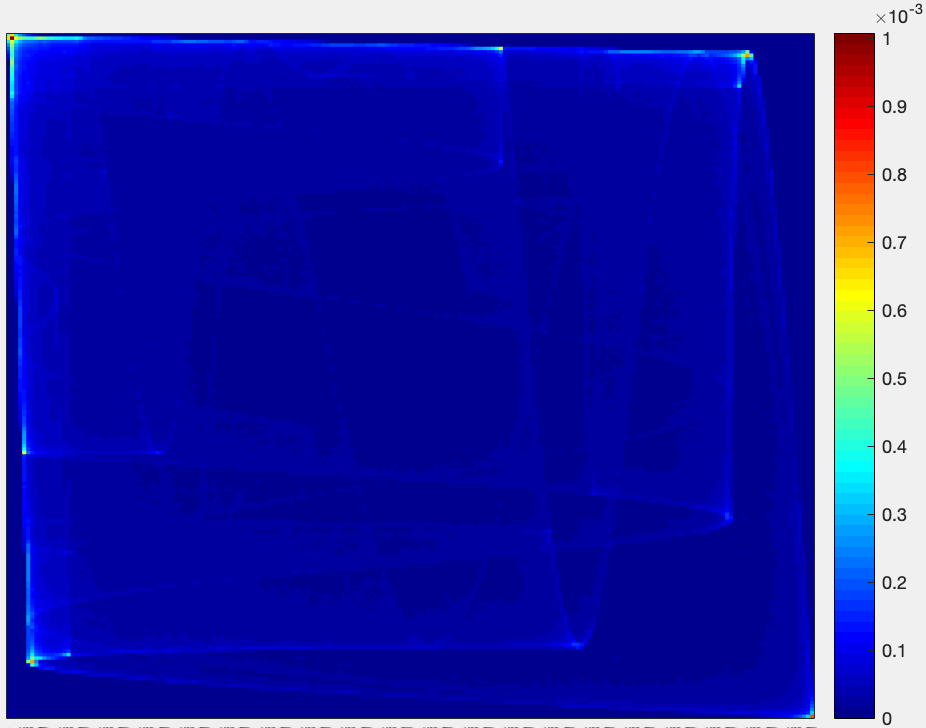}
}
\subfloat[\footnotesize $c = 0.1$]{
\includegraphics[width = 0.24\textwidth]{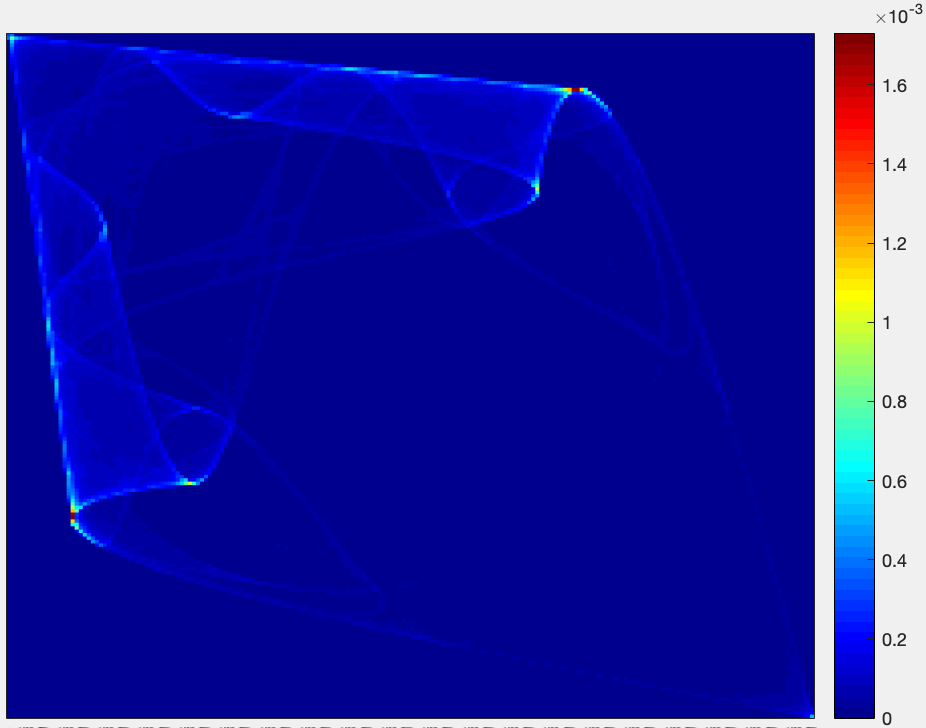}
}
\subfloat[\footnotesize $c = 0.11$]{
\includegraphics[width = 0.24\textwidth]{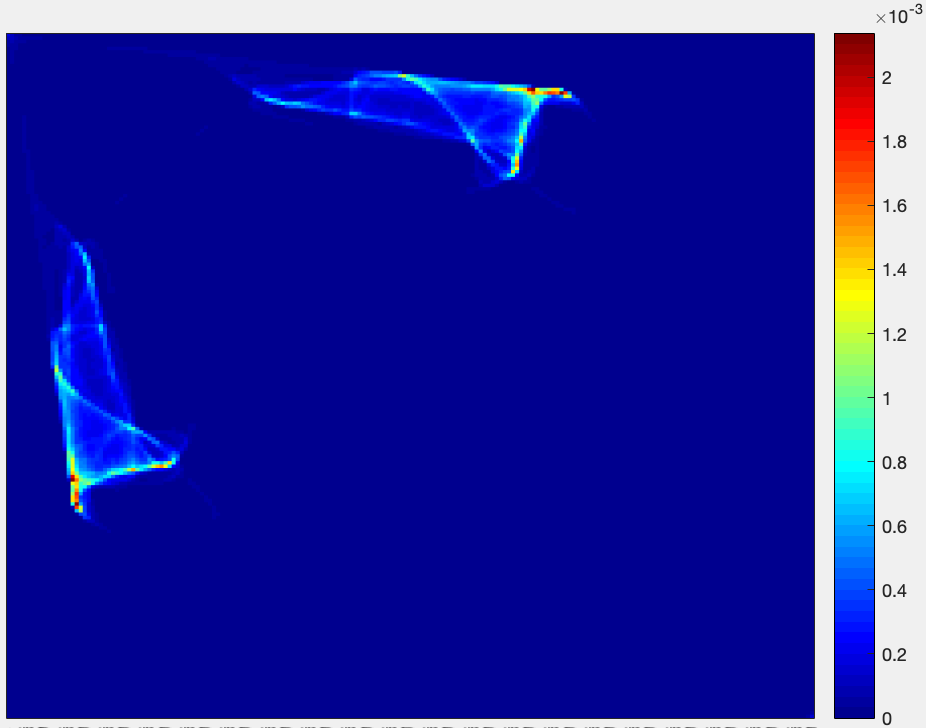}
}
\subfloat[\footnotesize $c = 0.201$]{
\includegraphics[width = 0.24\textwidth]{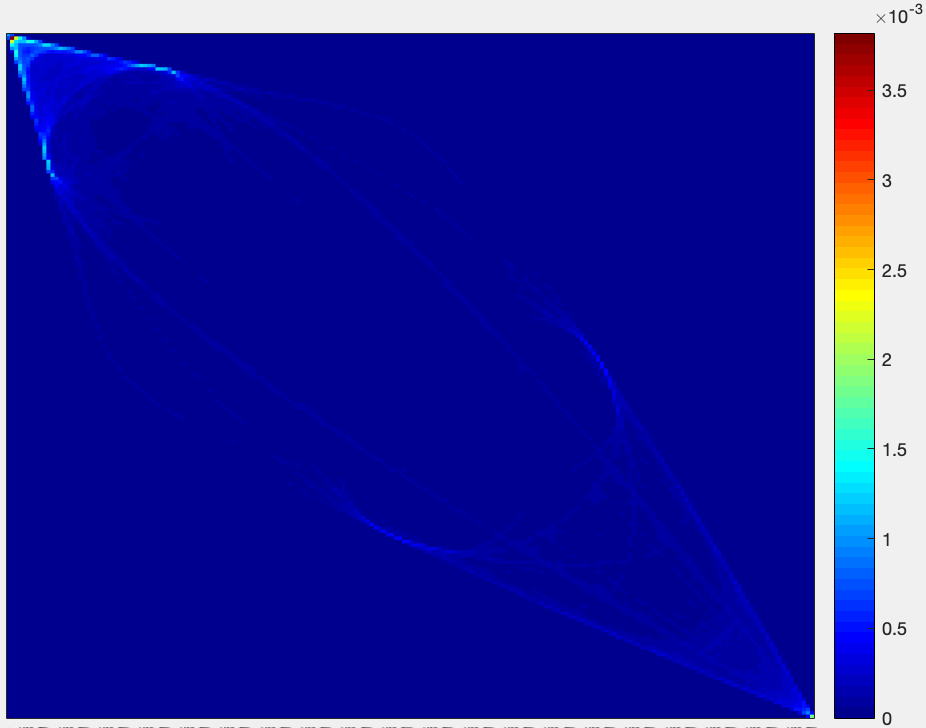}
}
\caption{Heatmaps of the invariant distribution for CMLs of two sites, with different couplings indicated in each caption. The $x$ and $y$ axes are $x^{(1)}$ and $x^{(2)}$, respectively (interchangeable). Colour codes the averaged probability that $10,000$ trajectories (randomly chosen initially) visit a given region over $10,100 (-100)$ iterations (dark blue $=$ low probability, $200\times 200$ bins). }
\label{2cml-invden-plots}
\end{figure}

\subsection{CMLs of many sites}
For applications in quantum field theory and high energy physics it is meaningful to extend to a CML of many sites. Let $T$ be a one-dimensional chaotic map, and consider the following four types of coupling \cite{CBblue, chaotic-string}: 

\underline{Type $A$}: forward, diffusive coupling
\begin{equation*}
x_{n+1}^{(i)} = (1-c)T(x_n^{(i)}) + \frac{c}{2}\left( T(x_n^{(i-1)}) + T(x_n^{(i+1)})\right)
\end{equation*}

\underline{Type $A^-$}: forward, anti-diffusive
\begin{equation*}
x_{n+1}^{(i)} = (1-c)T(x_n^{(i)}) - \frac{c}{2}\left( T(x_n^{(i-1)}) + T(x_n^{(i+1)})\right)
\end{equation*}

\underline{Type $B$}: backward, diffusive
\begin{equation*}
x_{n+1}^{(i)} = (1-c)T(x_n^{(i)}) + \frac{c}{2}\left( x_n^{(i-1)} + x_n^{(i+1)}\right)
\end{equation*}

\underline{Type $B^-$}: backward, anti-diffusive
\begin{equation*}
x_{n+1}^{(i)} = (1-c)T(x_n^{(i)}) - \frac{c}{2}\left( x_n^{(i-1)} + x_n^{(i+1)}\right)
\end{equation*}
The superscripts $i$ and subscripts $n$ of the dynamical variable $x$ denote the spatial position of the lattice site and the number of iterations, respectively; $c \in [0, 1]$ is the coupling strength, and $T$ is the local map to be specified. 

In \cite{CBblue} the local map $T$ was chosen as an ordinary Chebyshev map. To start with, consider the 2nd-order ordinary Chebyshev map $T_{2, 0} (x) = \cos (2\arccos x)$, $x \in [-1, 1]$. 
Fig.\ref{spatio-temporal-T2} below shows some spatio-temporal patterns of the four types of CMLs for this local $T = T_{2, 0}$. 
As a comparison, similar plots for local $T = T_{3, -\frac{\pi}{2}}$ are shown in Fig.\ref{spatio-temporal-T3sine}. 

\begin{figure}[H]
\centering 
\subfloat[\footnotesize Type $2A$: $c = 0$]{
\includegraphics[width = 0.32\textwidth]{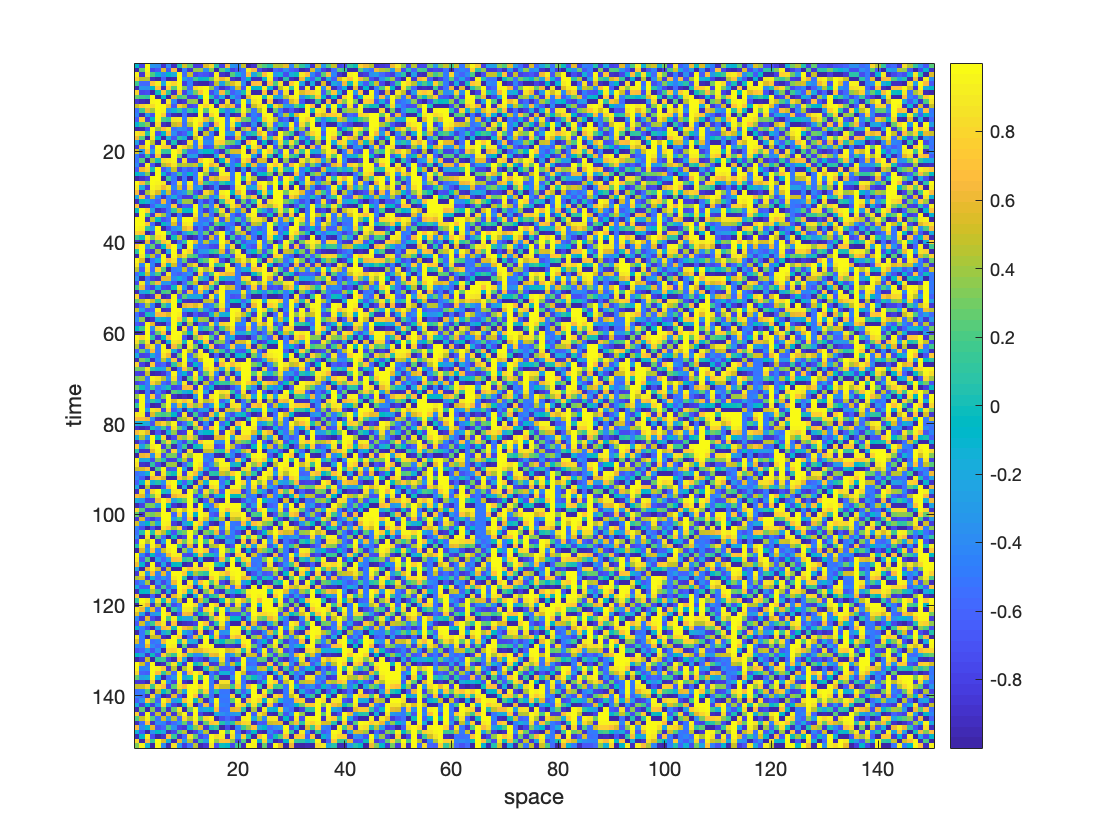}
}
\subfloat[\footnotesize Type $2A$: $c = 0.45$]{
\includegraphics[width = 0.32\textwidth]{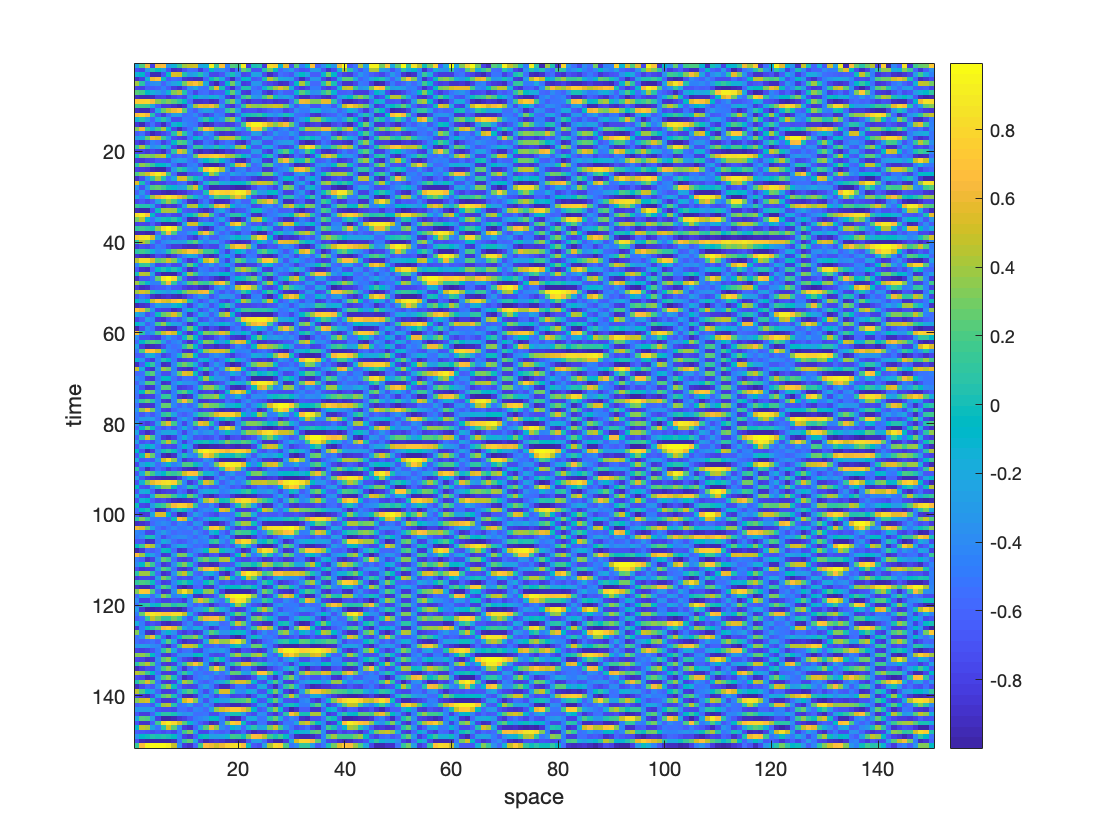}
}
\subfloat[\footnotesize Type $2A$: $c = 0.95$]{
\includegraphics[width = 0.32\textwidth]{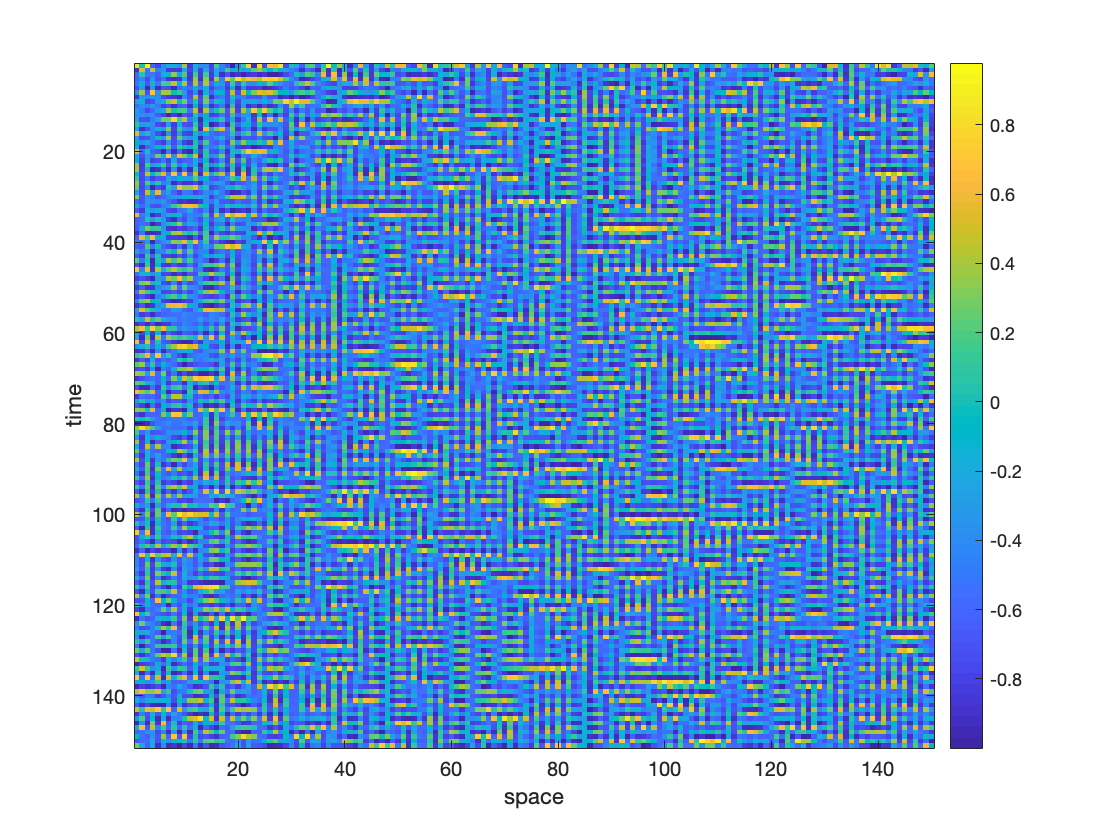}
}

\subfloat[\footnotesize Type $2A^-$: $c = 0.1$]{
\includegraphics[width = 0.32\textwidth]{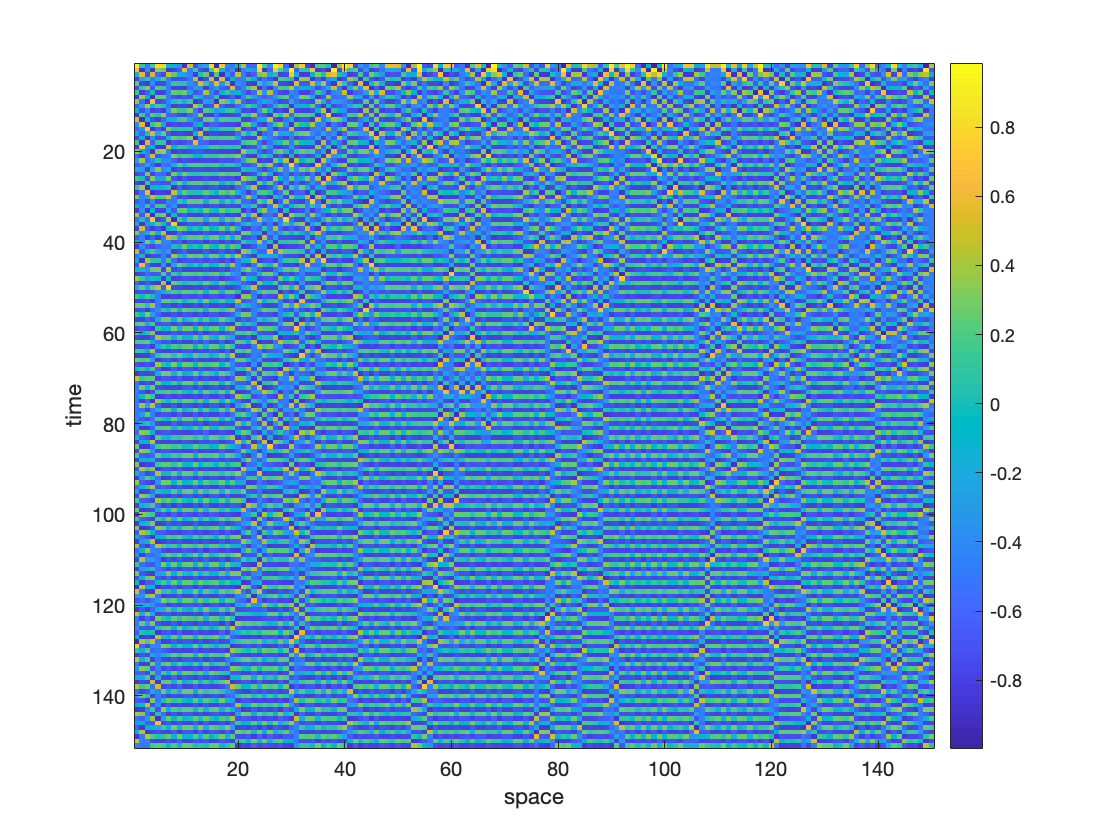}
}
\subfloat[\footnotesize Type $2A^-$: $c = 0.76$]{
\includegraphics[width = 0.32\textwidth]{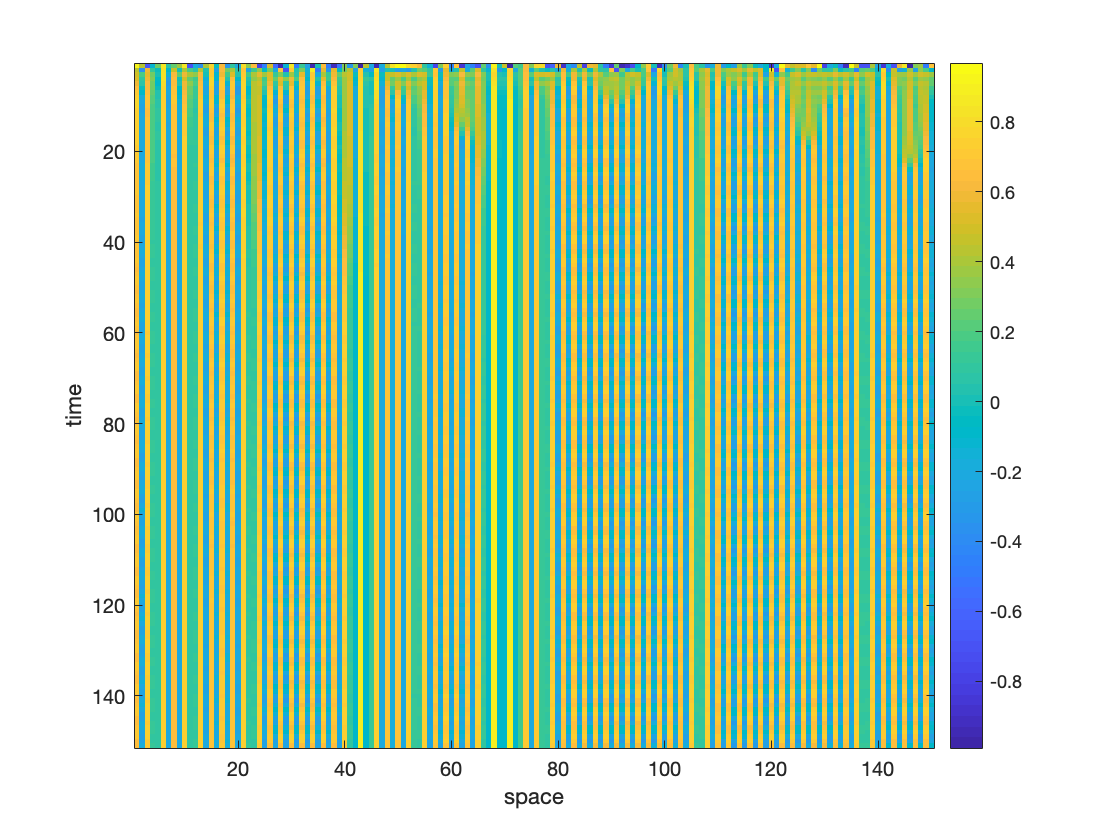}
}
\subfloat[\footnotesize Type $2A^-$: $c = 0.98$]{
\includegraphics[width = 0.32\textwidth]{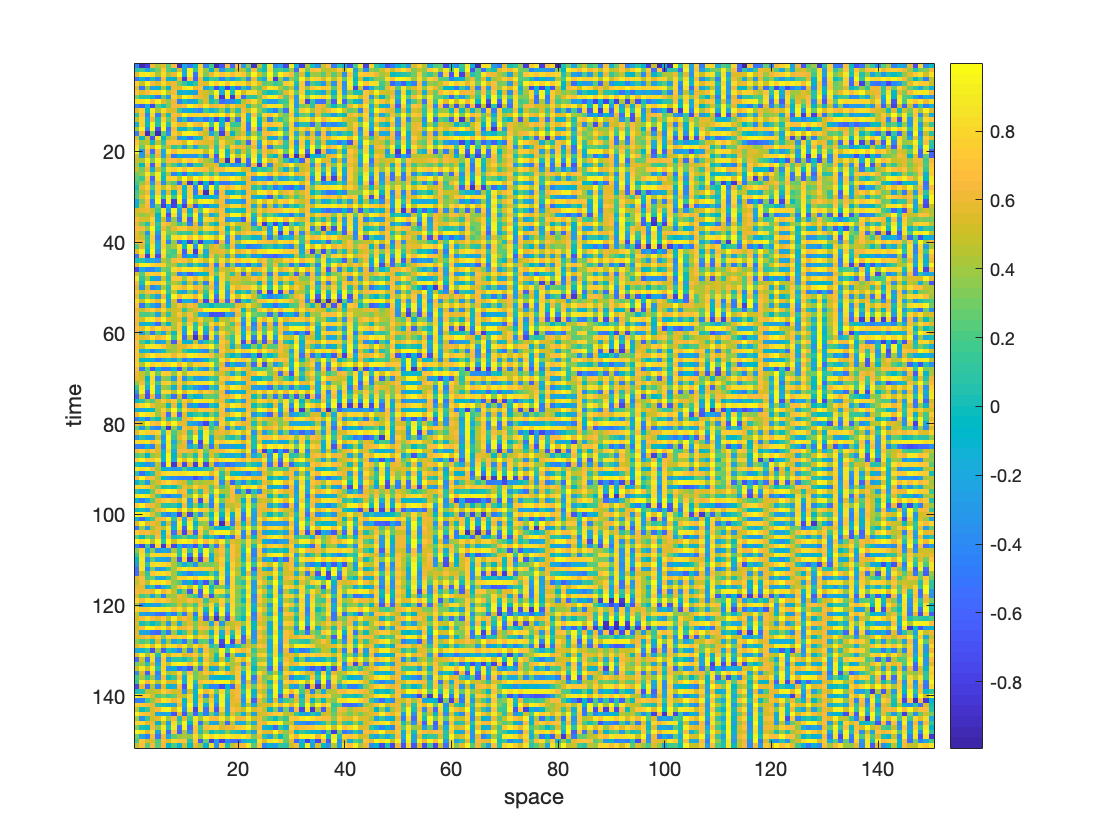}
}

\subfloat[\footnotesize Type $2B$: $c = 0.13$]{
\includegraphics[width = 0.32\textwidth]{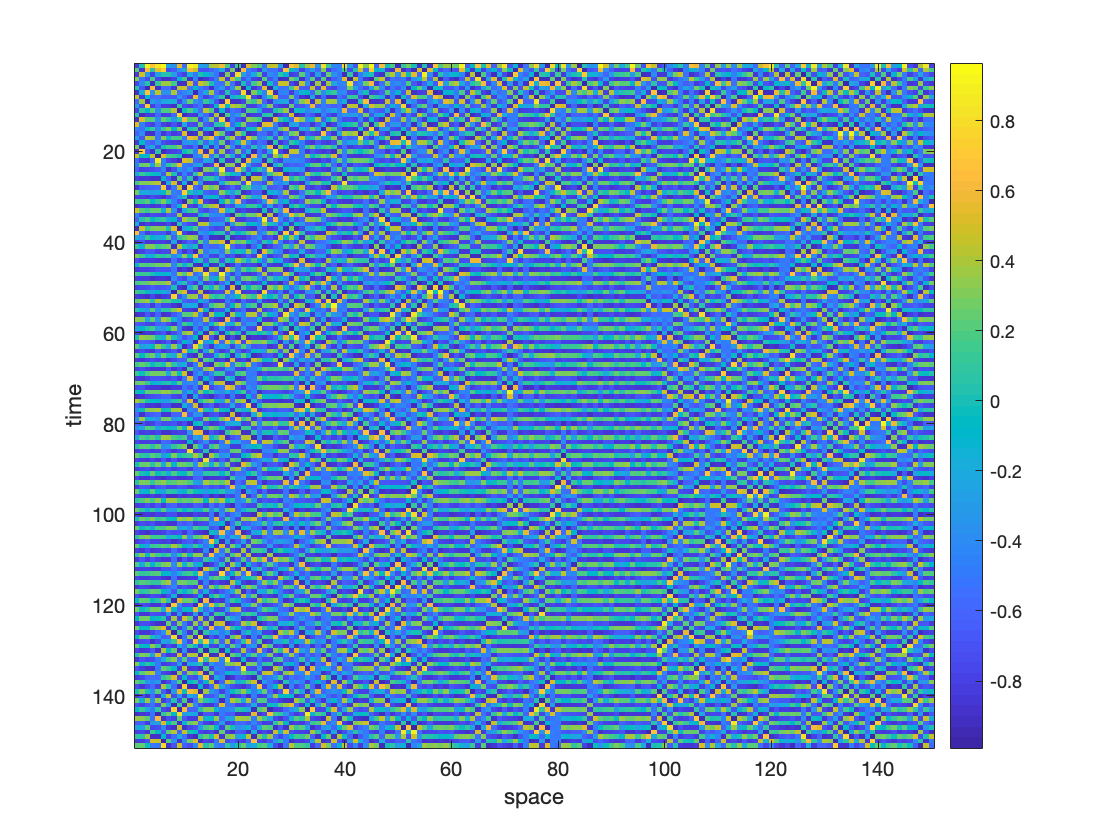}
\label{turb}
}
\subfloat[\footnotesize Type $2B$: $c = 0.79$]{
\includegraphics[width = 0.32\textwidth]{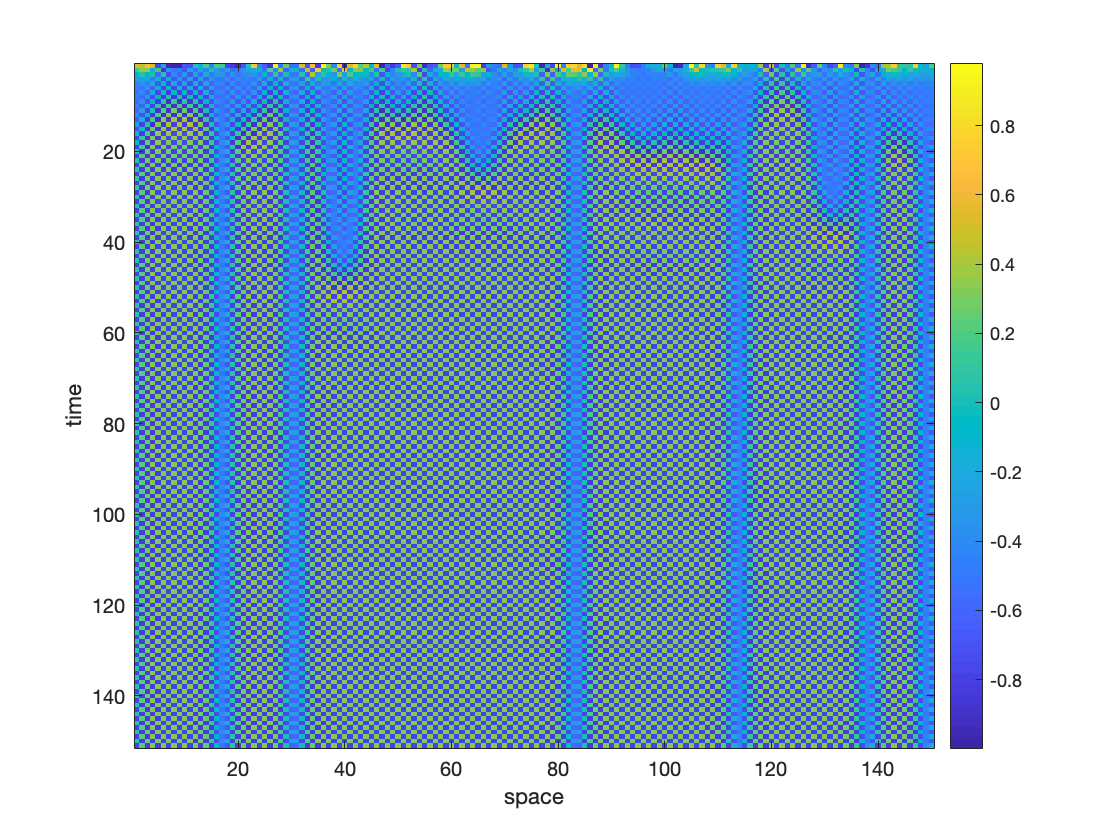}
}
\subfloat[\footnotesize Type $2B$: $c = 0.97$]{
\includegraphics[width = 0.32\textwidth]{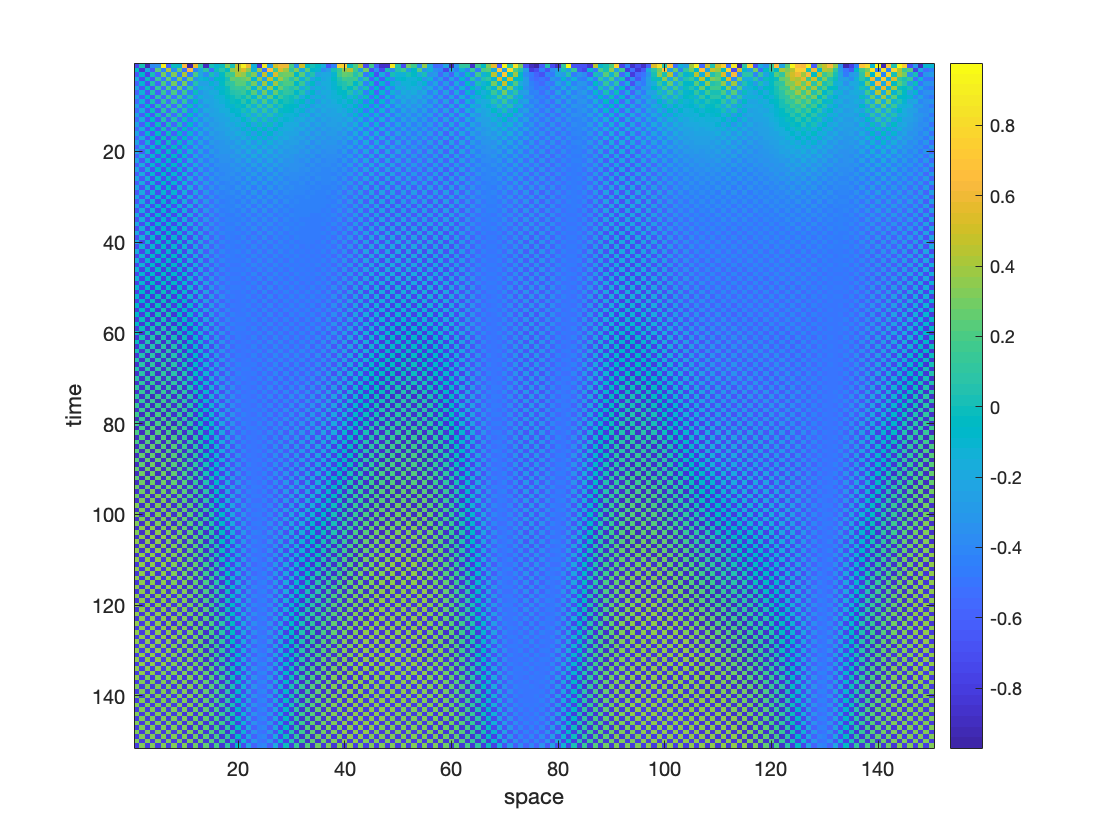}
}

\subfloat[\footnotesize Type $2B^-$: $c = 0.23$]{
\includegraphics[width = 0.32\textwidth]{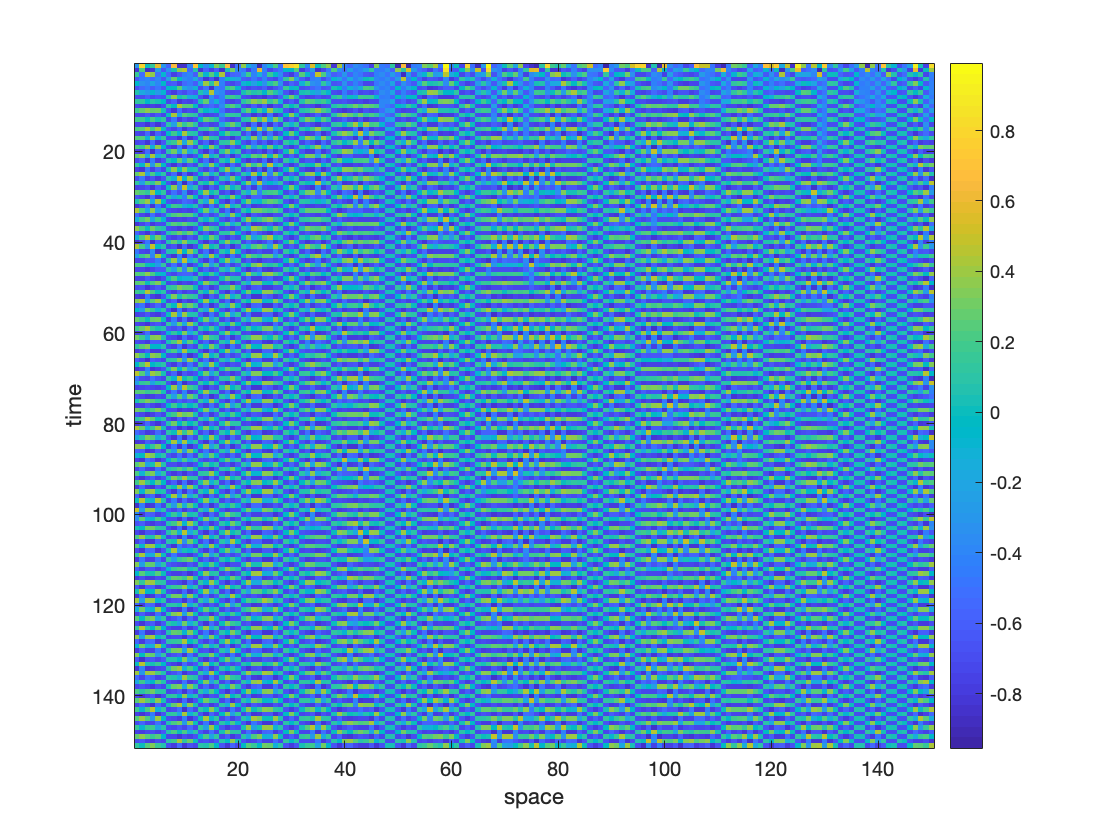}
\label{patt-select}
}
\subfloat[\footnotesize Type $2B^-$: $c = 0.47$]{
\includegraphics[width = 0.32\textwidth]{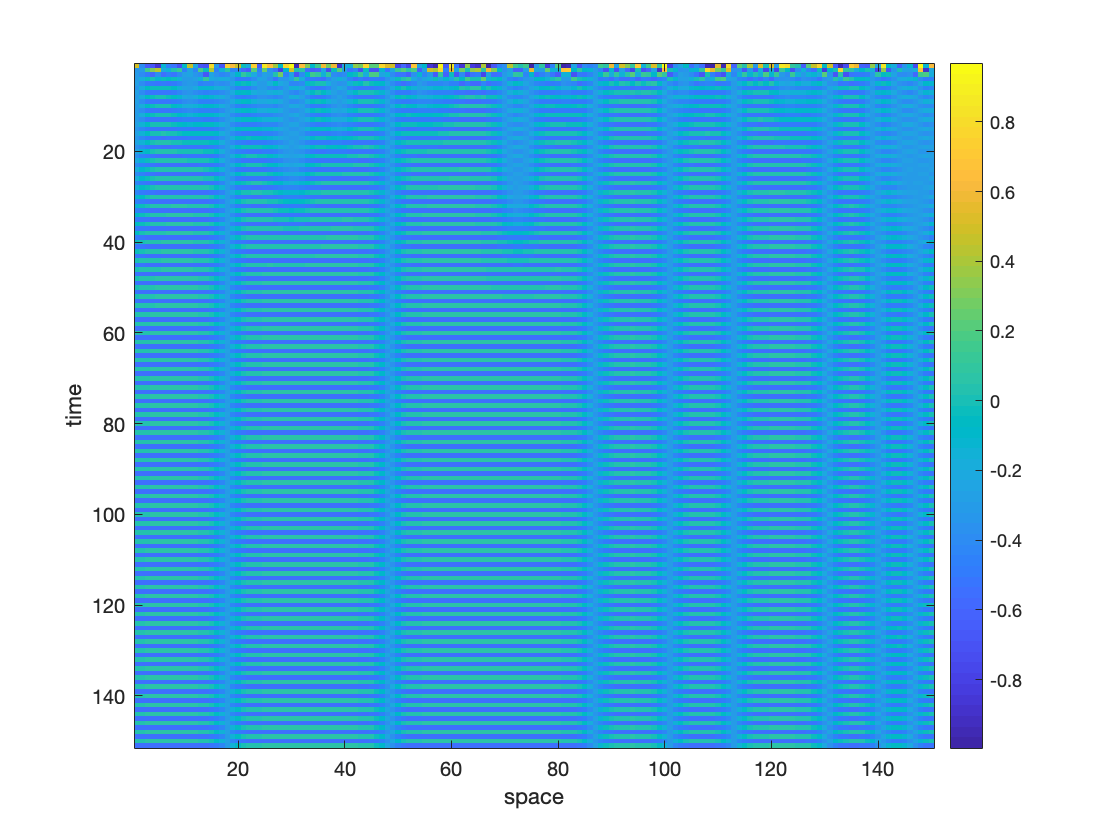}
}
\subfloat[\footnotesize Type $2B^-$: $c = 0.9$]{
\includegraphics[width = 0.32\textwidth]{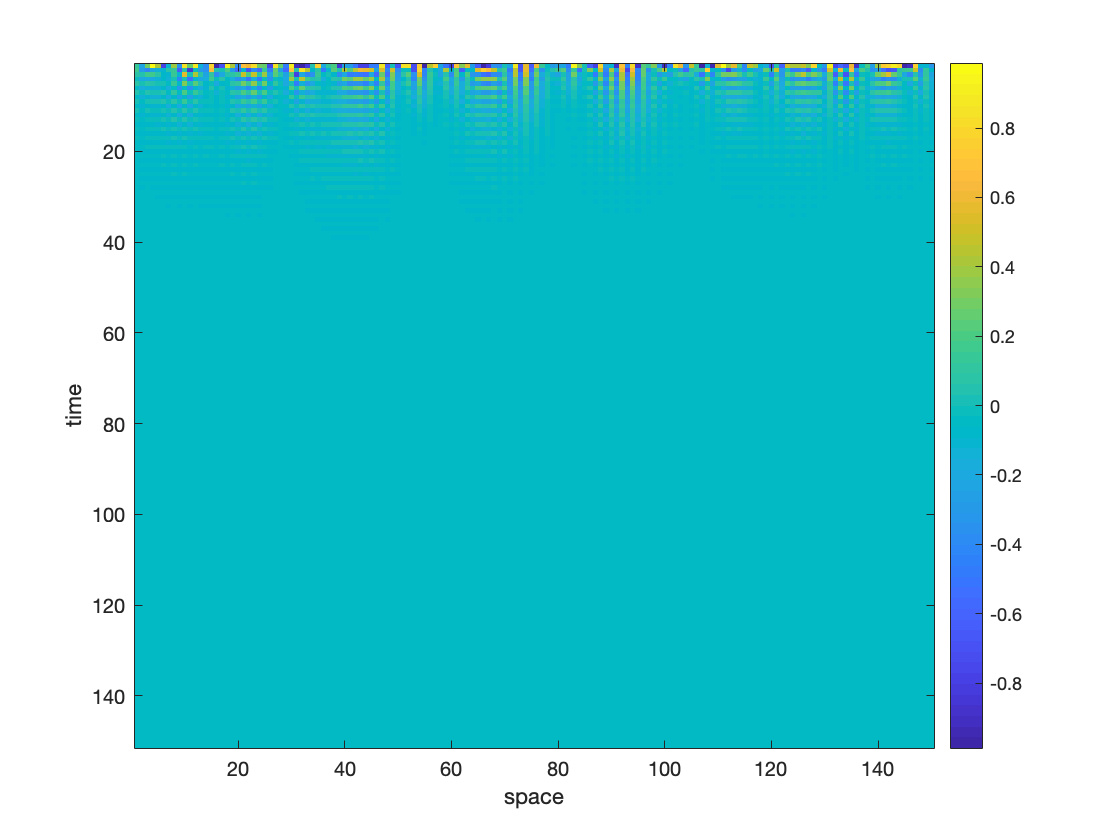}
}
\caption{Spatio-temporal patterns of CMLs with local $T_{2, 0}$ for increasing coupling strength. Colour encodes the value of the dynamical variable $x_n^{(i)}$. The space$\times$time size is $x\times y = 150\times 150$; on each lattice site $i$ the initial value is randomly chosen as $x_0^{(i)} \in \text{Uni}(-1, 1)$, with periodic boundary conditions.} 
\label{spatio-temporal-T2}
\end{figure}

\begin{figure}[H]
\centering 
\subfloat[\footnotesize Type $3A$: $c = 0$]{
\includegraphics[width = 0.32\textwidth]{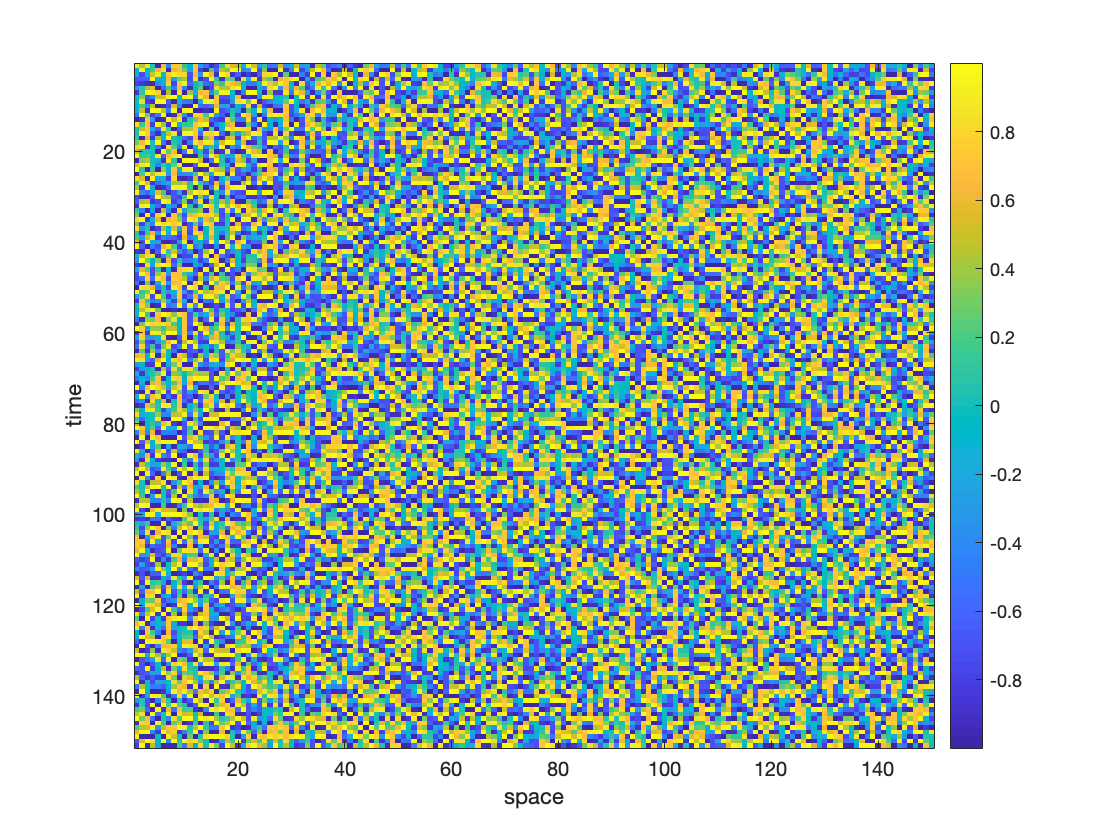}
}
\subfloat[\footnotesize Type $3A$: $c = 0.45$]{
\includegraphics[width = 0.32\textwidth]{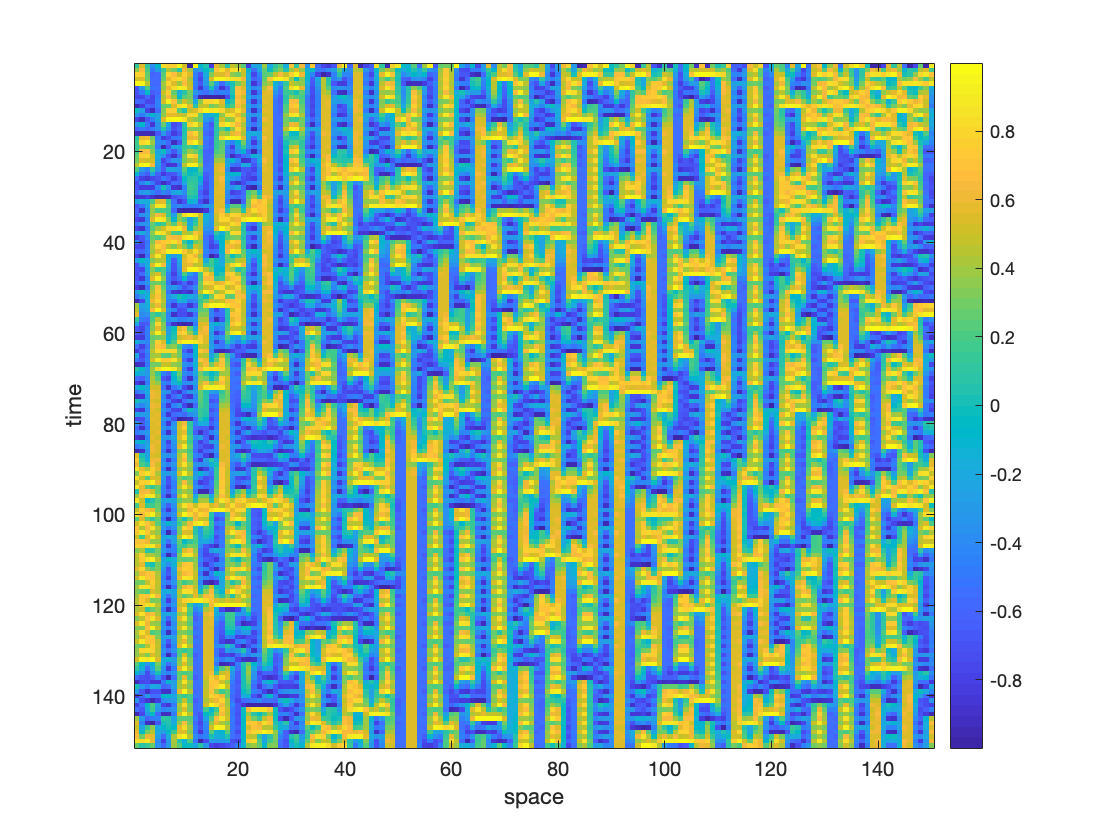}
\label{interm}
}
\subfloat[\footnotesize Type $3A$: $c = 0.95$]{
\includegraphics[width = 0.32\textwidth]{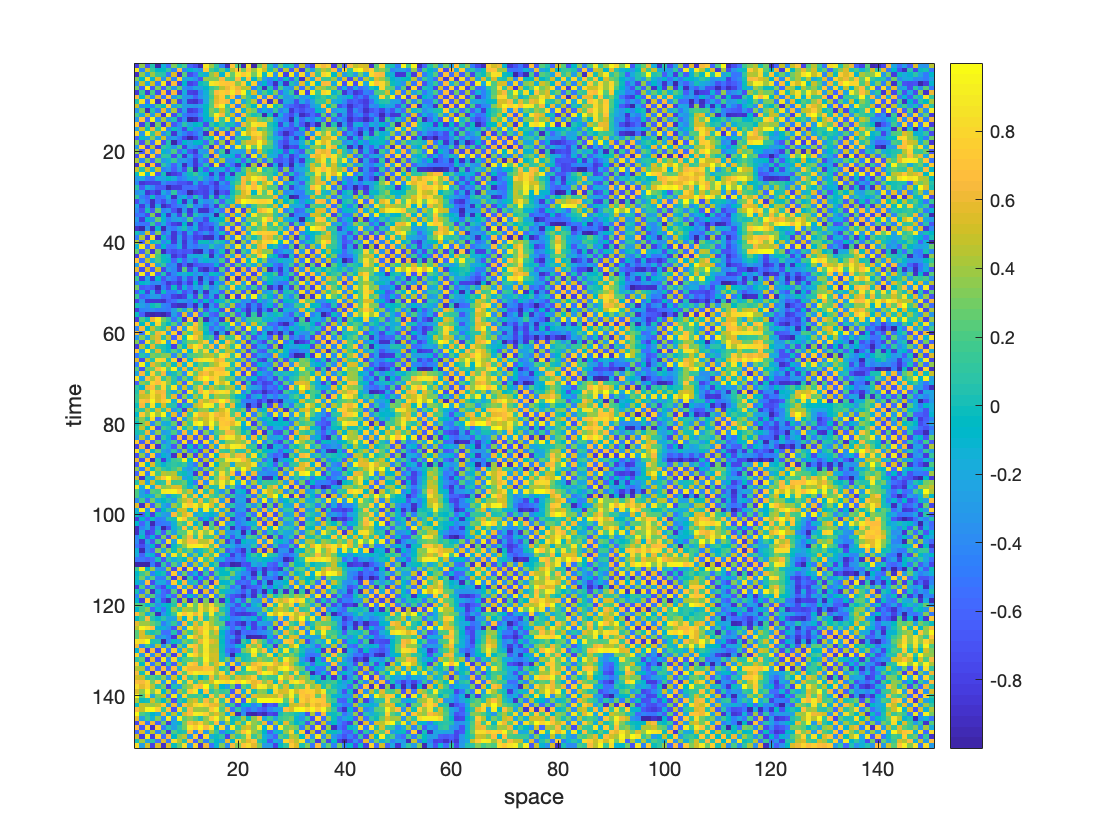}
\label{more1}
}

\subfloat[\footnotesize Type $3A^-$: $c = 0.11$]{
\includegraphics[width = 0.32\textwidth]{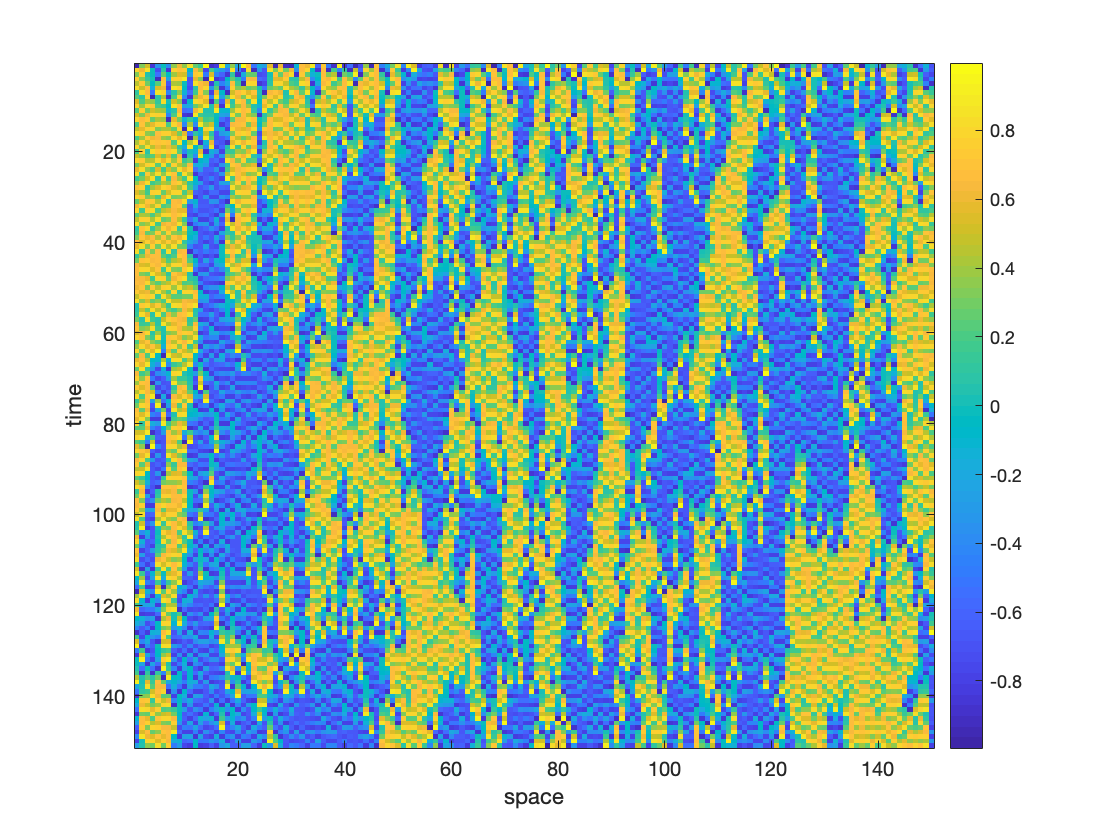}
\label{more2}
}
\subfloat[\footnotesize Type $3A^-$: $c = 0.36$]{
\includegraphics[width = 0.32\textwidth]{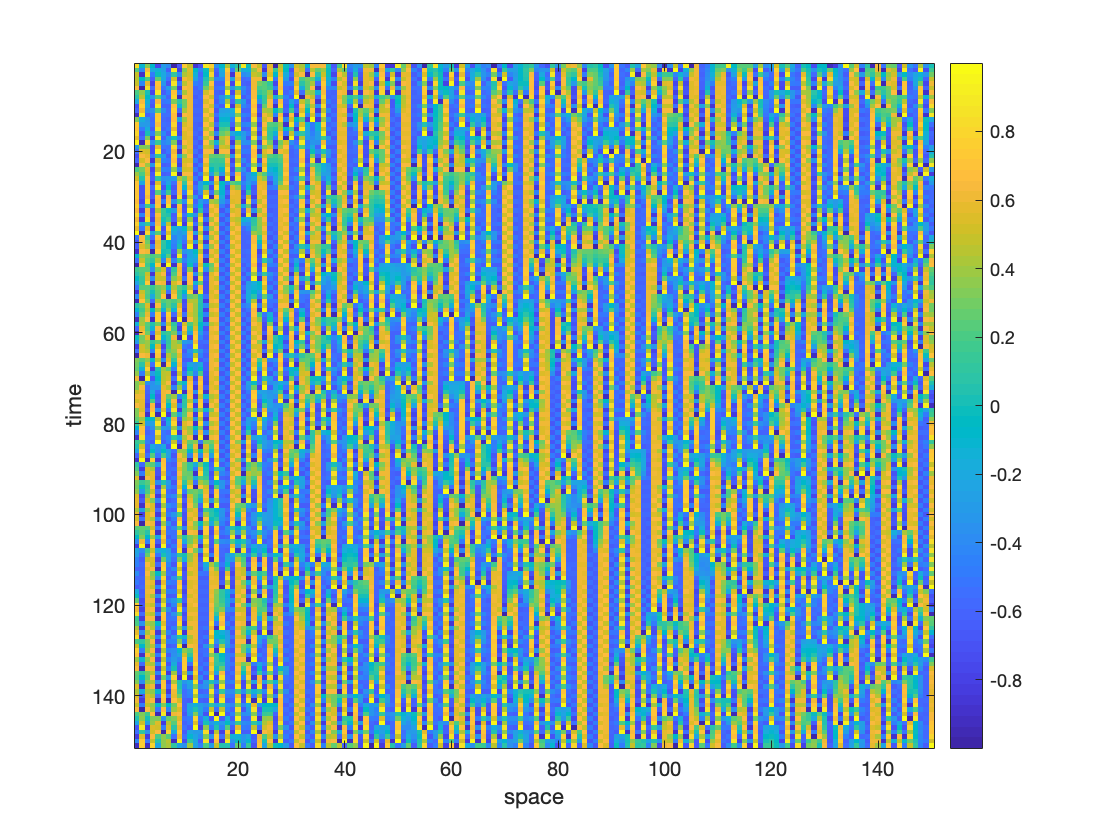}
}
\subfloat[\footnotesize Type $3A^-$: $c = 0.91$]{
\includegraphics[width = 0.32\textwidth]{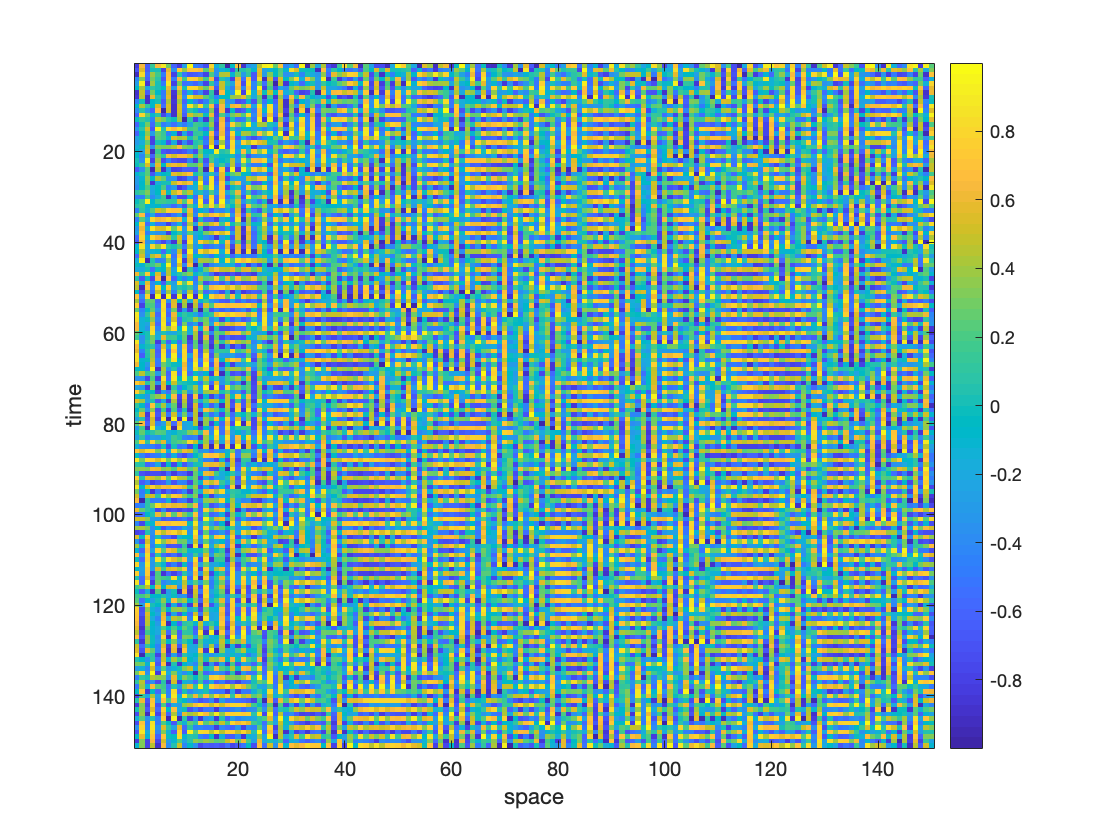}
\label{chaosBro}
}

\subfloat[\footnotesize Type $3B$: $c = 0.19$]{
\includegraphics[width = 0.32\textwidth]{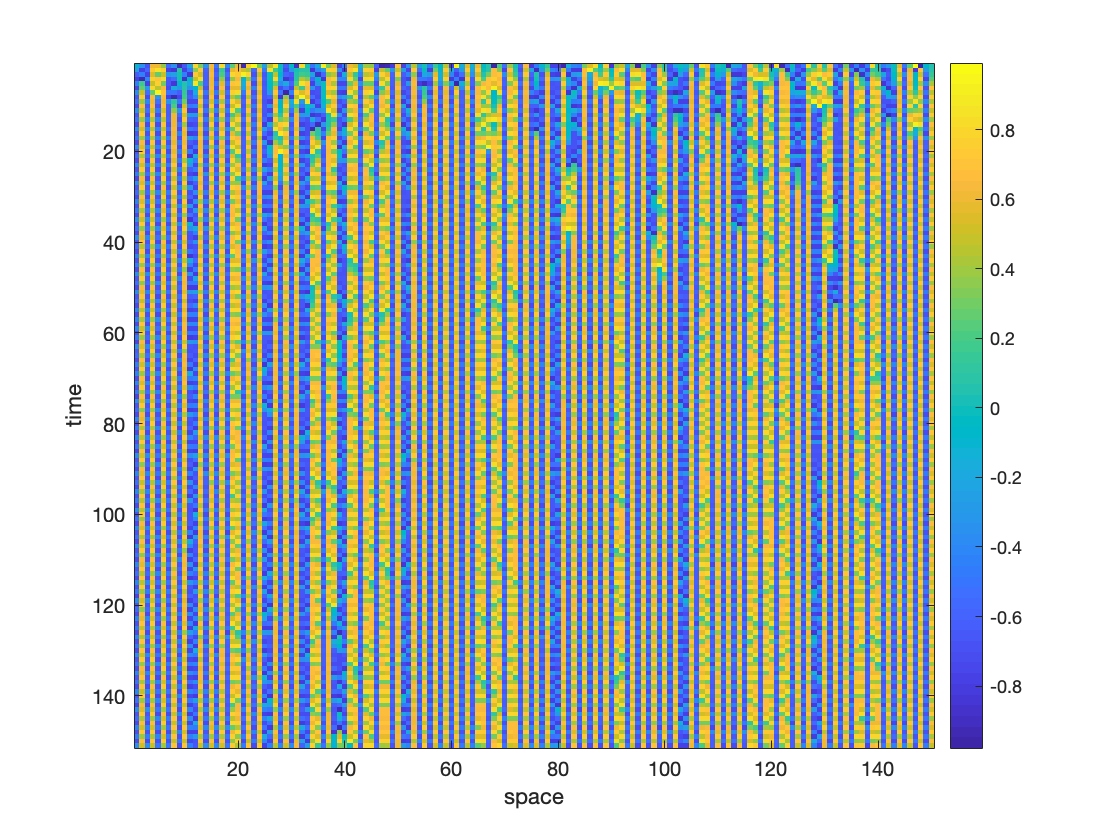}
}
\subfloat[\footnotesize Type $3B$: $c = 0.52$]{
\includegraphics[width = 0.32\textwidth]{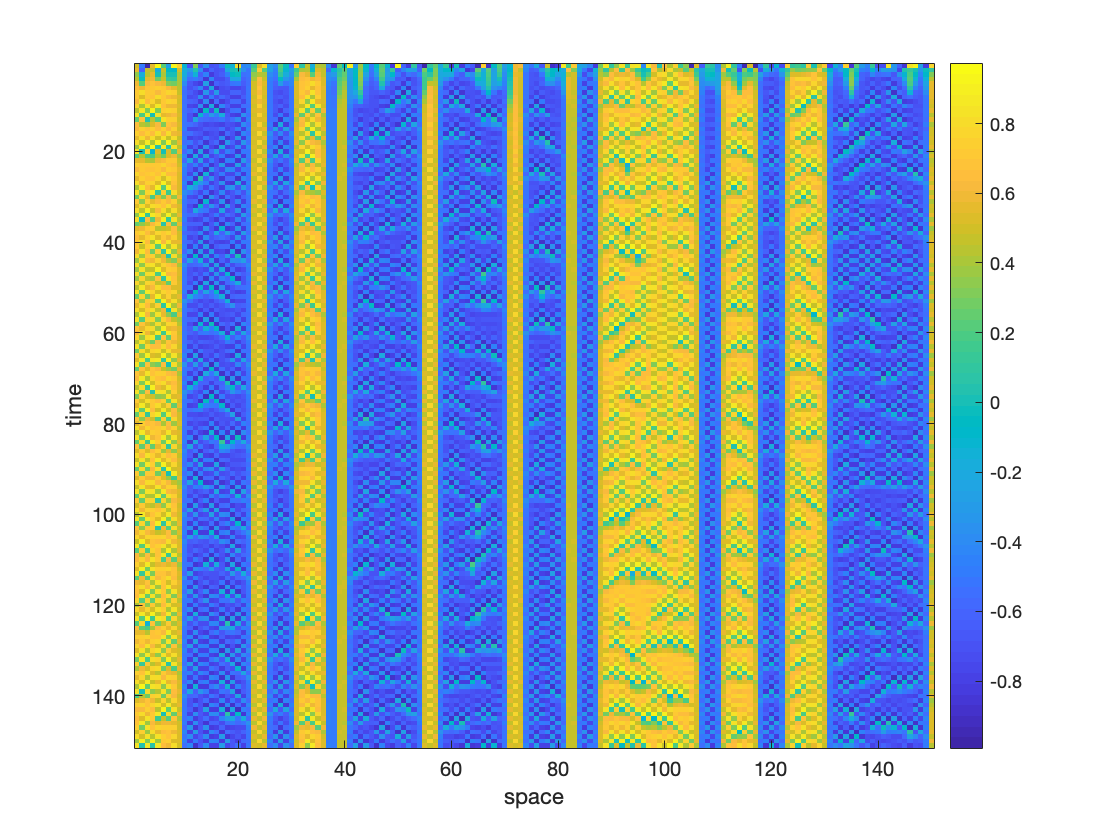}
}
\subfloat[\footnotesize Type $3B$: $c = 0.92$]{
\includegraphics[width = 0.32\textwidth]{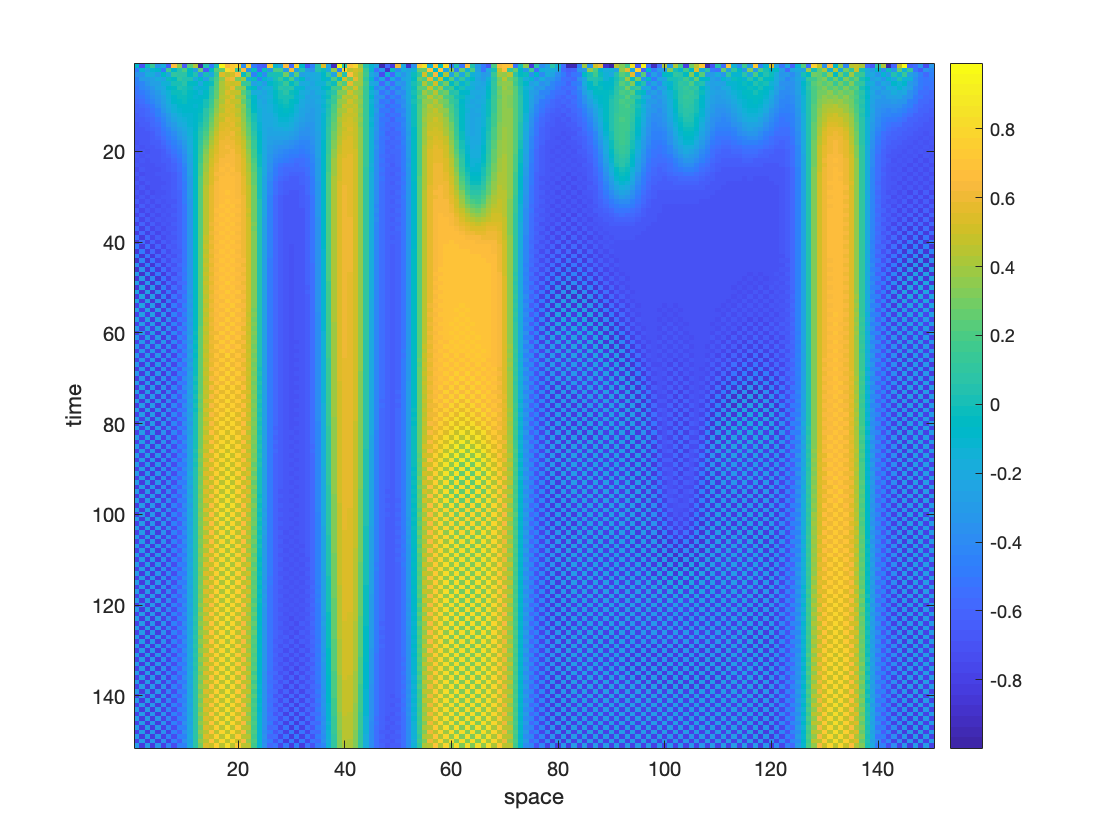}
}

\subfloat[\footnotesize Type $3B^-$: $c = 0.24$]{
\includegraphics[width = 0.32\textwidth]{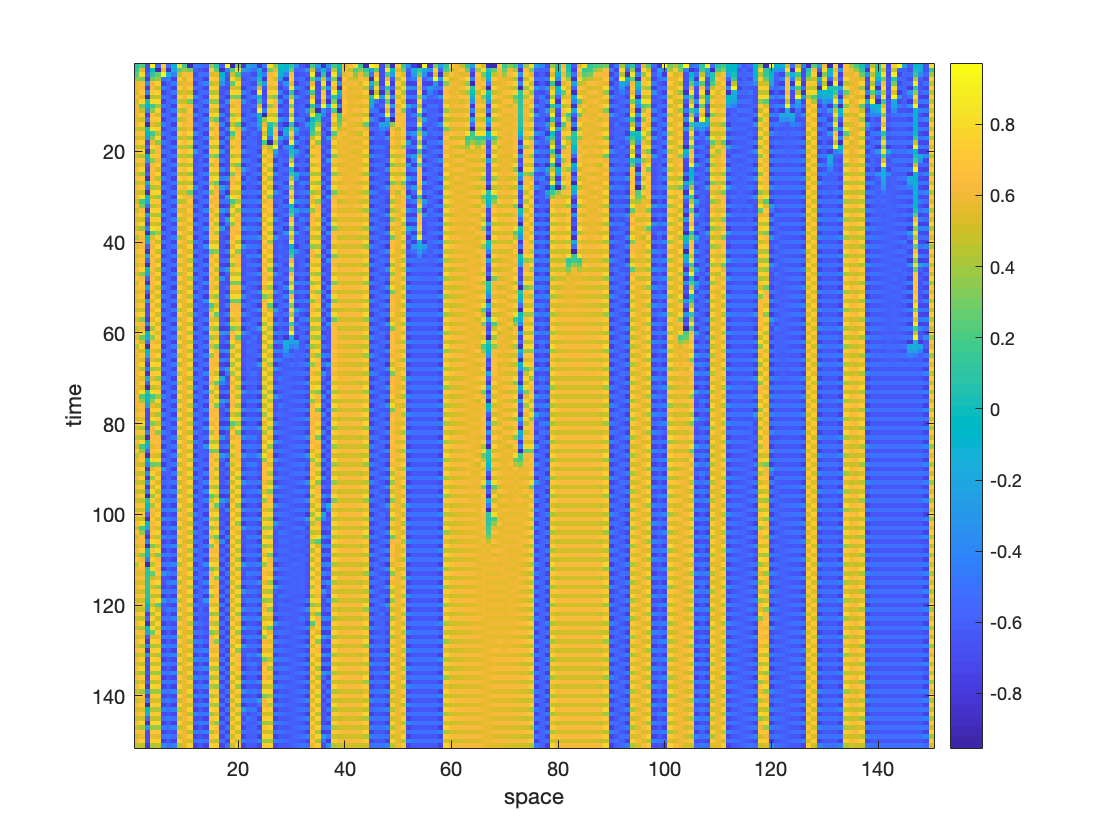}
}
\subfloat[\footnotesize Type $3B^-$: $c = 0.49$]{
\includegraphics[width = 0.32\textwidth]{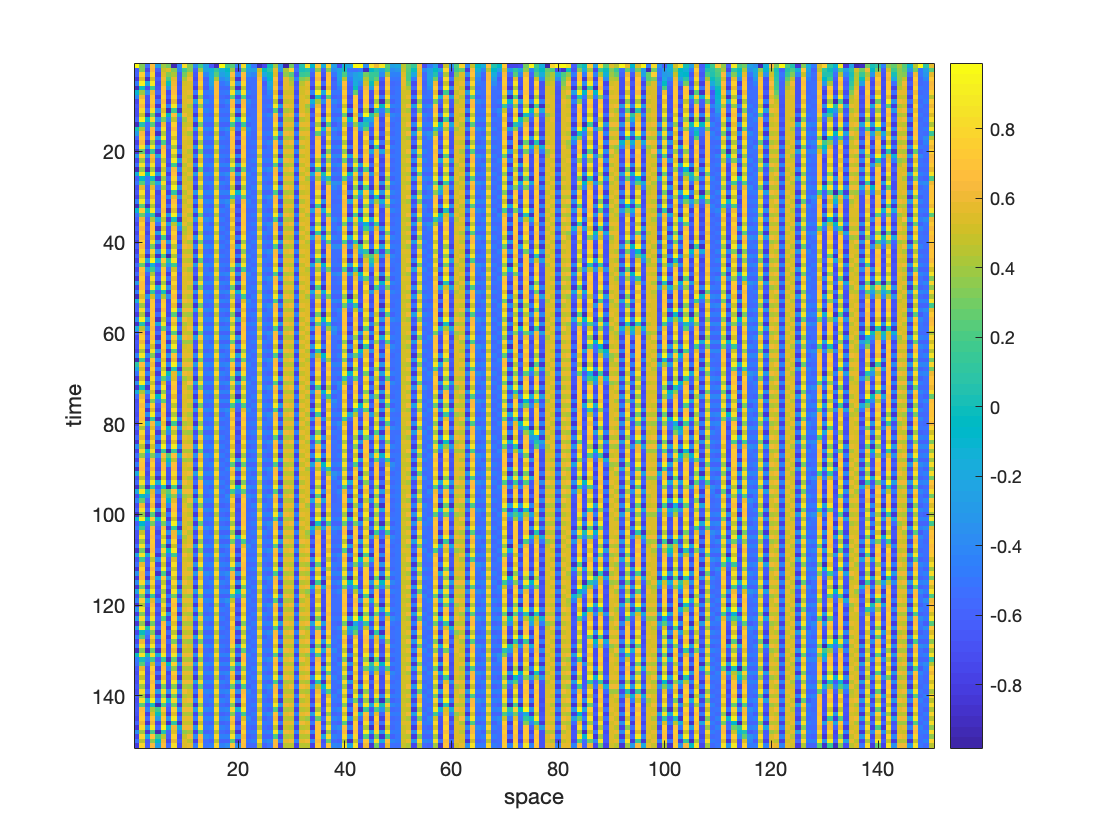}
\label{fro-chaos}
}
\subfloat[\footnotesize Type $3B^-$: $c = 0.93$]{
\includegraphics[width = 0.32\textwidth]{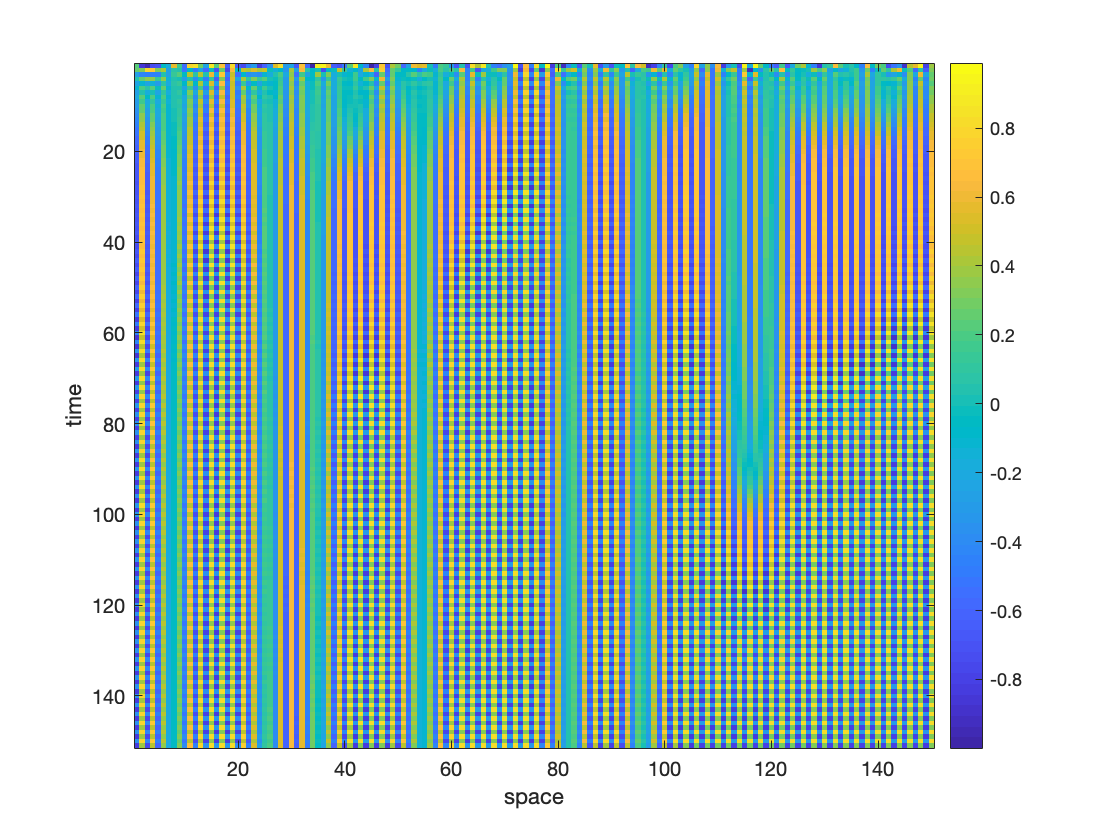}
}
\caption{Spatio-temporal patterns of CMLs with local $T_{3, -\frac{\pi}{2}}$; same iteration parameters as in Fig.\ref{spatio-temporal-T2}.} 
\label{spatio-temporal-T3sine}
\end{figure}

A variety of interesting patterns is generated, such as turbulent-like behaviour (Fig.\ref{turb}), intermittency (Fig.\ref{interm}), and frozen chaos (Fig.\ref{fro-chaos}); see e.g.\cite{kaneko} for a characterisation of different phenomena. 
Some of the patterns remind us of biological structures: Figs.\ref{more1}, \ref{more2} resemble some self-organised patches in nature or ecosystems, such as kelps and duckweeds.

\subsection{Spatial and temporal correlation functions} 
Our main interest in this paper is to quantify the correlation structure. For CMLs, two different types of correlations occur, namely in the spatial direction and in the temporal direction. Let us define the $SNNC$ (\textit{spatial nearest-neighbour correlation}) and the $TNNC$ (\textit{temporal nearest-neighbour correlation}) as the following averages:
\begin{equation}
\begin{split}
SNNC &= \lim_{K \rightarrow \infty} \frac{1}{KJ}\sum_{n = 1}^K\sum_{j = 1}^J x_n^{(j)}x_n^{(j+1)},\\
TNNC &= \lim_{K \rightarrow \infty} \frac{1}{JK}\sum_{j = 1}^J\sum_{n = 1}^K x_n^{(j)}x_{n+1}^{(j)}.
\end{split} 
\label{def-snnc-tnnc}%
\end{equation}
For finite $K$ and $J$ the order of the two sums is interchangeable. In practice, we are interested in very large values of $K$ and $J$, and in particular in the limit $K \to \infty$, corresponding to the long-term iteration limit, where the system (if ergodic) may again approach an invariant density, which for the CML is a function of $J$ different variables. In the following, we show some numerical results for the above observables, $SNNC$ and $TNNC$, as a function of the parameters $c$ and $a$.

(i) CML of Type $2A$: 
\begin{figure}[H]
\centering 
\subfloat[\footnotesize $SNNC$ (blue) and $TNNC$ (red) at $a = 0$]{
\includegraphics[width = 0.48\textwidth]{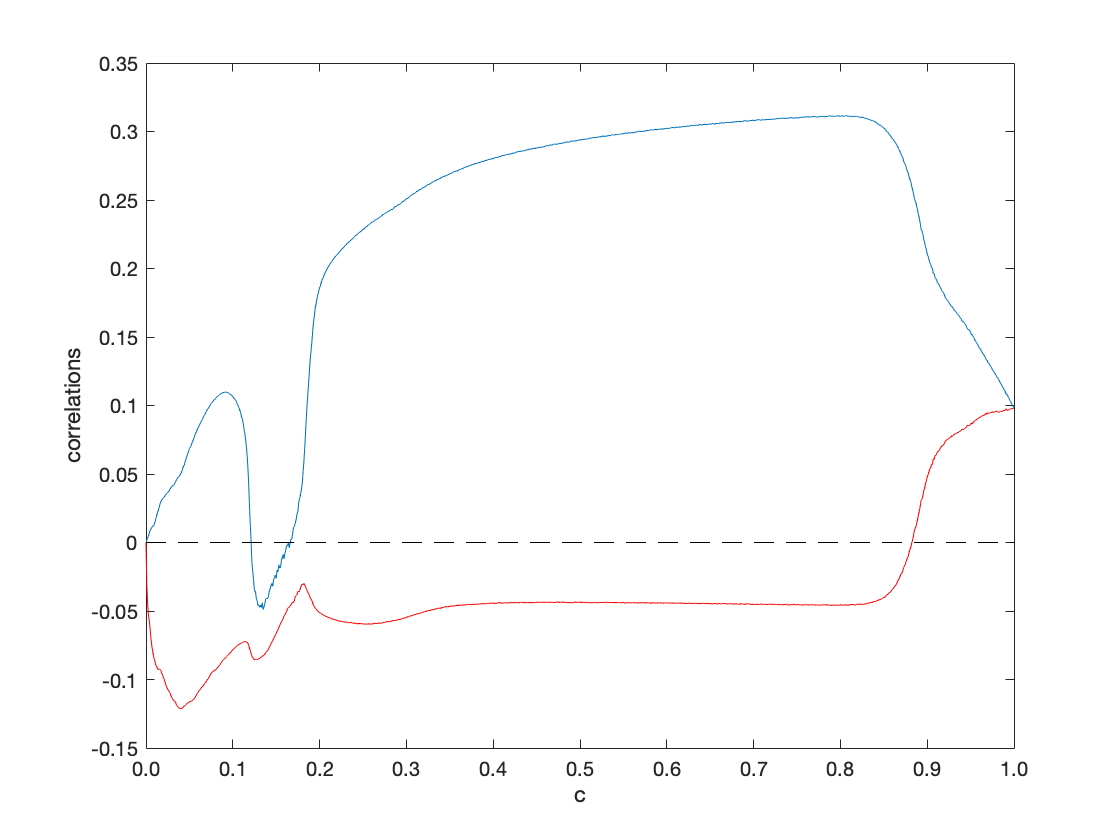}
}
\subfloat[\footnotesize $SNNC$ (z-height) and $TNNC$ (coloured)]{
\includegraphics[width = 0.48\textwidth]{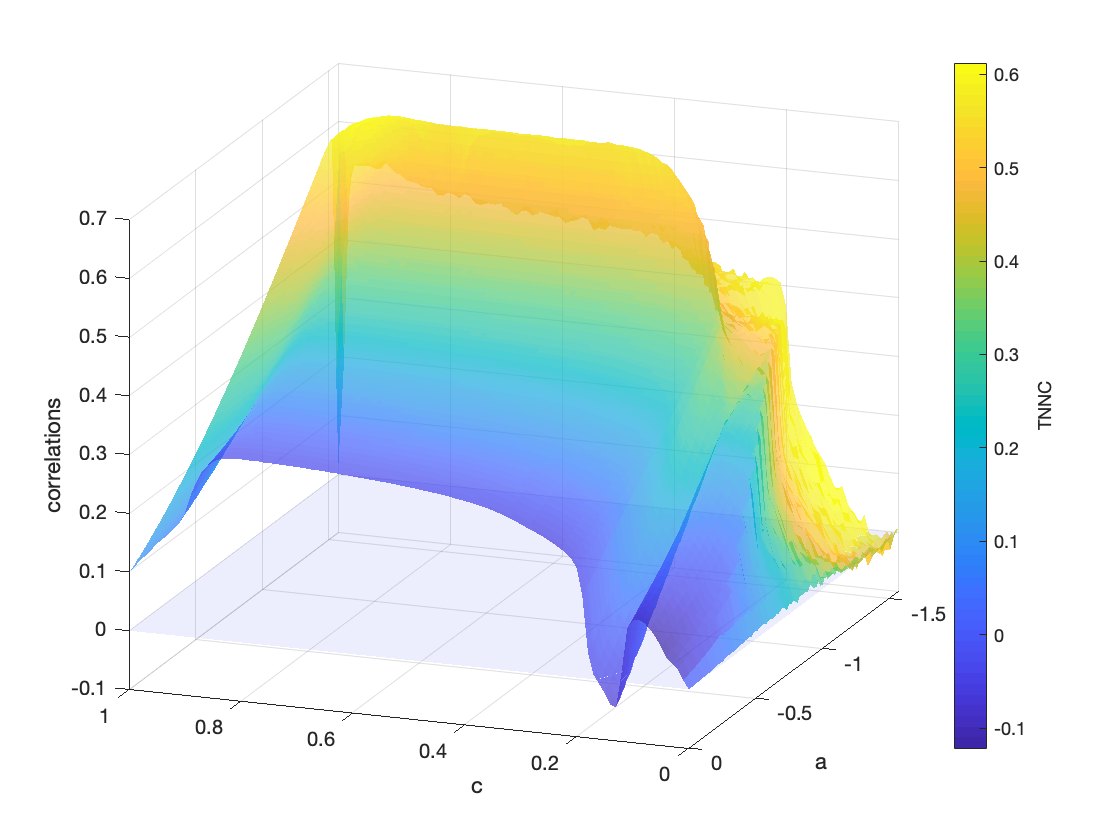}
}
\caption{$SNNC$ and $TNNC$ for Type $2A$ CML: $c \in [0, 1]$, $a \in \left[ -\frac{\pi}{2}, 0\right]$. The space$\times$time size is $J\times K = 5,000 \times 1,100(-100)$; initial points randomly chosen from a uniform distribution Uni$(-1, 1)$.} 
\label{2A}
\end{figure}

Notice that in the correlation surface plot, there is a rapid folding around $c \in [0.1, 0.2]$, creating two distinct zeros of $SNNC$ depending on $a$. Fluctuations as $c \rightarrow 0$ with $a \rightarrow -\frac{\pi}{2}$ indicate that the system loses its ergodic property and that the averaged correlation fails to stabilise, which could be induced by non-mixing of the underlying system. These fluctuation regions, depending on both the coupling $c$ and the shift parameter $a$, are very intricate, and occur
for other coupling forms as well.

For $a = 0$ we observe special coupling parameter values (such as $c \sim 0.12$) where $SNNC = 0$, see \cite{CBblue} for a physical interpretation in a quantum field theoretical setting. There are also special values (such as $c \sim 0.88$) where $TNNC=0$. This means that although there is non-trivial spatial coupling, some features of the uncoupled (most random-looking) local Chebyshev dynamics are restored for these special coupling constants.

(ii) CML of Type $2A^-$:
\begin{figure}[H]
\centering 
\subfloat[\footnotesize $SNNC$ (blue) and $TNNC$ (red) at $a = 0$]{
\includegraphics[width = 0.48\textwidth]{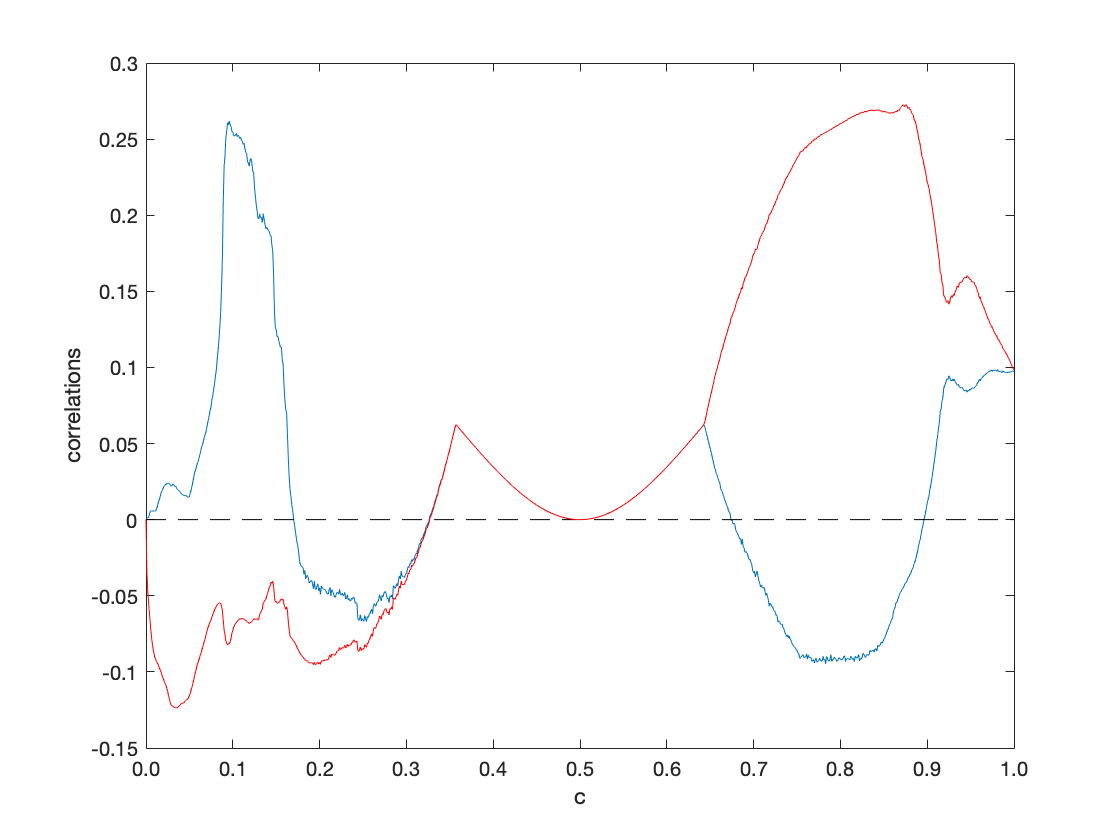}
}
\subfloat[\footnotesize $SNNC$ (z-height) and $TNNC$ (coloured)]{
\includegraphics[width = 0.48\textwidth]{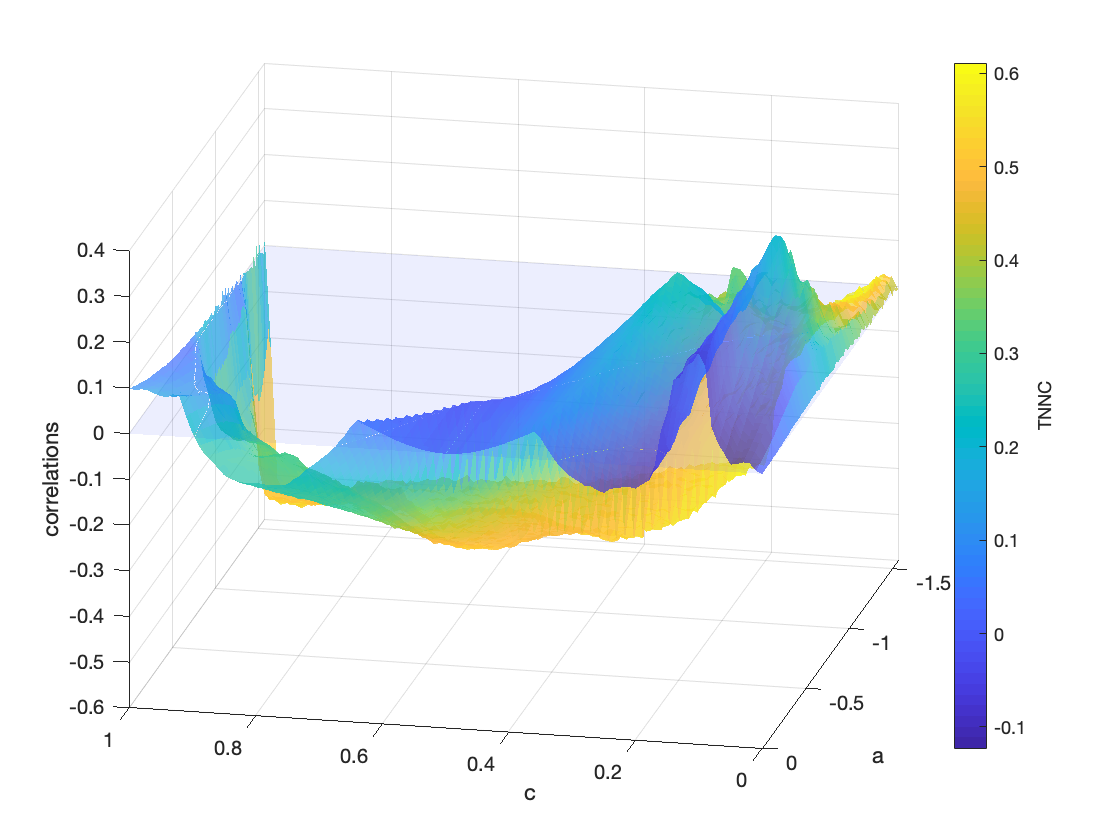}
}
\caption{Same as in Fig.\ref{2A} except for Type $2A^-$.}
\end{figure}

\begin{figure}[H]
\centering 
\subfloat[\footnotesize $a = -\pi/16$]{
\includegraphics[width = 0.32\textwidth]{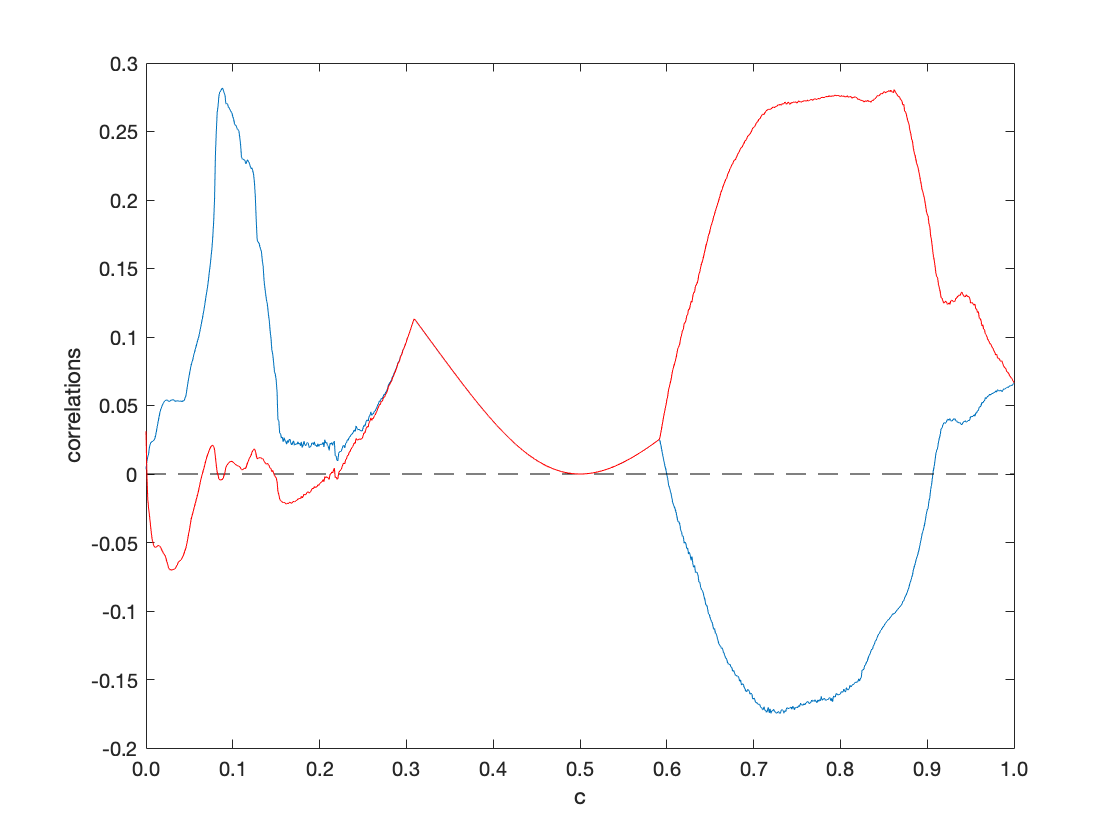}
}
\subfloat[\footnotesize $a = -\pi/8$]{
\includegraphics[width = 0.32\textwidth]{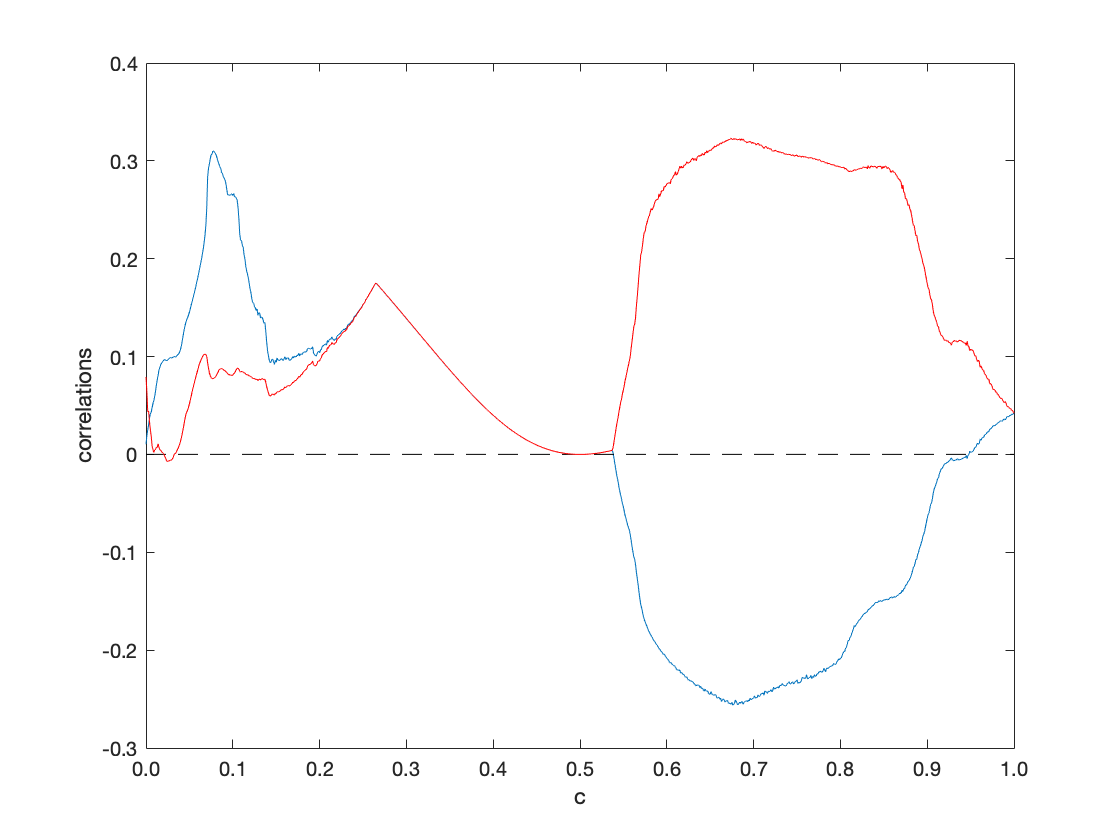}
}
\subfloat[\footnotesize $a = -3\pi/16$]{
\includegraphics[width = 0.32\textwidth]{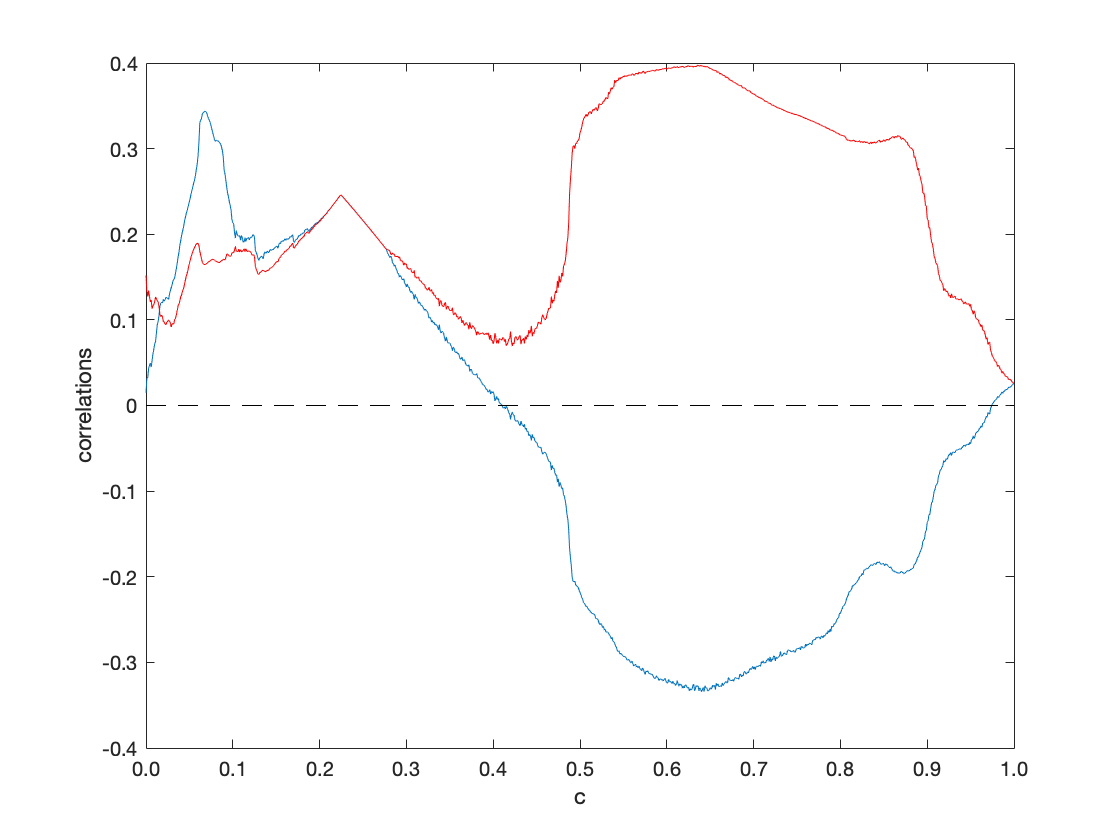}
}
\caption{$SNNC$ (blue) and $TNNC$ (red) for Type $2A^-$ at $a = -\pi/16$, $-\pi/8$ and $-3\pi/16$, $c \in [0, 1]$. Other parameters are the same as in Fig.\ref{2A}.}
\label{parabolic-zeros}
\end{figure}

Here in most parts the $SNNC$ is negative while $TNNC$ keeps growing as $a \rightarrow -\frac{\pi}{2}$. The  parabola-like curve starting at $(c, a) = (0.5, 0)$ is due
to a stable synchronised fixed point of the CML, as a simple stability analysis shows, and it disappears when $a > 0.45$, see Fig.\ref{parabolic-zeros} for more detail.

(iii) CML of Type $2B$: 
\begin{figure}[H]
\centering 
\subfloat[\footnotesize $SNNC$ (blue) and $TNNC$ (red) at $a = 0$]{
\includegraphics[width = 0.48\textwidth]{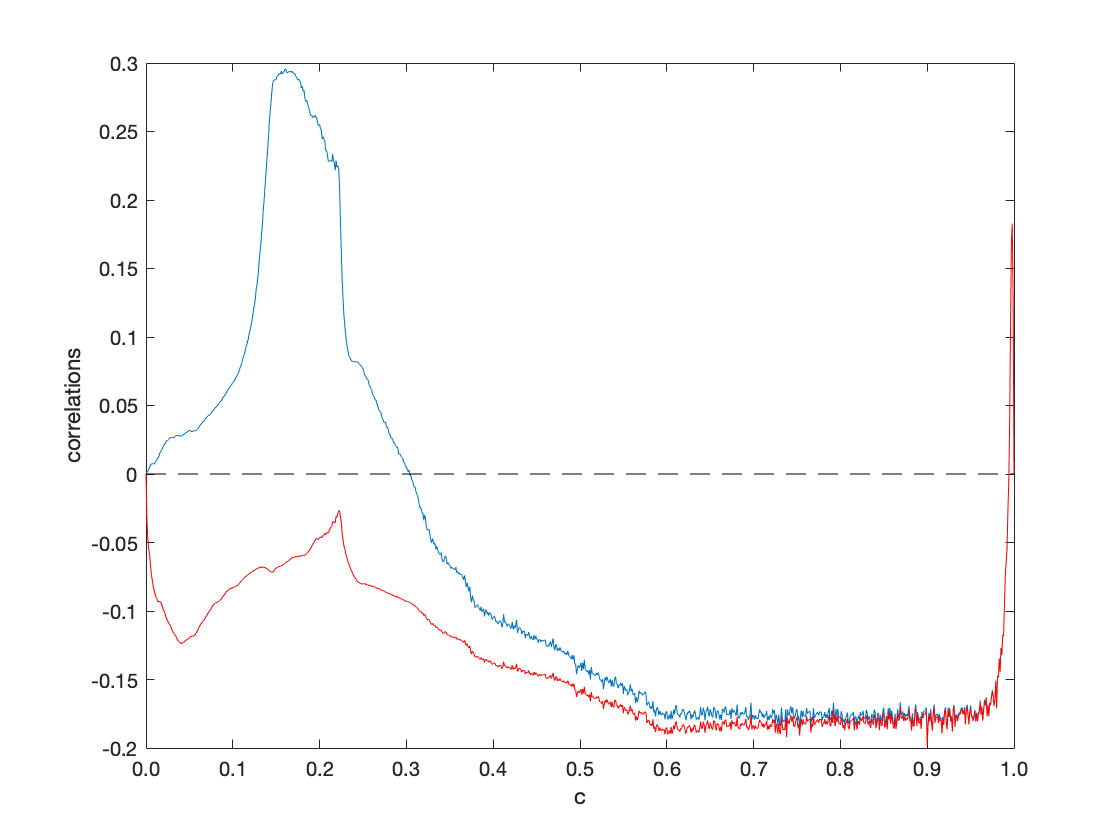}
}
\subfloat[\footnotesize $SNNC$ (z-height) and $TNNC$ (coloured)]{
\includegraphics[width = 0.48\textwidth]{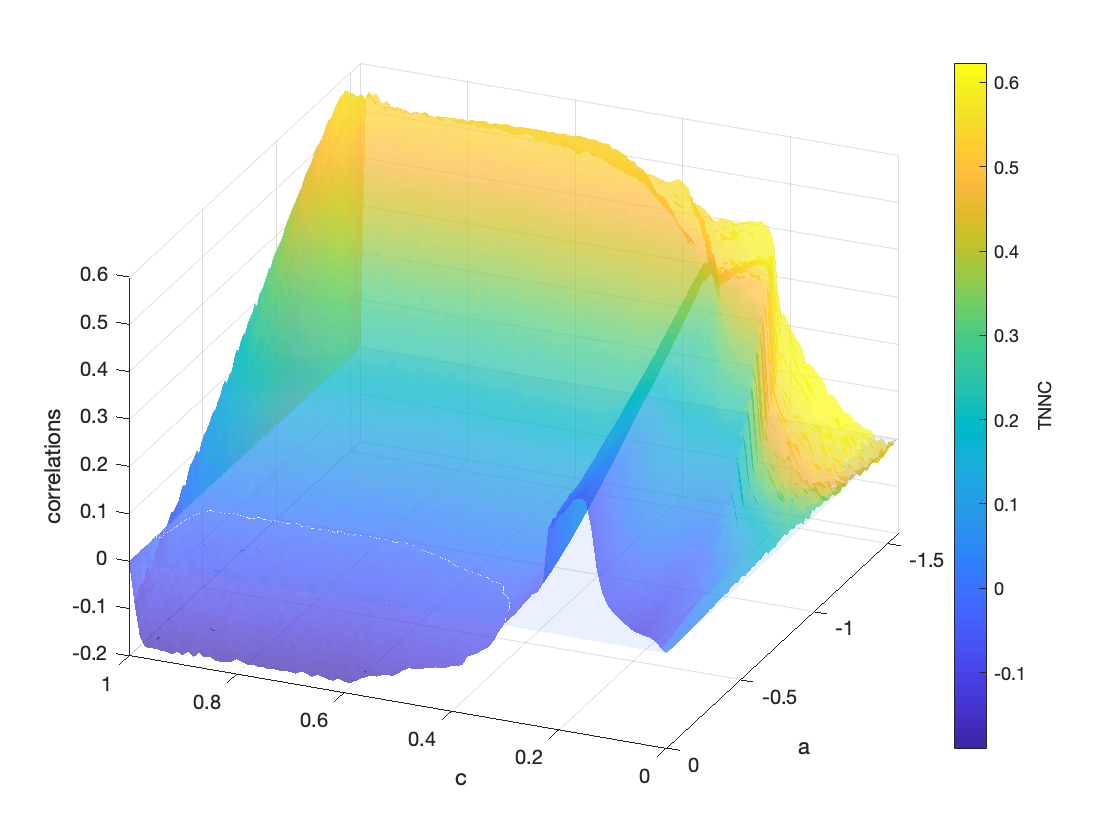}
}
\caption{Same as in Fig.\ref{2A} except for Type $2B$.}
\end{figure}

For small $c$ the correlation surface is similar as for Type $2A$: there also exists a rapid folding near $c=0.15$. In contrast to types $2A^-$ and $2B^-$, there exists no stable synchronised fixed point state for larger $c$.

(iv) CML of Type $2B^-$:
\begin{figure}[H]
\centering 
\subfloat[\footnotesize $SNNC$ (blue) and $TNNC$ (red) at $a = 0$]{
\includegraphics[width = 0.48\textwidth]{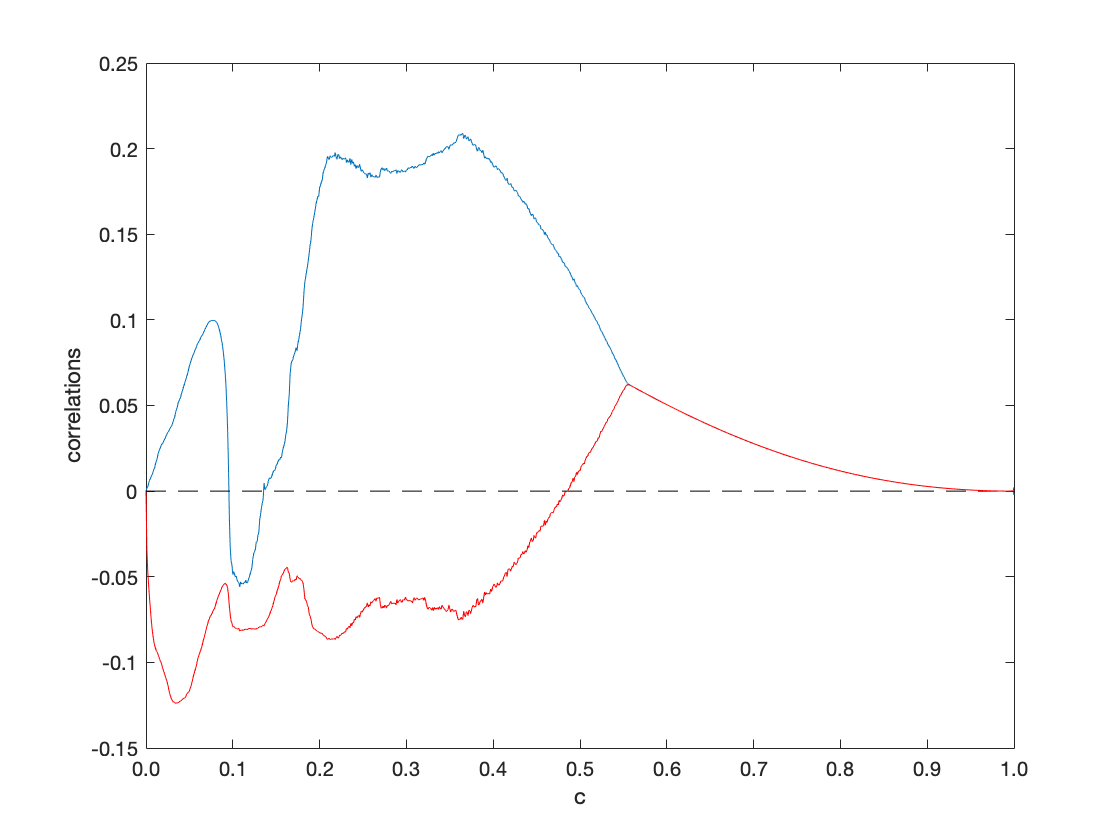}
}
\subfloat[\footnotesize $SNNC$ (z-height) and $TNNC$ (coloured)]{
\includegraphics[width = 0.48\textwidth]{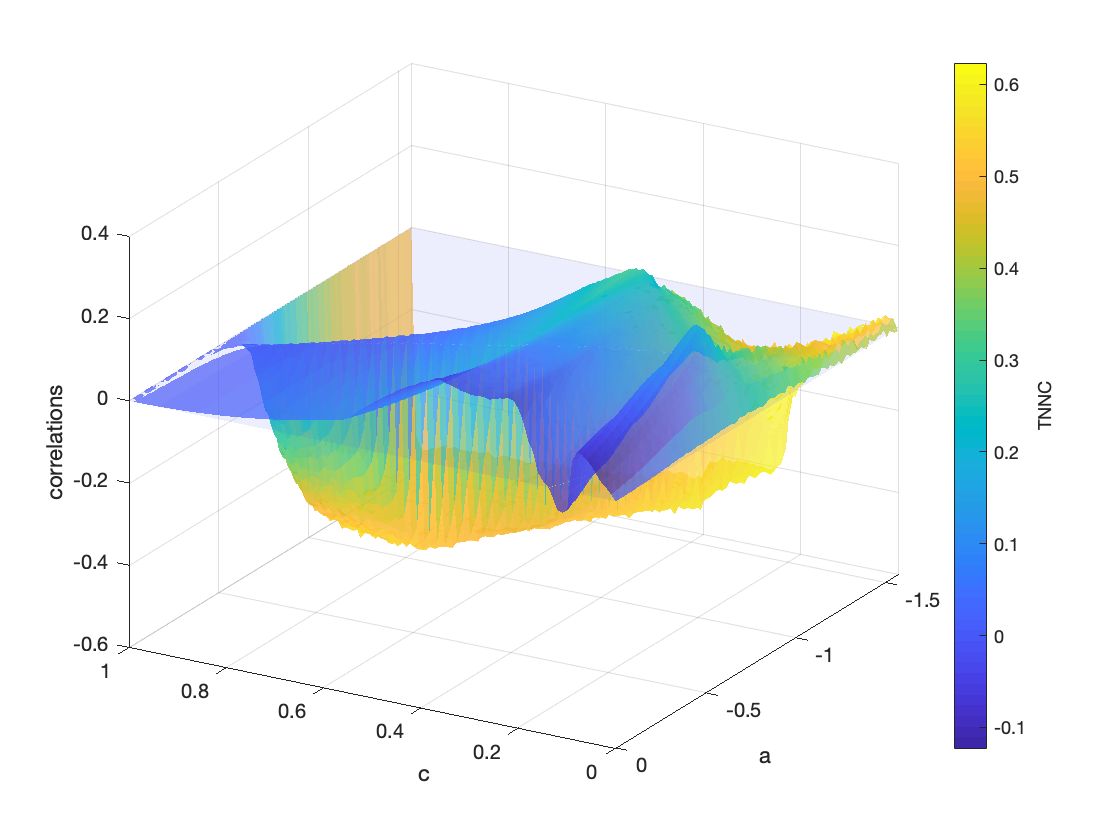}
}
\caption{Same as in Fig.\ref{2A} except for Type $2B^-$.}
\end{figure}

Yet another type of coupling form relevant in quantum field theoretical applications. For $c \to 1$ and small $a$ we observe a stable synchronised fixed point. In general, as outlined in more detail in \cite{CBblue}, physical interaction states correspond to zeros of $SNNC$.

Finally, let us also consider coupled Chebyshev maps with odd $N$, which have a different symmetry behaviour: for odd $N$ coupled Chebyshev systems, the average of the iterates is always zero, whereas for even $N$ it is not.

The coloured plots in Fig.\ref{3A3B} below show the $SNNC$ (z-height) and $TNNC$ (encoded in colour) for types $3A$ and $3B$ CMLs (i.e., the local map is $T = T_{3, 0}$). For Type $3A$ there is a dip in the spatial correlation ($SNNC$) when $a$ is close to $0$ with weak coupling $c$; while for Type $3B$ a jump in $SNNC$ occurs around the same parameter region, and when $c > 0.5$ again fluctuations caused by non-mixing are prevalent. The temporal correlation ($TNNC$) is increasing when $|a|$ becomes larger.

Overall, we notice that spatial coupling destroys the simple, distinguished correlation properties of Chebyshev maps, as visible already for the nearest neighbour-correlation in space and time, which was identical to zero for $a = 0$ with $c=0$. Still, some distinguished non-trivial parameter values exist in the $(c,a)$-plane which generate uncorrelated nearest-neighbour behaviour. As an example, subplots in Fig.\ref{3A3B} indicate the special curves in the $(c, a)$-plane where $SNNC=0$ (blue) and $TNNC=0$ (red) for both $3A$ and $3B$ coupling forms.

\begin{figure}[H]
\centering 
\subfloat[\footnotesize Type $3A$]{
\includegraphics[width = 0.48\textwidth]{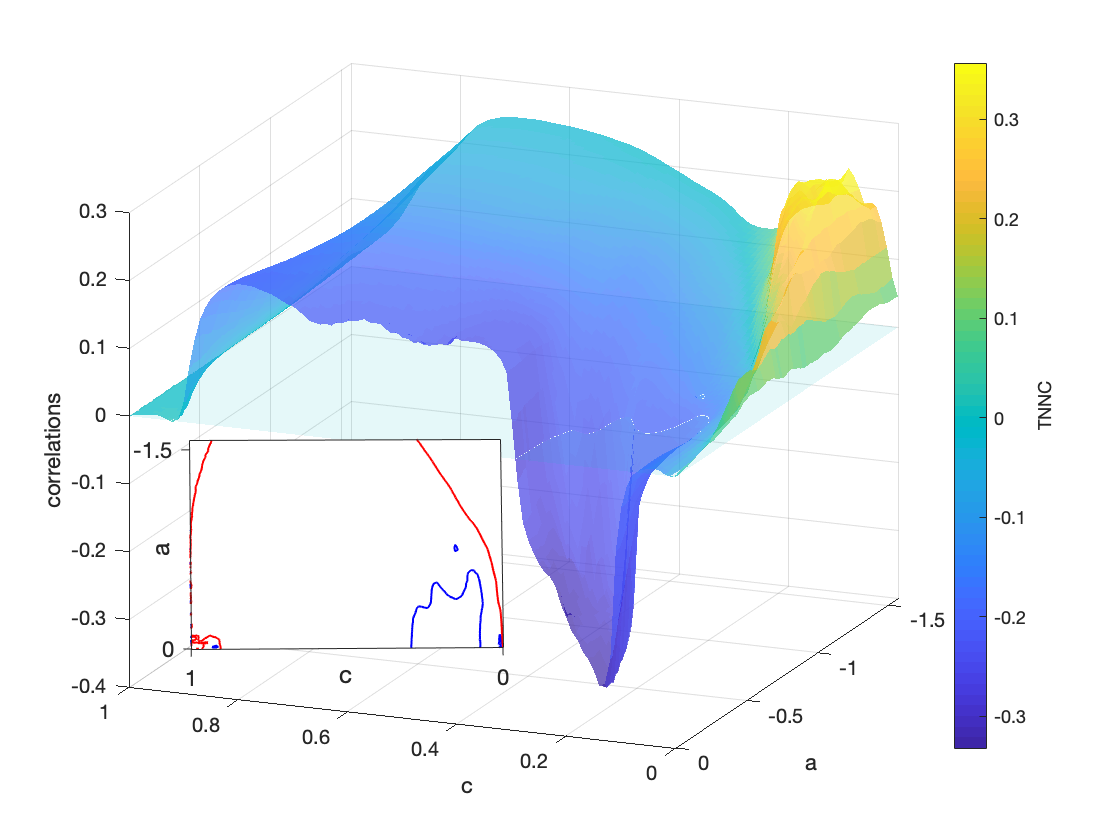}
}
\subfloat[\footnotesize Type $3B$]{
\includegraphics[width = 0.48\textwidth]{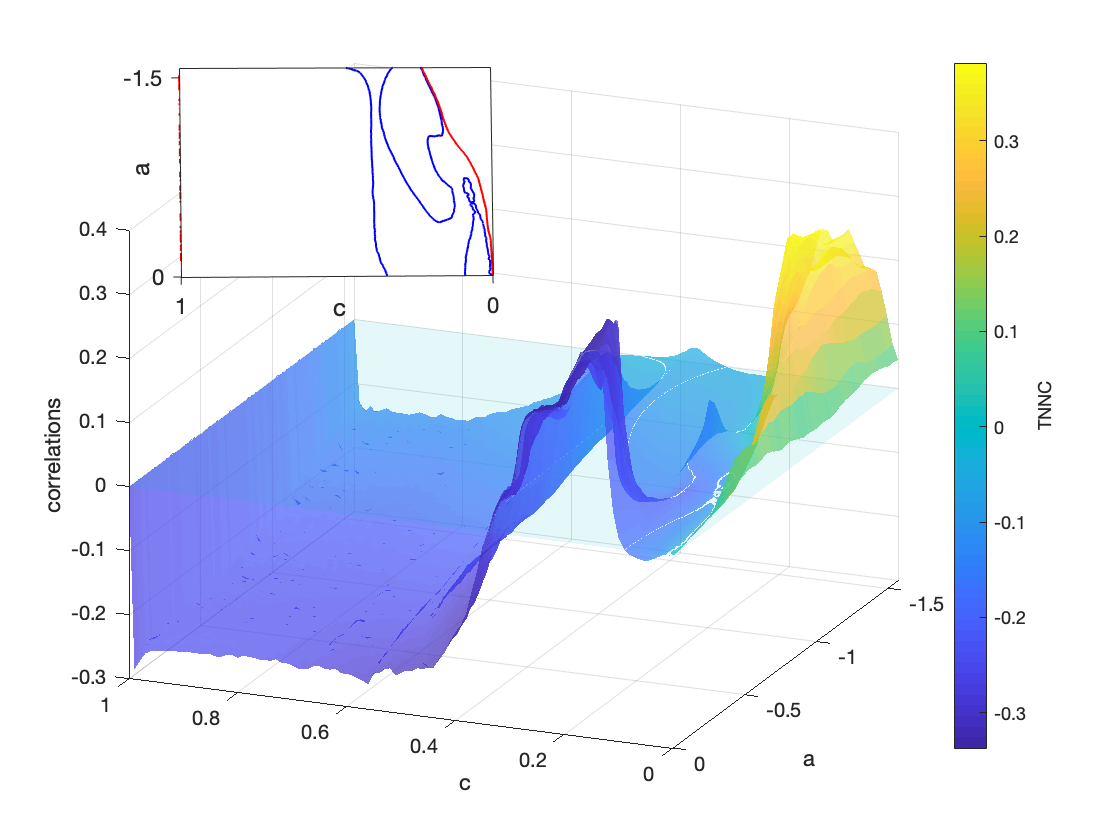}
}
\caption{$SNNC$ (z-height) and $TNNC$ (coloured) for types $3A$ and $3B$: $c \in [0, 1]$, $a \in \left[ -\frac{\pi}{2}, 0\right]$, with subplots showing zeros of $SNNC$ (blue) and $TNNC$ (red). Other parameters are the same as in Fig.\ref{2A}.}
\label{3A3B}
\end{figure}

\section{Conclusion}
In this paper we have shown that Chebyshev maps, among all maps conjugated to a shift dynamics of $N$ symbols, are very distinguished dynamical systems. They generate least higher-order correlations, and can thus be regarded as the most random-looking systems, among all possible dynamical systems, assuming a deterministic evolution dynamics as given given by a smooth differentiable one-dimensional map conjugated to an $N$-ary shift. We generalised the concept and were able to determine invariant densities for shifted Chebyshev maps as well, which are a modification of the ordinary Chebyshev map described by a translation parameter $a$. For particular values of the parameters $N$ and $a$ we were able to determine the entire (discrete) set of eigenfunctions of the Perron-Frobenius operator.

There are several aspects of evidence that the (ordinary, $a = 0$) Chebyshev maps are distinguished among other deterministic chaotic systems under consideration. First of all, for iterates of an uncoupled Chebyshev map, the two-point correlation decays very fast --- described by a Kronecker delta function --- indicating that it vanishes immediately instead of exponentially such as for an $N$-ary shift map with subtracted mean; zeros of higher-order correlations can be determined by diophantine equations with the simplest possible spin configurations, namely up-down spin space; other dynamical systems correspond to a more complicated integer-spin space that is embedded in the generalised diophantine equations we considered in this paper.

It is also remarkable that for Chebyshev maps an orthogonal set of eigenfunctions of the Perron-Frobenius operator exists, and that the eigenvalues are real, in spite of the fact that the operator itself is not Hermitean. We were able to find the eigenfunctions also for the generalised shifted Chebyshev maps, in the case of even $N$. The observation here is that a conjugating function $h_1$ can be found which relates the dynamics of shifted Chebyshev maps to that of ordinary Chebyshev maps, in particular, $h_1$ is simply a negative (ordinary) Chebyshev polynomial whose order depends on the translation parameter $a$.

When a spatial coupling $c \not= 0$ is introduced, the simple Bernoulli shift properties of shifted Chebyshev maps are destroyed in the corresponding coupled map lattice. However, numerically still some distinguished coupling constants can be found where the spatial nearest-neighbour correlation of the coupled map lattice vanishes. This means that in spite of the non-zero coupling, still a random-like correlation state can be achieved. These types of states have physical meaning in chaotically quantized field theories,  fixing allowed types of coupling constants, and making contact to observed coupling constants in the standard model of elementary particle physics \cite{CBblue, chaotic-string}. Our consideration in this paper has shown that indeed Chebyshev maps are the most distinguished dynamical systems, in the sense that they allow to construct a most random-like looking chaotic field theory, which at the microscopic level is purely deterministic.

\newpage
\appendix

\scriptsize
\section{Proof of topological conjugation between shifted Chebyshev maps and piecewise-linear maps}\label{A1}
Let us denote $X = [-1, 1]$, $Y = [0, 1]$. The shifted Chebyshev maps $T_{N, a}: X \rightarrow X$ are defined as $T_{N, a}(x) = \cos (N \arccos x + a)$, with the two parameters $N \in \mathbb{N}_{\geq 2}$ and $a \in \left [-\frac{\pi}{2}, 0\right]$. The piecewise-linear maps $g_{N, a}: Y \rightarrow Y$ are defined like a multi-tent or multi-upside-down tent map with piecewise constant slope $\pm N$ (cf. Fig.\ref{multi2}). 
Here we only consider $N = 2$ as an example; for a general $N \in \mathbb{N}_{\geq 2}$ the proof follows the same method (but one needs to distinguish odd and even $N$ cases). 
Denote 
\begin{equation*}
\begin{split}
f(x) &= T_{2, a}(x) = \cos (2 \arccos x + a), \quad x \in X = [-1, 1]\\
g(y) &= g_{2, a}(y) = \begin{cases}
2y + \alpha, \quad y \in \left[ 0, \frac{1 - \alpha}{2}\right) =: Y_1\\
-2y + 2 - \alpha, \quad y \in \left[ \frac{1 - \alpha}{2}, \frac{2 - \alpha}{2}\right) =: Y_2\\
2y + \alpha - 2, \quad y \in \left[ \frac{2 - \alpha}{2}, 1 \right] =: Y_3
\end{cases}
\end{split}
\end{equation*}
where $\alpha \in [0, 1]$ depends on the shift parameter $a$ and is to be determined. 
\\
\textit{Claim:} the function $h: X \rightarrow Y$, $h(x) = \frac{1}{\pi}\arccos (-x)$ is the topological conjugation between $f$ and $g$, such that $h \circ f = g \circ h$ for all $x \in X$. 

\begin{proof}
Consider the inverse $h^{-1}(y) = -\cos (\pi y)$, then we need to show that $h^{-1}(g(y)) = f(h^{-1}(y))$ for all $y \in Y$. 
First, 
\begin{equation*}
f(h^{-1}(y)) = \cos (2\arccos (-\cos (\pi y)) + a) = \cos (2(\pi y \pm \pi) + a) = \cos (2\pi y + a). 
\end{equation*}

Then, 
\begin{equation*}
\begin{cases}
y \in Y_1, \quad h^{-1}(g(y)) = -\cos(\pi (2y + \alpha)) = \cos(2\pi y + \pi \alpha \pm \pi); \\
y \in Y_2, \quad h^{-1}(g(y)) = -\cos(\pi (-2y + 2 - \alpha)) = -\cos (-(2\pi y + \pi \alpha)) = \cos(2\pi y + \pi \alpha \pm \pi); \\
y \in Y_3, \quad h^{-1}(g(y)) = -\cos(\pi (2y + \alpha - 2)) = \-\cos (2\pi y + \pi \alpha) = \cos(2\pi y + \pi \alpha \pm \pi).
\end{cases} 
\end{equation*}

Equating $f(h^{-1}(y))$ and $h^{-1}(g(y))$ gives $a = \pi \alpha \pm \pi$, or $\alpha = \frac{a}{\pi} \pm 1$. Since $\alpha \in [0, 1]$ we choose the positive sign and conclude
\begin{equation*}
\alpha = \frac{a}{\pi} + 1. 
\end{equation*} 

\end{proof}

\section{An example of a non-trivial invariant density for a piecewise-linear map $g_{N, a}$: $N = 3, a = -\frac{\pi}{9}$}\label{Amarkov}
\begin{figure}[H]
\centering 
\subfloat{
\includegraphics[width = 0.5\textwidth]{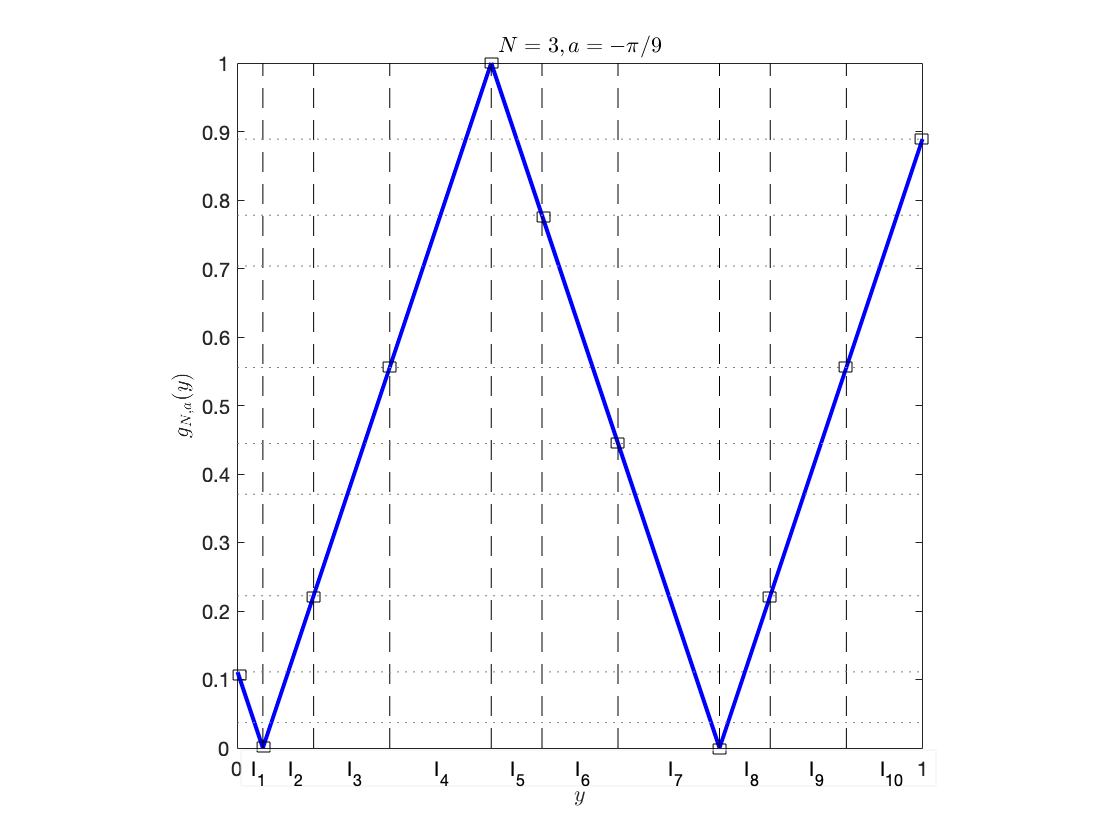}
}
\caption{\scriptsize Markov partition of the piecewise-linear map $g_{3, -\pi/9}(y)$ with partition points $y = 0, \frac{1}{27}, \frac{1}{9}, \frac{2}{9}, \frac{10}{27}, \frac{4}{9}, \frac{5}{9}, \frac{19}{27}, \frac{7}{9}, \frac{8}{9}, 1$, indicated by dashed vertical lines.} 
\label{N3a-pi9-markov}
\end{figure}

A Markov partition is defined such that the slope is constant on each subinterval and boundary points map to boundary points. For $g_{3, -\frac{\pi}{9}}$ such a partition, given by the eleven edge points of the partition indicated in Fig.\ref{N3a-pi9-markov}, consists of ten subintervals $I = [0, 1] = \cup_{i = 1}^{10} I_i$.

We have $g(I_1) = I_1 \cup I_2$, $g(I_2) = I_1 \cup I_2 \cup I_3$, $g(I_3) = I_4 \cup I_5 \cup I_6$, etc. 
The transition matrix $A$ is therefore 
\begin{equation*}
A = \begin{pmatrix}
1& 1& 0& 0& 0& 0& 0& 0& 0& 0\\
1& 1& 1& 0& 0& 0& 0& 0& 0& 0\\
0& 0& 0& 1& 1& 1& 0& 0& 0& 0\\
0& 0& 0& 0& 0& 0& 1& 1& 1& 1\\
0& 0& 0& 0& 0& 0& 0& 0& 1& 1\\
0& 0& 0& 0& 0& 1& 1& 1& 0& 0\\
1& 1& 1& 1& 1& 0& 0& 0& 0& 0\\
1& 1& 1& 0& 0& 0& 0& 0& 0& 0\\
0& 0& 0& 1& 1& 1& 0& 0& 0& 0\\
0& 0& 0& 0& 0& 0& 1& 1& 1& 0
\end{pmatrix}. 
\end{equation*}

The invariant density $\rho$ is determined by the (right-)eigenvector of the transfer matrix $B = \frac{1}{3}A^t$ ($t$ denotes transpose) associated with the unit eigenvalue: $\rho = B\rho$, 
\begin{equation*}
\rho = \begin{cases}
\alpha_1 \\
\alpha_2 = \alpha_1 \\
\alpha_3 = \frac{2}{3}\alpha_1 \\
\alpha_4 = \frac{21}{38}\alpha_1 \\
\alpha_5 = \frac{21}{38}\alpha_1 \\
\alpha_6 = \frac{33}{57}\alpha_1 \\
\alpha_7 = \frac{1}{2}\alpha_1 \\
\alpha_8 = \frac{1}{2}\alpha_1 \\
\alpha_9 = \frac{28}{57}\alpha_1 \\
\alpha_{10} = \frac{7}{19}\alpha_1 
\end{cases}
\end{equation*}

Using the normalisation condition $1 = \sum_{i = 1}^{10} \alpha_i$ we get $\alpha_1 = \frac{19}{118}$. Hence, the invariant density of $g_{3, -\frac{\pi}{9}}$ is
\begin{equation} 
\rho(y) = \begin{cases} 
\frac{19}{118} \approx 0.161, & y \in \left[0, \frac{1}{9}\right)\\
\frac{19}{177} \approx 0.107, & y \in \left[ \frac{1}{9}, \frac{2}{9}\right)\\
\frac{21}{236} \approx 0.089, & y \in \left[ \frac{2}{9}, \frac{4}{9}\right)\\
\frac{11}{118} \approx 0.093, & y \in \left[ \frac{4}{9}, \frac{5}{9}\right)\\
\frac{19}{236}\approx 0.081, & y \in \left[ \frac{5}{9}, \frac{7}{9}\right)\\
\frac{14}{177} \approx 0.079, & y \in \left[ \frac{7}{9}, \frac{8}{9}\right)\\
\frac{7}{118} \approx 0.059, & y \in \left[ \frac{8}{9}, 1\right]
\end{cases}
\label{invden-N3m9}%
\end{equation}
which coincides with the numerical result in Fig.\ref{pwlinear-N3-a-pi9-invden}.

\section{Proof of eigenfunctions of the PF operator for multi-upside-down tent and multi-tent maps}\label{A2} 
(i) For multi-upside-down tent maps, $g_N$ with even $N$, the general piecewise-linear map $g_N$ illustrated in Fig.\ref{multi-even} is defined explicitly as 
\begin{equation*}
g_N(x) = \begin{cases}
1 - Nx, \quad x \in [0, \frac{1}{N}]\\
Nx - 1, \quad x \in [\frac{1}{N}, \frac{2}{N}]\\
3 - Nx, \quad x \in [\frac{2}{N}, \frac{3}{N}]\\
\vdots \\
2k + 1 - Nx, \quad x \in [\frac{2k}{N}, \frac{2k + 1}{N}]\\
Nx - (2k + 1), \quad x \in [\frac{2k + 1}{N}, \frac{2k + 2}{N}]\\
\vdots \\
Nx - (N - 1), \quad x \in [\frac{N - 1}{N}, 1]
\end{cases}.
\end{equation*}

From definition of the PF operator $\mathcal{L}$ we have
\begin{equation*}
\mathcal{L}\rho(y) = \sum_{x \in g_N^{-1}(y)}\frac{\rho(x)}{|g'_N (x)|} = \frac{1}{N}(\rho(x_1) + \rho(x_2) + ... + \rho(x_N)),
\end{equation*}
where $x_i$ ($i = 1, ..., N$) are the preimages of $y \in [0, 1]$ under the map $g_N$. Assuming the eigenvalues are $\lambda^{(n)} = N^{-2n}$, then the eigenvalue equation $\mathcal{L}\rho^{(n)}(y) = \lambda^{(n)} \rho^{(n)}(y)$ gives 
\begin{equation*}
\frac{1}{N}(\rho^{(n)}(x_1) + \rho^{(n)}(x_2) + ... + \rho^{(n)}(x_N)) = N^{-2n}\rho^{(n)}(y),
\end{equation*}
or explicitly 
\begin{equation}
N^{2n - 1}\left[\rho^{(n)}\left( \frac{1-y}{N}\right) + \rho^{(n)}\left( \frac{1+y}{N}\right) + ... + \rho^{(n)}\left(\frac{y + (N - 1)}{N}\right)\right] = \rho^{(n)}(y),
\label{eigenvalue-fun-Neven}%
\end{equation}
for $n \in \mathbb{N}_0$. 
\\
\textit{Claim 1:} $\rho^{(n)}(x) = B_{2n}(\frac{x}{2} + \frac{1}{2})$, where $B_n (x)$ are the Bernoulli polynomials. 

\begin{proof}
Using the Multiplication Theorem for $B_n$ with $n = 2n, m = N$ (even) and $x = \frac{1}{N}\left( \frac{x}{2} + \frac{1}{2}\right)$ we have
\begin{equation}
\begin{split}
B_{2n}\left( \frac{x}{2} + \frac{1}{2}\right) &= B_{2n} \left( N\cdot \frac{x+1}{2N}\right)\\
&= N^{2n-1}\sum_{k = 0}^{N-1}B_{2n} \left( \frac{x+1}{2N} + \frac{k}{N}\right)\\
&= N^{2n-1}\left[ B_{2n}\left( \frac{x}{2N} + \frac{1}{2N}\right) + B_{2n}\left( \frac{x}{2N} + \frac{3}{2N}\right) + ... + B_{2n}\left( \frac{x}{2N} + \frac{2N-1}{2N}\right)\right].
\label{multi-thm-Neven}%
\end{split}
\end{equation} 

Now we want to show that each term in the last expression corresponds to a term in LHS of \eqref{eigenvalue-fun-Neven}. 

If \textit{Claim 1} is true, we have 
\begin{equation*}
\rho^{(n)}(x_1) = \rho^{(n)}\left(\frac{1-y}{N}\right) = B_{2n}\left( \frac{1}{2}\cdot \frac{1-y}{N} + \frac{1}{2}\right) = B_{2n} \left( \frac{1}{2} - \frac{1-y}{2N}\right) = B_{2n} \left( \frac{y}{2N} + \frac{N-1}{2N}\right),
\end{equation*}
where in the second equality we have used a symmetry property of $B_{2n}$ (basically whenever there is a minus sign in front of a term involving $y$ we need to use this symmetry), and also 
\begin{equation*}
\rho^{(n)}(x_2) = \rho^{(n)}\left( \frac{1+y}{N}\right) = B_{2n}\left( \frac{1}{2}\cdot \frac{1+y}{N} + \frac{1}{2}\right) = B_{2n} \left( \frac{y}{2N} + \frac{N+1}{2N}\right).
\end{equation*}

By induction, we have
\begin{equation*}
\begin{split}
\rho^{(n)}(x_{2k-1}) &= \rho^{(n)}\left( \frac{(2k-1)-y}{N}\right) = B_{2n}\left( \frac{y}{2N} + \frac{N + (2k-1)}{2N}\right), \\
\rho^{(n)}(x_{2k}) &= \rho^{(n)}\left( \frac{2k-y}{N}\right) = B_{2n}\left( \frac{y}{2N} + \frac{N - (2k-1)}{2N}\right),
\end{split}
\end{equation*}
for all integers $1 \leq k \leq \frac{N}{2}$. So the last two preimages $x_{N-1}$ and $x_N$ correspond to
\begin{equation*}
\begin{split}
\rho^{(n)}(x_{N-1}) &= \rho^{(n)}\left( \frac{(N-1)-y}{N}\right) = B_{2n}\left( \frac{y}{2N} + \frac{2N-1}{2N}\right), \\
\rho^{(n)}(x_{N}) &= \rho^{(n)}\left( \frac{N-y}{N}\right) = B_{2n}\left( \frac{y}{2N} + \frac{1}{2N}\right),
\end{split}
\end{equation*}
for which reason we group pairs of consecutive preimages together in Fig.\ref{multi-even}. Therefore, \eqref{eigenvalue-fun-Neven} and \eqref{multi-thm-Neven} are equivalent, and \textit{Claim 1} together with the eigenvalue assumption hold. 

\end{proof}

(ii) For multi-tent maps (cf. Fig.\ref{multi-odd}), similarly, the piecewise-linear map is defined as 
\begin{equation*}
g_N(x) = \begin{cases}
Nx, \quad x \in [0, \frac{1}{N}]\\
2 - Nx, \quad x \in [\frac{1}{N}, \frac{2}{N}]\\
Nx - 2, \quad x \in [\frac{2}{N}, \frac{3}{N}]\\
\vdots \\
Nx - 2k, \quad x \in [\frac{2k}{N}, \frac{2k + 1}{N}]\\
2k + 2 - Nx, \quad x \in [\frac{2k + 1}{N}, \frac{2k + 2}{N}]\\
\vdots \\
Nx - (N - 1), \quad x \in [\frac{N - 1}{N}, 1]
\end{cases}. 
\end{equation*}

Again, by definition of the transfer operator $\mathcal{L}$ and eigenvalue equation, assuming $\lambda^{(n)} = N^{-2n}$, we get 
\begin{equation}
N^{2n-1}\left[ \rho^{(n)}\left( \frac{y}{N}\right) + \rho^{(n)}\left( \frac{2-y}{N}\right) + ... + \rho^{(n)}\left( \frac{y}{N} + \frac{N-1}{N}\right) \right] = \rho^{(n)}(y),
\label{eigenvalue-fun-Nodd}%
\end{equation}
for $n \in \mathbb{N}_0$. 
\\
\textit{Claim 2a:} $\rho^{(n)}(x) = B_{2n}(x)$, where $B_n (x)$ are the Bernoulli polynomials. 

\begin{proof}
By the Multiplication Theorem for $B_n$ with $n = 2n, m = N$ (odd) and $x = \frac{x}{N}$ we have
\begin{equation*}
\begin{split} 
B_{2n}(x) &= B_{2n}\left( N\cdot \frac{x}{N}\right) \\
&= N^{2n-1}\sum_{k = 0}^{N-1}B_{2n}\left( \frac{x}{N} + \frac{k}{N}\right) \\
&= N^{2n-1}\left[ B_{2n}\left( \frac{x}{N}\right) + B_{2n}\left( \frac{x}{N} + \frac{1}{N}\right) + ... + B_{2n} \left( \frac{x}{N} + \frac{N-1}{N}\right) \right].
\end{split}
\end{equation*}

This expression is equivalent to \eqref{eigenvalue-fun-Nodd} since there is a term-to-term correspondence \footnote{\scriptsize This is a re-ordering the $N$ (odd) preimages $\{ x_1, x_2, ..., x_N\}$ in the way that 
\\
i) if $N = 4q +1$, the index of the middle point of natural ordering is odd then 
\begin{equation*}
\begin{split}
x_1 &\rightarrow x_{N-1} \rightarrow x_3 \rightarrow x_{N-3} \rightarrow x_5 \rightarrow x_{N-5} \rightarrow ... \\
&\rightarrow x_{\frac{N+1}{2} + 1} \rightarrow x_{\frac{N+1}{2}} \rightarrow x_{\frac{N+1}{2} - 1} \rightarrow x_{\frac{N+1}{2} + 2} \rightarrow x_{\frac{N+1}{2} - 3} \rightarrow x_{\frac{N+1}{2} + 4} \rightarrow x_{\frac{N+1}{2} - 5} \rightarrow ... \rightarrow x_2 \rightarrow x_N; 
\end{split}
\end{equation*}
and ii) if $N = 4q - 1$, the index of the middle point of natural ordering is even then 
\begin{equation*}
\begin{split}
x_1 &\rightarrow x_{N-1} \rightarrow x_3 \rightarrow x_{N-3} \rightarrow x_5 \rightarrow x_{N-5} \rightarrow ... \\
&\rightarrow x_{\frac{N+1}{2} - 1} \rightarrow x_{\frac{N+1}{2}} \rightarrow x_{\frac{N+1}{2} + 1} \rightarrow x_{\frac{N+1}{2} - 2} \rightarrow x_{\frac{N+1}{2} + 3} \rightarrow x_{\frac{N+1}{2} - 4} \rightarrow x_{\frac{N+1}{2} + 5} \rightarrow ... \rightarrow x_2 \rightarrow x_N.
\end{split}
\end{equation*}
\\
Compare to the even $N$ case, where the re-ordering reads 
\begin{equation*}
x_N \rightarrow x_{N-2} \rightarrow x_{N-4} \rightarrow ... \rightarrow x_4 \rightarrow x_2 \rightarrow x_1 \rightarrow x_3 \rightarrow x_5 \rightarrow ... \rightarrow x_{N-3} \rightarrow x_{N-1}.
\end{equation*}
}

\begin{equation*}
\begin{split}
\rho^{(n)}(x_{2k-1}) &= \rho^{(n)}\left( \frac{y}{N} + \frac{2k-2}{N}\right) = B_{2n}\left( \frac{y}{N} + \frac{2k-2}{N}\right), \\
\rho^{(n)}(x_{2k}) &= \rho^{(n)}\left( \frac{2k}{N} - \frac{y}{N}\right) = B_{2n}\left( \frac{2k}{N} - \frac{y}{N}\right) = B_{2n}\left( \frac{1}{2} - (\frac{y}{N} + \frac{1}{2} - \frac{2k}{N})\right) \\
&\quad = B_{2n}\left( \frac{1}{2} + (\frac{y}{N} + \frac{1}{2} - \frac{2k}{N})\right) = B_{2n} \left( \frac{y}{N} + \frac{N-2k}{N}\right), 
\end{split}
\end{equation*}
for all integers $1 \leq k \leq \frac{N-1}{2}$, for which one can easily verify that 
\begin{equation*}
\begin{split}
\rho^{(n)}(x_1) &= B_{2n}\left( \frac{y}{N}\right), \\
\rho^{(n)}(x_2) &= B_{2n} \left( \frac{y}{N} + \frac{N-2}{N}\right), 
\end{split}
\end{equation*}
and 
\begin{equation*}
\begin{split}
\rho^{(n)}(x_{N-1}) &= B_{2n}\left( \frac{y}{N}+\frac{1}{N}\right), \\
\rho^{(n)}(x_N) &= B_{2n}\left( \frac{y}{N}+\frac{N-1}{N}\right). 
\end{split}
\end{equation*} 
Hence the \textit{Claim 2a}. 

\end{proof}

\textit{Claim 2b:} $\rho^{(n)}(x) = E_{2n - 1}(x)$, where $E_n(x)$ are the Euler polynomials. 

\begin{proof}
Indeed it can be verified that $E_{2n - 1} (x)$ is also an eigenfunction for the multi-tent map $g_N(x)$ with $N \geq 3$ odd, corresponding to the eigenvalue $N^{-2n}$. The proof uses the Multiplication Theorem for $E_n$
\begin{equation*}
E_n(mx) = m^n \sum_{k = 0}^{m-1} (-1)^k E_n\left( x + \frac{k}{m}\right), \quad m = 1, 3, ..., 
\end{equation*}
with $k = N$ (odd), $x = \frac{x}{N}$ and $m = 2n-1$, and an anti-symmetry property for odd Euler polynomials: $E_{2n-1}\left(\frac{1}{2} - x\right) = -E_{2n-1}\left(\frac{1}{2} + x\right)$. This indicates a degeneracy for odd $N$ multi-tent maps $g_N (x)$, therefore also for odd Chebyshev maps $T_{N, 0}(x)$. 
\end{proof}

\section{Discussion of eigenfunctions for odd-$N$ shifted Chebyshev maps}\label{AshiftedEfun}
In this case, a more complicated situation arises from the fact that, when $a \neq 0$, for an arbitrary $x \in [0, 1]$ the piecewise-linear map $g_{N, a}$ can have different numbers of preimages, for example, for $g_{3, a}$ (cf. Fig.\ref{counter-eg-of-full-branch}), there can be 1, 2, 3 or 4 preimages that separate the unit interval into several subintervals, and we need to take into account the values of $x$ in each of the subintervals, resulting in the eigenfunctions being piecewise defined. This is consistent with the fact from section 2.3 that the invariant densities for some $T_{N, a}$ are not always smooth but piecewise smooth (cf. last three columns in Fig.\ref{3egs-inv-den}). 

We therefore conjecture that in the case of $a = -\frac{\pi}{m}$ ($m = 2, 3, ...$), eigenfunctions of the PF operator for $T_{N, a}$ are of the form 
\begin{equation}
\rho^{(n)}_{T_{N, a}}(x) = \sum_{k = 1}^{mN^n} \alpha_{n,k} \chi_{I_{n, k}}(x) \frac{1}{\pi \sqrt{1 - x^2}}F_n(x), \quad x \in [-1, 1], 
\label{eigenfun-conjec2}%
\end{equation}
where $\alpha_{n, k} = \alpha_{n, k}(\lambda^{(n)})$ are the weights depending on the associated eigenvalue $0<\lambda^{(n)} \leq 1$, $\chi_I (x)$ is the indicator function such that $\chi_I(x) = 1$ if $x \in I$ and $0$ otherwise. The partition $I_{n, k} = \left[ \frac{k-1}{mN^n}, \frac{k}{mN^n} \right)$ may vary according to the order $n$ of the eigenvalue; the last subinterval includes the right-boundary point 1, $I_{n, mN^n} = \left[ \frac{mN^n - 1}{mN^n}, 1\right]$, so that $I_{n, k_1} \cap I_{n, k_2} = \emptyset, \forall k_1 \neq k_2$ and $\cup_k I_{n, k} = I = [0, 1]$ for any given $n \in \mathbb{N}_0$. This forms a hierarchically finer structure in partitioning the unit interval, initially defined by the Markov partition at the invariant density level ($n = 0$, with $\lambda^{(0)} = 1$ being the largest eigenvalue). $\{F_n\}$ are some appropriate \textit{smooth} functions with $F_0(x) = 1$ $\forall x \in [-1, 1]$. 

In particular, similarly to eq.\eqref{semi-conj1}, we have 
\begin{equation}
(-T_{2m, 0})\circ T_{N, -\frac{\pi}{m}} = T_{N, 0}\circ (-T_{2m, 0}) = -T_{2mN, 0}, \quad N \text{ odd},
\label{semi-conj2}%
\end{equation}
that is, the odd $N$ shifted Chebyshev map $T_{N, -\frac{\pi}{m}}$ is topologically semi-conjugated to its corresponding ordinary $T_{N, 0}$ via the semi-conjugacy $h_2(x) := -T_{2m ,0}(x)$, $m \in \mathbb{N}_{\geq 2}$. By a coordinate transformation one gets the shape of eigenfunctions for $T_{N, -\frac{\pi}{m}}$ ($N$ odd) as (cf. eq.\eqref{efun-TN})
\begin{equation*}
\begin{split}
|h_2'(x)|\rho_{T_{N, 0}}^{(n)}(h_2(x)) \propto & \frac{1}{\pi \sqrt{1 - x^2}}B_{2n}\left( \frac{2m}{\pi}\arccos x\right) \\
\text{or } \propto & \frac{1}{\pi \sqrt{1 - x^2}}E_{2n-1}\left( \frac{2m}{\pi}\arccos x\right).
\end{split}
\end{equation*}

So the set of smooth functions $\{ F_n \}$ in \eqref{eigenfun-conjec2} can be chosen to be Bernoulli or Euler polynomials, and the eigenvalues $\lambda^{(n)} = N^{-2n}$ remain the same.

\section{Derivation of higher-order correlations of iterates of shifted Chebyshev maps}\label{A3} 
Let $x_i$ be the $i$-th iterate of a shifted Chebyshev map $T_{N, 0}$. Using the change of variables described in Sec.2, 
\begin{equation*}
\begin{split}
x_0 &= \cos(\pi u),\\
x_n &= \cos \left( N^n\pi u + \frac{N^n - 1}{N - 1}a\right),
\end{split}
\end{equation*}
we have, for shifted Chebyshev maps $T$ in the category \eqref{star-cond},
\begin{equation*}
\begin{split}
\langle x_{n_1}\cdots x_{n_r} \rangle =& \int_{-1}^1 \frac{1}{\pi \sqrt{1 - \cos^2 \pi u}} \cos \left(N^{n_1}\pi u + \frac{N^{n_1}-1}{N-1}a\right)\cdot...\cdot \cos \left(N^{n_r}\pi u + \frac{N^{n_r}-1}{N-1}a\right) d(\cos \pi u)\\
=& \int_{1}^0 \frac{1}{\pi \sin \pi u} \cos \left(N^{n_1}\pi u + \frac{N^{n_1}-1}{N-1}a\right)\cdot...\cdot \cos \left(N^{n_r}\pi u + \frac{N^{n_r}-1}{N-1}a\right) \cdot (-\pi \sin \pi u) du\\
=& \int_0^1 \cos \left(N^{n_1}\pi u + \frac{N^{n_1}-1}{N-1}a\right)\cdot...\cdot \cos \left(N^{n_r}\pi u + \frac{N^{n_r}-1}{N-1}a\right) du\\
=& \int_0^1 \frac{\exp\left(i(N^{n_1}\pi u + \frac{N^{n_1} -1}{N-1}a)\right) + \exp \left(-i(N^{n_1}\pi u + \frac{N^{n_1} - 1}{N-1}a)\right)}{2}\cdot ... \\
& \quad \cdot \frac{\exp \left(i(N^{n_r}\pi u + \frac{N^{n_r} -1}{N-1}a)\right) + \exp \left(-i(N^{n_r}\pi u + \frac{N^{n_r} - 1}{N-1}a)\right)}{2} du\\
=& 2^{-r} \sum_{\sigma}\int_0^1 \exp \left(i \sum_{l=1}^r \sigma_l \cdot \left(N^{n_l}\pi u + \frac{N^{n_l} - 1}{N-1}a\right)\right)du\\
=& 2^{-r} \sum_{\sigma} \left[ \exp\left(ia \sum_{l=1}^r \sigma_l \frac{N^{n_l} - 1}{N-1}\right) \cdot \int_0^1 \exp \left(i \pi u\sum_{l=1}^r \sigma_l N^{n_l}\right) du\right] \\
=& 2^{-r} \sum_{\sigma} \left[ \exp\left(ia \sum_{l=1}^r \sigma_l \frac{N^{n_l} - 1}{N-1}\right) \cdot \delta \left( \sum_{l = 1}^r \sigma_l N^{n_l}, 0\right) \right], 
\end{split}
\end{equation*}
where in the fourth line we have used the Euler formula $\cos x = \frac{1}{2}(e^{ix} + e^{-ix})$; in the fifth line $\sigma_l$ denotes a choice of the spin configuration $\{ -1, +1\}$ and the sum over $\sigma$ is a summation over all possible configurations $\sigma := (\sigma_1, \sigma_2, ..., \sigma_r)$, $\sigma_l \in \{ -1, +1\}$; in the last line $\delta(x, 0)$ is the Kronecker delta defined as being 1 if $x=0$ and 0 else.

\section{Derivation of the two-point correlation function for the $N$-ary shift with subtracted mean}\label{A4}
Consider the $N$-ary shift $x_{n+1} = f(x_n) = Nx \mod 1$, $x_0 \in [0, 1]$, $N \in \mathbb{N}_{\geq 2}$. 

With the $N$-ary representation of the initial point $x_0 = \sum_{j = 1}^{\infty}a_j N^{-j}$, $a_j \in \{0, 1, ..., N-1 \}$ we have $x_n = \sum_{j = 1}^{\infty}a_{j+n}N^{-j}$. 

Define 
\begin{equation*}
\begin{split}
w_n :&= x_n - \langle x_n\rangle \\
&= \sum_{j = 1}^{\infty} a_{j+n}N^{-j} - \frac{1}{2}\\
&= \sum_{j = 1}^{\infty} a_{j+n}N^{-j} - \frac{N-1}{2}\sum_{j = 1}^{\infty}N^{-j}\\
&= \sum_{j = 1}^{\infty} b_{j+n}N^{-j}, 
\end{split}
\end{equation*}
where $b_{j+n} \in \left\{ -\frac{N-1}{2}, -\frac{N-3}{2}, ..., \frac{N-1}{2}\right\}$ so that $w_n$ has the average value zero. 

Then 
\begin{equation*}
\begin{split}
\langle b_{j+n}^k\rangle &= \frac{1}{N}\left[ \left( -\frac{N - 1}{2}\right)^k + \left( -\frac{N - 3}{2}\right)^k + ... + \left( \frac{N - 3}{2}\right)^k + \left( \frac{N - 1}{2}\right)^k\right]\\
&=\begin{cases}
0, \quad \text{if $k$ is odd}\\
\frac{2}{N}\left[ \left( \frac{N - 1}{2}\right)^k + \left( \frac{N - 3}{2}\right)^k + ... + \left( \frac{N - \ceil{\frac{N+1}{2}}}{2}\right)^k\right], \quad \text{if $k$ is even}.
\end{cases}
\end{split}
\end{equation*}

Let us denote 
\begin{equation*}
\begin{split}
S_{N, k} :&= \left( \frac{N - 1}{2}\right)^k + \left( \frac{N - 3}{2}\right)^k + ... + \left( \frac{N - \ceil{\frac{N+1}{2}}}{2}\right)^k\\
&=\begin{cases}
\frac{1}{2^k}(1 + 3^k + ... + (N-1)^k) \quad \text{if $N$ is even}\\
1 + 2^k + ... + \left( \frac{N-1}{2}\right)^k \quad \text{if $N$ is odd}.
\end{cases}
\end{split}
\end{equation*} 
Note: one can get an explicit formula for the sum of the $k$-th powers of the first $n$ natural numbers by Faulhaber's formula \cite{faulhaber}. We shall use it when computing some specific examples. 

Consider $b_{j_1}\cdot ... \cdot b_{j_r} = b_{m_1}^{k_1}\cdot ... \cdot b_{m_s}^{k_s}$ ($s \leq r$), then 
\begin{equation*}
\langle b_{j_1}\cdot ...\cdot b_{j_r}\rangle = \begin{cases}
0, \quad \text{if $\exists i \in \{ 1, ..., s\}$ s.t. $k_i$ is odd}\\
\frac{2^s}{N^s}S_{N, k_1}\cdot ... \cdot S_{N, k_s}, \quad \text{otherwise}.
\end{cases}
\end{equation*}

The $r$th-order correlation of iterates of $w_n$ is 
\begin{equation*}
\begin{split}
\langle w_{n_1}\cdot ... \cdot w_{n_r}\rangle &= \langle \left( \sum_{j_1 = 1}^{\infty}b_{j_1 + n_1}N^{-j_1}\right) \cdot ... \cdot \left( \sum_{j_r = 1}^{\infty}b_{j_r + n_r}N^{-j_r}\right)\rangle \\
&= \sum_{j_1 = 1}^{\infty}...\sum_{j_r = 1}^{\infty}\langle b_{j_1 + n_1}\cdot ... \cdot b_{j_r + n_r}\rangle N^{-(j_1+ ... + j_r)}\\
&\begin{cases}
= 0, \quad \text{if $r$ is odd}\\
\neq 0, \quad \text{if $r$ is even}.
\end{cases}
\end{split}
\end{equation*}

In particular, for $r = 2$ we have
\begin{equation}
\begin{split}
\langle w_{n_1}w_{n_2}\rangle &= \sum_{j_1 = 1}^{\infty}\sum_{j_2 = 1}^{\infty} \langle b_{j_1 + n_1} b_{j_2 + n_2}\rangle N^{-(j_1 + j_2)}\\
&= \sum_{j_1 = 1}^{\infty}\sum_{j_2 = 1}^{\infty}\delta(j_1 + n_1, j_2 + n_2) \cdot \frac{2^2}{N^2}S_{N, 2}\cdot \frac{1}{N^{j_1 + j_2}}\\
&= \frac{4}{N^2}S_{N, 2}\sum_{j_1 = 1}^{\infty}\frac{1}{2^{2j_1 + |n_1 - n_2|}}\\
&= \frac{4}{N^2}\cdot \left(\frac{1}{24}N(N-1)(N+1)\right) \cdot \frac{1}{N^2 - 1}\cdot \frac{1}{N^{|n_1 - n_2|}}\\
&= \frac{1}{6N}\cdot \frac{1}{N^{|n_1 - n_2|}},
\end{split}
\label{two-pt-corr-Nary}%
\end{equation}
where one can verify that $S_{N, 2}$ equals to $\frac{1}{24}N(N-1)(N+1)$ for both even and odd $N$'s: 

i) if $N$ is even (e.g., $N = 2$ the binary shift), we have 
\begin{equation*}
\begin{split}
S_{N, 2} &= \left(\frac{N-1}{2}\right)^2 + \left(\frac{N-3}{2}\right)^2 + ... + \left(\frac{1}{2}\right)^2\\
&= \left(\frac{1}{2}\right)^2 + \left(\frac{3}{2}\right)^2 + ... + \left(\frac{N-1}{2}\right)^2\\
&= \frac{1}{4}( 1^2 + 3^2 + ... + (N-1)^2)\\
&= \frac{1}{4}\cdot \frac{1}{3}\cdot \frac{N}{2}(N-1)(N+1)\\
&= \frac{1}{24}N(N-1)(N+1),
\end{split}
\end{equation*}
where we have used $1 + 3^2 + 5^2 + ... + (2k-1)^2 = \frac{1}{3}k(2k-1)(2k+1)$. 

ii) if $N$ is odd, we have
\begin{equation*}
\begin{split}
S_{N, 2} &= \left(\frac{N-1}{2}\right)^2 + \left(\frac{N-3}{2}\right)^2 + ... + 1^2 + 0^2\\
&= 1^2 + 2^2 + ... + \left(\frac{N-1}{2}\right)^2\\
&= \frac{1}{6}\cdot \frac{N-1}{2}\left(\frac{N-1}{2} + 1\right)\left( 2\cdot \frac{N-1}{2} + 1 \right)\\
&= \frac{1}{6}\cdot \frac{N-1}{2}\cdot \frac{N+1}{2}\cdot N\\
&= \frac{1}{24}N(N-1)(N+1), 
\end{split}
\end{equation*}
where we have used $1 + 2^2 + ... + k^2 = \frac{1}{6}k(k+1)(2k+1)$.

\end{document}